\documentclass[aop,MSNbibl,seceqn,dvips]{arximspdf}
\usepackage{graphicx}

\doi{10.1214/12-AOP801} %kopijuoti is PTS
\volume{41}
\issue{6}
\pubyear{2013}
\firstpage{3697}
\lastpage{3785}

\makeatletter

\newcommand{\inn}{\mathrm{in}}

\newcommand{\rright}{\right}
\newcommand{\lleft}{\left}

\newcommand{\cal}{\mathcal}

\newtheorem{theorem}{Theorem}[section]
\newtheorem{lemma}[theorem]{Lemma}
\newtheorem{proposition}[theorem]{Proposition}
\newproclaim{definition}[theorem]{Definition}
\newtheorem{corollary}[theorem]{Corollary}

\def\bpi{\overline{\pi}}
\def\btau{\tau}

\def\Z{{\mathbb Z}}
\def\R{{\mathbb R}}
\def\N{{\mathbb N}}
\def\P{P}

\def\cE{{\cal E}}
\def\mut{\widetilde{\mu}}
\def\cA{{\cal A}}
\def\tb{{\tau_\beta}}
\def\tic{{\widetilde{c}}}
\def\hic{{\widehat{c}}}
\def\tbot{\tau_{\mathrm{bottom}}}
\def\ocA{\overline{\cal A}}
\def\ocM{\overline{\cal M}}
\def\ot{\widetilde{\omega}}
\def\cR{{\cal R}}
\def\cC{{\cal C}}
\def\cG{{\cal G}}
\def\cD{{\cal D}}
\def\cX{{\cal X}}
\def\cB{{\cal B}}
\def\cM{{\cal M}}
\def\cV{{\cal V}}

\def\STC{\operatorname{\mbox{\textsc{stc}}}}
\def\SES{\mbox{\textsc{ses}}}
\def\Ga{\Gamma}
\def\st{\widetilde\sigma}
\def\mt{\widetilde\mu}
\def\mh{\widehat\mu}
\def\G{{\Gamma}}
\def\L{{\Lambda}}
\def\k{{\kappa}}
\def\a{{\alpha}}
\def\b{{\beta}}
\def\e{{\varepsilon}}
\def\h{{\eta}}
\def\s{{\sigma}}
\def\t{{\tau}}
\def\o{{\omega}}

\newcommand{\dinf}{\mathrm{d}_{\infty}}

\def\cA{{\cal A}}
\def\cB{{\cal B}}
\def\cC{{\cal C}}
\def\cD{{\cal D}}
\def\cE{{\cal E}}
\def\cG{{\cal G}}
\def\cM{{\cal M}}
\def\cR{{\cal R}}
\def\cT{{\cal T}}

\def\cV{{\cal V}}
\def\cX{{\cal X}}
\def\cY{{\cal Y}}

\newcommand{\diam}{\operatorname{diam}_\infty}
\def\tc{\dvt}
\def\plus{\mathbf{+1}}
\def\minus{\mathbf{-1}}
\def\per{\operatorname{perimeter}}

\def\R{{\mathbb R}}
\def\N{{\mathbb N}}

\newcommand{\height}{{\operatorname{height}}}
\newcommand{\depth}{\operatorname{depth}}
\newcommand{\bottom}{\operatorname{bottom}}
\def\uv{\underline{v}}
\def\ov{\overline{v}}

\makeatother

\begin{document}
\begin{frontmatter}

\title{Nucleation and growth for the Ising model in \lowercase{$d$}~dimensions at very low temperatures}
\runtitle{Nucleation and growth for the Ising model}

\begin{aug}
\author[A]{\fnms{Rapha\"el} \snm{Cerf}\corref{}\ead[label=e1]{rcerf@math.u-psud.fr}}
\and
\author[B]{\fnms{Francesco} \snm{Manzo}\ead[label=e2]{manzo.fra@gmail.com}}
\runauthor{R. Cerf and F. Manzo}
\affiliation{Universit\'e Paris Sud et IUF and Universit\`a di Roma Tre}
\address[A]{Math\'ematique, B\^atiment 425\\
Universit\'e Paris Sud\\
91405 Orsay Cedex\\
France\\
\printead{e1}}
%adresu isvedimo komanda gale!
\address[B]{Dipartmento di Matematica\\
Universit\'a ``Roma Tre''\\
Largo S. G. Murialdo 1\\
00100 Roma\\
Italy\\
\printead{e2}}
\end{aug}

% HISTORY:
\received{\smonth{2} \syear{2011}}
\revised{\smonth{9} \syear{2012}}

% ABSTRACT
%
\begin{abstract}
This work extends to dimension $d\geq3$ the main result of Dehghanpour
and Schonmann. We consider the stochastic Ising model on $\Z^d$
evolving with the Metropolis dynamics under a fixed small positive
magnetic field $h$ starting from the minus phase.
%As we scale the inverse temperature $\beta$ to $\infty$,
When the inverse temperature $\beta$ goes to $\infty$,
the
relaxation time of the system, defined as the time when
the plus phase has invaded the origin, behaves like
$\exp({\b\k_d})$. The value $\kappa_d$ is equal to
\[
\k_d = \frac{1}{d+1} (\G_1+\cdots+\G_d),
\]
where $\G_i$ is the energy of the $i$-dimensional critical droplet of
the Ising model at zero temperature and magnetic field $h$.
\end{abstract}

% KEYWORDS
\begin{keyword}[class=AMS]
\kwd{60K35}
\kwd{82C20}
\end{keyword}
\begin{keyword}
\kwd{Ising}
\kwd{Metropolis}
\kwd{metastability}
\kwd{nucleation}
\kwd{growth}
\end{keyword}
\vspace*{6pt}
\end{frontmatter}

\tableofcontents[level=2]

%s1 #&#
\section{Introduction}
We consider the kinetic Ising model in $\Z^d$ under a small positive
magnetic field in the limit of vanishing temperature, and
we study the relaxation of the system starting from the
metastable state where all the spins are set to minus.
%We present the background of the metastability problem
%for the Ising model
%in $d$ dimensions at very low temperatures in infinite volume
An introduction of the metastability problem is presented
in Section~\ref{back}. In Section~\ref{prob}, we explain the three
major problems we had to solve to extend the two-dimensional results
to dimension $d$.
The main results are stated in Section~\ref{mr}.
The strategy of the proof
is explained in Section~\ref{sp}.

%s1.1 #&#
\subsection{Background}
\label{back}
This work extends to dimension $d\geq3$ the main result of
Dehghanpour and Schonmann \cite{DS1}.
We consider the stochastic Ising model on $\Z^d$ evolving with
the Metropolis dynamics under a fixed small positive magnetic field~$h$.
We start the system in the minus phase. Let $\tau_d$ be the typical
relaxation time of the system, defined here as the time where
the plus phase has invaded the origin.
We will study the asymptotic behavior of $\tau_d$ when we
scale the temperature to $0$.
The corresponding problem in finite volume
(i.e., in a box $\Lambda$ whose size is fixed) has been previously
studied in arbitrary dimension by Neves \mbox{\cite{N2,N1}}.
% and in dimensions $3$ in \cite{BAC}.
In this situation, Neves proved that the relaxation time behaves as
$\exp(\beta\Gamma_d)$
where $\beta=1/T$ is the inverse temperature, and
$\Gamma_d$ is the energy barrier the system has to overcome
to go from the metastable state $\minus$ to the stable state $\plus$.
An explicit formula is available for $\Ga_d$; however the formula is quite
complicated.
The energy barrier $\Gamma_d$ is the solution of a minimax problem,
and it
is reached for configurations which are optimal saddles
between $\minus$ and $\plus$ in the energy landscape of the Ising model.
These results have been refined in dimension $3$ in \cite{BAC}.
In dimension $3$,
the optimal saddles are identified.
They are configurations called
critical droplets which
contain exactly one connected component
of pluses of cardinality~$m_3$, and their shape is an appropriate
union of a specific quasicube (whose sides depend on~$h$) and a two-dimensional
critical droplet.
In dimension $d\geq4$, the results of Neves yield
that the configurations
consisting of the appropriate union of a $d$-dimensional quasicube and a
$(d-1)$-dimensional critical droplet are optimal saddles, but it
is currently
not proved that they are the only ones. However, it is reasonable
to expect that the cases of equality in the discrete isoperimetric
inequality on the lattice can be analyzed in dimension $d\geq4$
in the same way they were studied in dimension $d=3$ \cite{AC},
so that the three-dimensional results could be extended to higher dimensions.

In infinite volume, instead of nucleating locally in a finite box
near the origin,
a~critical droplet of pluses might be created far from the origin,
and this droplet can grow, become supercritical and invade the origin.
It turns out that this is the most efficient mechanism to relax
to equilibrium.
This was shown by Dehghanpour and Schonmann in the two-dimensional case
\cite{DS1} and it required several new ideas and insights
compared to the finite volume
analysis.
Indeed, one has to understand the typical
birth place of the first critical droplets which are likely to invade
the origin, as well as their growth mechanism.
The heuristics given in \cite{DS1} apply in $d$ dimensions as well.
Suppose that nucleation in a finite box is exponentially
distributed with rate $\exp({-\b\G_d})$,
independently from other boxes, and
that the speed of growth of a large supercritical droplet is $v_d$.
The droplets which can reach the origin
at time $t$
are the droplets which are born inside the space--time cone
whose basis is a $d$-dimensional square with side length $v_dt$
and whose height is $t$.
The critical space--time cone is
such that its volume times the nucleation rate is of order one.
Let $\t_d$ be the typical relaxation time in dimension $d$,
that is, the time when the stable plus phase invades the origin.
From the previous heuristics, we conclude
that $\t_d$ satisfies
\[
\tfrac{1}{3}\t_d ( v_d \t_d
)^d \exp({-\b\G_d}) = 1.
\]
Solving this identity
and neglecting the factor ${1}/{3}$, we get
\[
\t_d = \exp \biggl( \frac{1}{d+1} (\beta\G_d-d\ln
v_d ) \biggr).
\]
%
%A first major difficulty is to estimate the speed of growth $v_d$ of
%large supercritical droplets.
Since the large supercritical droplets are
approximately parallelepipeds, the dynamics on one face behaves
like a $(d-1)$-dimensional stochastic Ising model, and the time
needed to fill a face with pluses is of order $\tau_{d-1}$.
Thus $v_d$ should behave like the inverse of $\tau_{d-1}$, and the
previous formula becomes
\[
\ln\t_d = %\exp\Big(
\frac{1}{d+1} (\beta\G_d+d\ln
\t_{d-1} ).
\]
In this computation, we take only into account the terms on
the exponential scale, of order $\exp(\beta \mbox{ constant})$.
Setting $\t_d=\exp({\b\k_d})$, the constant $\k_d$ satisfies
%the following recursive formula:
%
\[
\k_d = \frac{1}{d+1} (\G_d+d\k_{d-1} ).
\]
Solving the recursion, and using that $\k_0=0$, we get that
\[
\k_d = \frac{1}{d+1} (\G_1+\cdots+\G_d
).
\]
%
%s1.2 #&#
\subsection{Three major problems}
\label{prob}
Although these heuristics are rather convincing, it is a real
challenge to prove rigorously that the asymptotics of
the relaxation time are indeed of order
$\exp({\b\k_d})$.
Our strategy is to implement inductively
%with respect to the dimension
the scheme of Dehghanpour and Schonmann.
To do so, we had to overcome three major problems.\vspace*{9pt}

\textit{Speed of growth.}
A first major difficulty is to control
the speed of growth $v_d$ of
large supercritical droplets.
The upper bound on the speed of growth in \cite{DS2} was based on
a very detailed analysis of the growth of an infinite interface.
Using a combinatorial argument based on chronological paths,
first introduced by Kesten and Schonmann in the context of
a simplified growth model \cite{KS}, Dehghanpour and Schonmann
were able to prove that
%In dimensions $2$, it turned out that
$v_2$ is of order $\exp(-\beta\Ga_1/2)$.
Despite considerable efforts,
we never managed to extend this technique of analysis to higher dimension.
Here
we consider only
interfaces with a size that is exponential in
$\b$.
In order to control the growth of these interfaces, we use inductively
coupling techniques
introduced
to analyze
the finite-size scaling in the bootstrap percolation model
\cite{CeCi,CM}.
We apply successively these techniques
in two distinct ways, the first sequential
and the second parallel. This strategy has been elaborated first
in a simplified growth model
\cite{CM2}, yet its application in the context of the Ising model
is more troublesome.
Contrary to the case of the growth model, we did not manage to compare
the dynamics in a strip with a genuine $(d-1)$-dimensional dynamics,
and we perform the induction on the boundary conditions rather than
on the dimension.
An additional source of trouble is to control the configurations
in the metastable
regions.
We introduce an adequate hypothesis describing their law, which is
preserved until the arrival of supercritical droplets,
in order to
tackle this problem.
A~key result to control the speed of growth is Theorem
\ref{T2}.\vspace*{9pt}

\textit{Energy landscape.}
A second major difficulty is that
it is very hard to analyze the
energy landscape of the Ising model
in high dimension,
and
the results we are able to obtain are very weak compared to
the corresponding results in finite volume and in
dimensions two and three;
see \cite{NS1,NS2,BAC}.
%The geometry of the high-dimensional clusters is
%very complicated and
For instance we are not able to
determine whether a given cluster of pluses tends to shrink
or to grow.
Moreover, we do not know
some of the fine details of the energy landscape such as the
depth of the small cycles that could trap the process
and increase the relaxation time.
In other words,
we do not know how to compute
the inner resistance of the metastable cycle in $d$ dimensions,
that is,
the energy barrier that a subcritical
configuration has to overcome in order
to reach either the plus configuration or
the minus configuration in a finite box.
%Fortunately, the attractive features of the Ising model
%allow to avoid these difficulties.
%
This fact affects both strategies for
the upper as well as for the lower estimate of the relaxation time,
since in order
to approximate the distribution of the nucleation time
as an exponential law with rate $\exp(-\b\G_d )$,
one has to rule out the possibility that the process is trapped
in a deep well.
We are able to get the required bounds by
using the attractivity and the reversibility of the dynamics;
see
Lemma \ref{fugaup} and
Proposition \ref{nucl}.%\vspace*{9pt}

\textit{Space--time clusters.}
The third major difficulty to extend the analysis of Dehghanpour and
Schonmann is to control adequately the space--time clusters.
For instance, we cannot proceed as in \cite{DS1} to rule out
the possibility that a subcritical cluster crosses
a long distance.
This question turns out to be much more involved in higher dimension.
It is tackled in Theorem \ref{totcontrole},
which is a key of the whole analysis.
To control the diameters of the space--time clusters, we use
ideas of recurrence and a decomposition of the space into sets
called ``cycle compounds.''
A~cycle compound is a connected set of states $\ocA$
such that the communication
energy between two points
of $\ocA$
is less than or equal to
the communication
energy between
$\ocA$ and its complement. A cycle is a cycle compound, yet an appropriate
union of cycles might form a cycle compound without being a cycle.

%s1.3 #&#
\subsection{Main results}
\label{mr}
We now briefly describe the model, and we state next our main result.
We study the
$d$-dimensional nearest-neighbor stochastic Ising model at inverse
temperature $\beta$
with a fixed
small positive magnetic field $h$, that is,
%More precisely, we study
the continuous-time Markov process
$(\sigma_t)_{t\geq0}$
with state space
$\{-1,+1\}^{\Z^d}$ defined as follows.
In the configuration $\sigma$,
the spin at the site $x\in\Z^d$ flips at rate
\[
c\bigl(\sigma,\sigma^x\bigr) = \exp \bigl(-\beta \bigl(
\Delta_x H(\sigma) \bigr)^+ \bigr),
\]
where $(a)^+=\max(a,0)$ and
\[
\Delta_x H(\sigma) = \sigma(x) \biggl(\mathop{\sum
_{y \in\Z^d}}_{|x-y|=1} \sigma(y) + h \biggr).
\]
In other words,
the infinitesimal generator
of the process $(\sigma_t)_{t\geq0}$
%the Markov process
acts on a local observable $f$ as
\[
(Lf) (\sigma) = \sum_{x} c
\bigl(\sigma,\sigma^x\bigr) \bigl(f\bigl(\sigma ^x
\bigr)-f(\sigma)\bigr),
\]
where $\sigma^x$ is the configuration $\sigma$ in which the spin
at site $x$ has been turned upside down.
Formally, we have
\[
\Delta_x H(\sigma) = H\bigl(\sigma^x\bigr)-H(\sigma),
\]
where $H$ is the formal Hamiltonian given by
\[
H(\sigma) = - \frac{1}{2} \mathop{\sum
_{\{x,y\} \subset\Z^d}}_{|x-y|=1} \sigma(x) \sigma(y) - \frac h 2 \sum
_{x \in\Z^d} \sigma(x).
\]
More details on the construction of this process are given in
Sections~\ref{isisec}
and \ref{sigm}.
%Let us consider the process $(\s_t)_{t\geq0}$ starting from $\minus$.
We denote by
$(\s_t^\minus)_{t\geq0}$ the process starting from $\minus$, the
configuration in which
all the spins are equal to $-1$.
%We give our result for a sufficiently small magnetic field $h$.
A local observable is
a real valued function $f$ defined on the configuration space
which depends only on a finite set of spin variables.
%
%th1.1 #&#
\begin{theorem}
\label{main}
%Let $\G_i$ be the energy of the $i$-dimensional critical
%droplet
%(see Lemma \ref{hcritico} below) and
%let $\k_d$ be defined by \eqref{heur2}.
%Let us consider the Ising model in $\Z^d$ with
Let $f$ be a local observable.
If the magnetic field $h$ is positive and sufficiently small, then
there exists a
%positive
value $\k_d$
%depending only on the dimension $d$ and the magnetic field $h$
such that,
letting
%$$\tau(\beta)=\exp(\beta\k),$$
$\tb=\exp(\beta\k)$,
we have
\begin{eqnarray*}%\label{main1}
\lim_{\beta\to\infty}E\bigl(f\bigl(\s_{{\tb}}^\minus
\bigr)\bigr) &=& f(\minus)  \qquad\mbox{if } \k<\k_d,
\\
\lim_{\beta\to\infty}E\bigl(f\bigl(\s_{\tb}^\minus
\bigr)\bigr) &=& f(\plus)  \qquad\mbox{if } \k>\k_d.
\end{eqnarray*}
The value $\k_d$ depends only on the dimension $d$ and the magnetic
field $h$; in fact, if we denote by $\G_i$ the energy of the
$i$-dimensional critical droplet of the Ising model at zero temperature
and magnetic field $h$, then
\[
\k_d = \frac{1}{d+1} (\G_1+\cdots+\G_d
).
\]
\end{theorem}
%
%The proof yields a slightly precise estimate. Indeed,
%we show the following:
%for each $\e>0$ there exists $\d>0$
%such that
%&&\E(f(\s_t))=f(\minus) + o(\exp({-\b\d})) \mbox{ for any }
%t<\exp({\b(\k_d-\e))}, \\
%&&\E(f(\s_t))=f(\plus) + o(\exp({-\b\d})) \mbox{ for
%any } t>\exp({\b(\k_d+\e))}.
%Theorem \ref{main} (as well as our proof)
%holds also in a box whose side-length
%grows with $\b$ faster than $\exp(\b(\k_d-\k_{d-1}))$.
%For smaller boxes, the time needed to cross the volume
%is smaller than the time needed to see a nucleus,
%so that the typical relaxation time is $\exp({\b\G_d})$
%divided by the volume of the box.
%%\color{red}
%In Section \ref{geom} we introduce some definitions in order to
%In Section \ref{energyestmates} we compute the
Besides the aforementioned technical difficulties,
our proof is
basically an inductive implementation of the scheme of \cite{DS1},
combined with the strategy of \cite{CM}.
%Let us give some insight into the scheme of the proof.
The first step of the proof consists in reducing the problem to
a process defined in a finite exponential volume.
%The results in finite volume imply in a standard way the following
%result for the process in infinite volume.
Let $\kappa>0$, and let $\tau_\beta=\exp(\beta\k)$.
Let $L>\kappa$, and let $\L_\beta=\L(\exp(\beta L))$ be a cubic
box of side length $\exp(\beta L)$.
%Indeed,
%if
%$\kappa<L$ and we set
%,\]
We have that
%then
%
\[
\lim_{\beta\to{\infty}} \P \bigl(f\bigl(\s_{\tau_\beta}^\minus
\bigr) = f\bigl(\s_{\L_\beta,\tau
_\beta}^{-,\minus}\bigr) \bigr) = 1,
\]
where $(\s_{\L_\beta,t}^{-,\minus})_{t\geq0}$ is the process in
the box $\L_\beta$ with minus boundary
conditions starting from $\minus$.
This follows from
a standard large deviation estimate
based on the fact that the maximum rate in the model is $1$;
see Lemmas 1, 2 of \cite{Sch}
for the complete proof.
We state next the finite volume results
that we will prove.
%We recall that
%$(\s_{\L_\beta,t}^{-,\minus})_{t\geq0}$ is the process in the box $\L_
%conditions starting from $\minus$.
%corresponding to this discussion.
%Theorem \ref{T1iv} follows then from
%the finite volume results which are stated below.
%
%th1.2 #&#
\begin{theorem}\label{mainfv}%{(\bf Exponential volume.)}
Let $L>0$, and let $\L_\beta=\L(\exp(\beta L))$ be a cubic
box of side length $\exp(\beta L)$.
Let $\kappa>0$, and let $\tau_\beta=\exp(\beta\k)$.
There exists $h_0>0$ such that, for any $h\in\ ]0,h_0[$, the
following holds:

$\bullet$ If $\kappa<\max(\Gamma_d-dL,\kappa_d)$, then
\[
\lim_{\beta\to\infty} \P \bigl(\s_{\L_\beta,\tau_\beta}^{-,\minus}(0)=1
\bigr) = 0.
\]
%
%and this probability is exponentially small in $\beta$.
%$$\limsup_{\beta\to\infty}

$\bullet$ If $\kappa>\max(\Gamma_d-dL,\kappa_d)$, then
\[
\lim_{\beta\to\infty} \P \bigl(\s_{\L_\beta,\tau_\beta}^{-,\minus}(0)=-1
\bigr) = 0. %\exists  x\in\L_\beta
\]
%
%and the probability of the complement event
%is super--exponentially small in $\beta$.
%and this probability is exponentially small in $\beta$.
%$$\lim_{\beta\to\infty}
\end{theorem}

Recall that $\Ga_d$ and $\kappa_d$ depend on the magnetic field $h$.
Explicit formulas are available for $\Ga_d$ and $\kappa_d$; however,
they are quite complicated. An important point is that $\Ga_d$ and
$\kappa_d$ are continuous functions of the magnetic field $h$ (this is
proved in Lemma \ref{contga}), and this will allow us to reduce the
study to irrational values of $h$. An explicit bound on $h_0$ can also
be computed. In dimension $d$, the proof works if $h_0\leq1$ and Lemma
\ref{control} holds. Let us denote by $m_{d}$ the volume of the
critical droplet in dimension $d$. Lemma \ref{control} holds as soon as
\[
\forall n \leq d\qquad (\Ga_{n-1})^{{n}} \leq (m_{n-1})^{n-1}.
\]
%
%and these inequalities are true for $h$ small enough.
We next shift our attention to finite volumes, and we try to perform
simple computations to understand why the critical constant appearing
in Theorem \ref{mainfv}
is equal to
$\max(\Gamma_d-dL,\kappa_d)$.
%the nucleation and the growth
%of supercritical droplets.
We have two possible scenarios for the relaxation to equilibrium in a
finite cube.
If the cube is small, then
the system relaxes via the formation of a single critical droplet
that grows until covering the entire volume.
If the cube is large, then we have a more efficient mechanism,
creating many critical droplets that grow and eventually coalesce.
The critical side length of the cubes separating these two mechanisms
%The threshold between these two mechanisms
scales exponentially
with $\b$ as $\exp(\b L_d)$, where
\[
L_d = \frac{\Gamma_d-\kappa_{d}}{d}.
\]
%
% =
% \min \Big(
This value is the result of the computations, and we do not have a simple
heuristic explanation for it.
There are three main factors controlling the relaxation time,
which correspond to the heuristics explained previously:\vspace*{9pt}
%The essential facts behind the proof
%of Theorem \ref{T1iv} are the following:

\textit{Nucleation.}
Within a
box of side length $\exp(\beta K)$, the typical time when
the first critical droplet appears is of order
$\exp(\beta(\Gamma_d-dK))$.\vspace*{9pt}

\textit{Initial growth.}
The typical time to grow from a critical droplet (which has a
diameter of order $2d/h$)
into a supercritical
droplet [which has a
diameter of order
$\exp(\b L_d)$]
traveling
at the asymptotic speed
$\exp(-\beta\kappa_{d-1})$
is
$\exp(\beta\Gamma_{d-1})$.\vspace*{9pt}

\textit{Asymptotic growth.}
%The typical time for such a droplet
In a time
$\exp(\beta(K+\kappa_{d-1}))$,
a supercritical
droplet
having a
diameter larger than
$\exp(\b L_d)$
and traveling
at the asymptotic speed
$\exp(-\beta\kappa_{d-1})$
covers a distance
%region of diameter
$\exp(\beta K)$
in each axis direction and its diameter increases by %a factor
$2\exp(\beta K)$.\vspace*{9pt}

The statement concerning the nucleation time contains no mystery.
Let us try to explain the statements on the growth of the droplets.
Once a critical droplet is born, it starts to grow at speed
$\exp(-\beta\Gamma_{d-1})$. As the droplet grows, the speed of growth
increases because the number of choices for the creation of a new
$(d-1)$-dimensional critical droplet
attached to the face of the droplet is of order the surface of the droplet.
Thus the speed of growth of a droplet of size
$\exp(\beta K)$ is
\[
\exp\bigl(\beta\bigl(K(d-1)-\Gamma_{d-1}\bigr)\bigr).
\]
When $K$ reaches the value $L_{d-1}$, the speed of growth is limited
by the inverse of the time needed for the
$(d-1)$-dimensional critical droplet
to cover an entire face of the
droplet. This time corresponds to the $(d-1)$-dimensional
relaxation time in infinite volume,
and the droplet
reaches its asymptotic speed, of order
$\exp(-\beta\kappa_{d-1})$.
The time needed
to grow a critical droplet into a supercritical droplet traveling
at the asymptotic speed
is
\[
\sum_{1\leq i\leq
\exp(\beta L_{d-1})} \exp \biggl(\beta \biggl(
\Gamma_{d-1}- \frac{d-1}{\beta}\ln i \biggr) \biggr)
\]
and,
for $d\geq2$,
this is still of order
$\exp(\beta\Gamma_{d-1})$.
With the help of the above facts, we can estimate
the relaxation time in a
box
%$\Lambda_\beta$
of side length $\exp(\beta L)$.
%Let us now try
%to obtain a lower bound
%on the relaxation time.
%Suppose that we examine the state of the origin
%at a time
%$\exp(\beta\kappa)$.
%The origin becomes occupied when it is covered
%by a droplet. This droplet can result either from the growth of a
%single
%critical droplet or from the coalescence of several droplets.
%Since the speed of propagation of the effects is finite, the state of
%the origin
%at time
%$\exp(\beta\kappa)$
%is unlikely to have been influenced by any event occurring
%outside the box
%of sidelength $\exp(2\beta\kappa)$.
%Thus all the subsequent computations
%can be restricted to this box.
%In particular, a droplet which covers
%the origin
%before time
%$\exp(\beta\kappa)$
%has to be born inside this box, meaning that the oldest site
%of the droplet belongs to this box.
%%$(\sigma_{\L,t})_{t\geq0}$
%Let us consider
%the box $\Lambda_\beta$
%of sidelength $\exp(\beta L)$.
Suppose that the origin is covered by a large supercritical droplet at
time $\exp(\beta\kappa)$. If this droplet is born at distance
$\frac{1}{2}\exp(\beta K)$, then nucleation has occurred inside the box
$\L( \exp(\beta K))$, and the initial critical droplet has grown into a
droplet of diameter $\frac{1}{2}\exp(\beta K)$ in order to reach the
origin. This scenario needs a time
\begin{eqnarray*}
&&
\pmatrix{ %\mbox{time for creating a nucleus}\\
\mbox{time for nucleation}
\cr
\mbox{in the box $\L\bigl(
\exp(\beta K)\bigr)$}} %of sidelength $\exp(\beta L)$
+ %\left(
%time to grow a critical droplet
%}\\
%into a supercritical droplet
%travelling
%}\\
%at the asymptotic speed
%%of sidelength $\exp(\beta L)$
%}
% \right)
% +
\pmatrix{ \mbox{time to cover}
\cr
\mbox{the box $\L\bigl(\exp(
\beta K)\bigr)$}} %of sidelength $\exp(\beta L)$
\\
&&\qquad\sim \exp\bigl(\beta(\Gamma_{d}-dK)\bigr) + \exp(\beta
\Gamma_{d-1}) + \exp\bigl(\beta(K+\kappa_{d-1})\bigr),
\end{eqnarray*}
%
%Thus
%the relaxation time in
%the box $\Lambda_\beta$
%of sidelength $\exp(\beta L)$
%is smaller than
which is of order
\[
\exp \bigl(\beta \max (\Gamma_{d}-dK, \Gamma_{d-1}, K+
\kappa_{d-1} ) \bigr).
\]
To find the most efficient scenario, we
optimize
over $K<L$,
and we
conclude that
the relaxation time in the box
$\L(\exp(\b L))$
%$\L_\beta=\L(\exp(\b L))$
is of order
\[
\exp \Bigl( \beta\inf_{K\leq L} \max (
\Gamma_{d}-dK, \Gamma_{d-1}, K+\kappa_{d-1} ) \Bigr).
\]
%
%Thus the relaxation time is larger than
%$$\displaylines{
%%\min\left(
%%\mbox{time for creating a nucleus}\\
%$\L( 2\exp(\beta K))$
%%of sidelength $\exp(\beta L)$
%}
% +
%time to grow a critical droplet
%into}\\
%a droplet
%of diameter
%$\exp(\beta K)$
%}\\
%%\right)
% \sim
%%\min\Big(
%%\max\Big(
% +
%%\frac{1}{2}
%%\Big)
%%.
%}$$
%%and the relaxation time in
%%the box $\Lambda_\beta$
%%of sidelength $\exp(\beta L)$
%%is larger than
%which is of order
%$$\exp\Big(\beta
%L+\kappa_{d-1}
%)\big)\Big)
%.
%$$
%By optimizing
%over the size of the box $\Lambda_\beta$, we
%conclude that
%the relaxation time in infinite volume
%satisfies
%L+\kappa_{d-1}
%)\big)
%.
%%\begin{eqnarray*}
%%\kappa_d &\leq \inf_L \max\big(\Gamma_{d}-dL,
%%\Gamma_{d-1},
%%L+\kappa_{d-1}
%%\big),\cr
%%\kappa_d &\geq
%%\sup_L
%%\min\big(\Gamma_{d}-dL,
%%\max(
%%\Gamma_{d-1},
%%L+\kappa_{d-1}
%%)\big)
%%.
%%\end{eqnarray*}
It turns out that, for $h$ small,
the above quantity is equal to
\[
\exp \bigl( \beta \max(\Gamma_d-dL,\kappa_d)
\bigr).
\]
In particular, the time needed to grow a critical droplet into
a supercritical droplet is not a limiting
factor for the relaxation whenever $h$ is small.
%Thus the relaxation time is larger than
%$$\displaylines{
%%\min\left(
%Since the optimal value of $L$ solves
%$\G_d-dL=L+\k_{d-1}$,
%the two constants appearing in the exponential
%in the lower and upper bounds for the relaxation time
%coincide, they are equal to
%$$\kappa_d = \max\Big(
%$$
%The easiest
%part of Theorem \ref{T1fv} is the lower bound on the relaxation
%time, i.e., the second case
%where $\kappa>
%
%Let us give a sketch of the proof of
%s1.4 #&#
\subsection{Strategy of the proof}
\label{sp}
The upper bound on the relaxation
time, that is, the second case
where $\kappa>
\max(\Gamma_d-dL,\kappa_d)$, is done in Section~\ref{relaxa}.
The ingredients involved in the upper bound are known since
the works of Neves, Dehghanpour and Schonmann,
and this part is considerably
easier than the lower bound.
The hardest part of Theorem \ref{mainfv} is the lower bound on the relaxation
time, that is, the first case
where $\kappa<
\max(\Gamma_d-dL,\kappa_d)$.
%We prove first that, up to a SES event, a large droplet is not formed
%before a time
%$\exp(\beta\Gamma_{d-1})$.
The lower bound is done in Sections~\ref{stc} and \ref{metare}.
Let us explain the strategy of the proof of
the lower bound,
%Theorem \ref{mainfv},
without stating precisely the definitions
and the technical results.

Let $L>0$ and let $\L_\beta=\L(\exp(\beta L))$ be a cubic
box of side length $\exp(\beta L)$.
Let $\kappa>0$ and let $\tau_\beta=\exp(\beta\k)$.
We want to prove that it is unlikely that the spin at the origin is equal
to $+1$ at time $\tb$ for the process
$(\s_{\L_\beta,t}^{-,\minus})_{t\geq0}$.
Throughout the proof, we use in a crucial way the notion of
space--time cluster.
A space--time cluster of the
trajectory $(\sigma_{\L,t},{0\leq t\leq\tb})$
is a maximal connected component of space--time points
for
the following relation:
two space--time points $(x,t)$ and $(y,s)$
are connected if $\sigma_{\L,t}(x) = \sigma_{\L,s}(y) = +1$
and
%$\bullet$
either ($s=t$ and $|x-y| \le1$)
%$\bullet$
or [$x=y$ and $\sigma_{\L,u}(x)=+1$ for $s\leq u\leq t$].
%or ($x=y$ and $\sigma_{\L,u}(x)=+1$ for $u\in[\min(s,t),\max(s,t)]$).
%A space--time cluster of the trajectory $(\sigma_{\L,t})_{t\geq0}$
%is a maximal connected component of space--time points.
%For $u\leq s\in\R^+$, we denote by $\STC(u,s)$ the space--time clusters
%of the trajectory restricted to the time interval $[u,s]$.
%
%The space--time cluster ($\STC$) at time $t$,
With the
space--time clusters, we record the influence of the plus spins
throughout the evolution. We can then compare the status of a spin in
dynamics associated to different boundary conditions with the help
of the graphical construction
(described in Section~\ref{sigm}).
The diameter
$\diam \cC$ of a space--time cluster $\cC$
is the diameter of its spatial projection.
We argue as follows.
If
$\s_{\L_\beta,\tb}^{-,\minus}(0)=+1$, then
the space--time point $(0,\tb)$ belongs to
a nonvoid space--time cluster, which we denote by $\cC^*$.
We discuss then according to the diameter of~$\cC^*$.

$\bullet$ If
$\diam \cC^*<\ln\ln\beta$, then $\cC^*$ is also
a space--time cluster of the process
$ (\s_{\L(\ln\beta),t}^{-,\minus}, 0\leq t\leq\tb )$,
and the spin\vspace*{1pt} at the origin is also equal to $+1$ in this process
at time $\tb$.
%${\s_{\L(\ln\beta),\tb}^{-,\minus}}$.
%By a straightforward extension of
The finite volume estimates obtained for
fixed boxes can be readily extended to boxes of side length $\ln\beta$,
and we obtain that the probability of the above event is
exponentially small if $\kappa<\Ga_d$, because the
entropic contribution to the free energy is negligible with
respect to the energy.

$\bullet$
If $\diam \cC^*>\exp(\beta L_d)$ (this case can occur only when $L>L_d$),
then we use the main technical estimate of the paper,
Theorem \ref{T2}, which states roughly the following:
for $\kappa<\kappa_d$,
the probability that,
%for the process
%{$(\s^{n\pm,\xi}_{\Sigma,t})_{t\geq0}$},
in the trajectory \mbox{$ (\s_{\L_\beta,t}^{-,\minus}, 0\leq t\leq\tb
)$,}
%a $\STC$ in
%$\STC_\xi(0,\tb)$
%a space--time cluster
there exists
a space--time cluster
of diameter
larger than
$\exp({\b L_d})$
%is created
%%before time
%in the process
%{$(\s^{n\pm,\xi}_{\Sigma,t}, 0\leq t\leq
%)$}
is a super exponentially small function of $\beta$ (denoted by $\SES$
in the following), and it can be neglected.
%of Section \ref{fvs}

$\bullet$
If $\ln\ln\beta\leq\diam \cC^*\leq\exp(\beta L_d)$, then
$\cC^*$ is also
a space--time cluster of the process
restricted to the box
$\L(3\exp(\beta L_d))\cap\L_\beta$.
%$\big(
%{\s_{\L(\ln\beta),t}^{-,\minus}}, 0\leq t\leq\tb\big)$,
A space--time cluster
is said to
be large
if its
diameter is larger than or equal to $\ln\ln\beta$.
A box is said to be small if its
sides have a length larger than $\ln\ln\beta$ and smaller than
$d\ln\beta$.
The diameters
of the space--time clusters increase with time when they coalesce
because of a spin flip. This implies that
if a large space--time cluster is created
in the box $\L_\beta$, then it has to be created first locally in a
small box.
The number of small boxes included in $\L_\beta$ is of order
\[
\bigl| \L\bigl(3\exp(\beta L_d)\bigr)\cap\L_\beta \bigr| = \exp
\bigl(\beta d\min(L_d,L) \bigr).
\]
For the dynamics restricted to a small box, we have
\begin{eqnarray*}
P \pmatrix{ %\mbox{for }i\in I,
\mbox{a large $\STC$ is}
\cr
\mbox{created before $
\tb$} %\mbox{in the process }
%(\s^{-,\minus}_{
%%S_\beta
%%\L(3\exp(\beta L_d))
%Q
%,t}, t\geq0)
}
&\leq& P \pmatrix{ %\mbox{for }i\in I,
\mbox{a large $\STC$ is created}
\cr
\mbox{before
nucleation} %\mbox{in the process }
%(\s^{-,\minus}_{
%%S_\beta
%%\L(3\exp(\beta L_d))
%Q
%,t}, t\geq0)
} \\
&&{}+ P \pmatrix{ %\mbox{for }i\in I,
\mbox{nucleation occurs}
\cr
\mbox{before
$\tb$} %\mbox{in the process }
%(\s^{-,\minus}_{
%%S_\beta
%%\L(3\exp(\beta L_d))
%Q
%,t}, t\geq0)
}.
\end{eqnarray*}
The main result of Section~\ref{diamostc},
Theorem \ref{totcontrole}, yields that
the first term of the right-hand side is $\SES$.
The finite volume estimates in fixed boxes obtained in the
previous studies of metastability can be readily extended to small
boxes.
By Lemma~\ref{fugaup},
we have that, up to corrective factors,
\[
P \pmatrix{ %\mbox{for }i\in I,
\mbox{nucleation occurs}
\cr
\mbox{before $\tb$}
%(\s^{-,\minus}_{
%%S_\beta
%%\L(3\exp(\beta L_d))
%Q
%,t}, t\geq0)
} \leq
\tau_\beta \exp(-\b\G_{d}). % + \SES,
\]
%
%|Q|^{m_d+1}\tau_\beta
% \exp(-\b\G_{d})
% + \SES
Finally, we have
\begin{eqnarray*}
\P \bigl(\diam\cC^*\geq\ln\ln\beta \bigr) &\leq& \exp \bigl(\beta d
\min(L_d,L) \bigr) \bigl( \tau_\beta \exp(-\b
\G_{d}) + \SES \bigr)
\\
&\leq& \exp \bigl(\beta \bigl(d\min(L_d,L) +\kappa-\G_{d}
\bigr) \bigr) + \SES
\\
&=& \exp \bigl(\beta \bigl( \kappa- \max(
\G_d-dL_d,\G_d-dL) \bigr) \bigr) + \SES,
\end{eqnarray*}
and the desired result follows easily.

From this quick sketch of proof, we see that the most difficult
intermediate results are
Theorems \ref{totcontrole} and \ref{T2}.
The remainder of the paper is mainly devoted to the proof of these results.
%It follows that
%$$\limsup_{\b\to\infty}
%dL_d+\k-\G_d
% < 0
%,$$
%and we are done!!
In Section~\ref{metro}, we consider a general Metropolis dynamics on a
finite state space, we recall the formulas for the law of exit in
continuous time and we introduce the notions of cycle and cycle
compound in this context. Section~\ref{feat} is devoted to the study of
some specific features of the cycle compounds of the Ising model. In
Section~\ref{energyestmates}, we state several discrete isoperimetric
results from \mbox{\cite{N2,N1,AC}}.
%and the fundamental estimate for the nucleation time in a finite box.
Apart from the notion of cycle compound, the definitions and the
results presented in Sections~\ref{metro}, \ref{feat} and
\ref{energyestmates} come from the previous literature on
metastability, with some rewriting and adaptation to fit the
continuous-time framework and our specific $n\pm$ boundary conditions.
The main technical contributions of this work are presented in
Sections~\ref{stc} and \ref{metare}. In Section~\ref{stc}, we prove the
key estimate on the diameters of the space--time clusters (Theorem
\ref{totcontrole}). Section~\ref{metare} is devoted to the proof of
Theorem \ref{T2}. The proof of the lower bound on the relaxation time
is completed in Section~\ref{concl}. The final section,
Section~\ref{relaxa}, contains the proof of the upper bound on the
relaxation time.

%s2 #&#
\section{The Metropolis dynamics}
\label{metro}

A very efficient tool for describing the meta\-stable behavior of a
process in the low temperature regime is a hierarchical decomposition
of the state space known as the cycle decomposition. In the context of
a Markov chain with finite state space evolving under a Metropolis
dynamics, the cycles can be defined geometrically with the help of the
energy landscape. Our context of infinite volume is much more
complicated, but since the system is attractive, we will end up with
some local problems that we handle with the finite volume techniques.
We start by reviewing these techniques. Here we recall some basic facts
about the cycle decomposition. For a complete review we refer to
\cite{S4,OS3,OS2,CaCe,OS1,OV}. Since we are working here with a
continuous-time process defined with the help of transition rates, as
opposed to a discrete-time Markov chain defined with transition
probabilities, we feel that it is worthwhile to present the exact
formulas giving the law of exit of an arbitrary subset in this slightly
different framework. This is the purpose of Section~\ref{lawexit}. In
Section~\ref{metrosec}, we define the Metropolis dynamics, and we show
how to apply the formulas of Section~\ref{lawexit} to this specific
dynamics. In Section~\ref{cycles}, we recall the definitions of a
cycle, the communication energy, the height of a set, its bottom, its
depth and its boundary. We introduce also an additional concept, called
cycle compound, which turns out to be useful when analyzing the energy
landscape of the Ising model. Apart from the notion of cycle compounds,
the definitions and the results presented in this section come from the
previous literature on metastability and simulated annealing, they are
simply adapted to the continuous-time framework.

%s2.1 #&#
\subsection{Law of exit}
\label{lawexit}
We will not derive in detail all the
results used in this paper concerning the behavior of a
Markov process with exponentially vanishing transition rates because
the proofs are essentially the same as in the discrete-time setting.
These proofs can be found in the book of Freidlin and Wentzell
(\cite{FW}, Chapter $6$, Section $3$), or in the lecture
notes of Catoni (\cite{Catoni}, Section $3$).
However, for the sake of clarity, we present the two basic formulas
in continuous time giving the law of the exit from an arbitrary set.
Let $\cX$ be a finite state space.
Let $c\dvtx \cX\times\cX\to\R$ be a matrix of transition rates
on $\cX$,
that is,
\begin{eqnarray*}
&\displaystyle \forall x,y \in\cX,\qquad x\neq y,\qquad c(x,y)\geq0,&
\\
&\displaystyle \forall x\in\cX\qquad \sum_{y\in\cX}c(x,y) = 0.&
\end{eqnarray*}
%
%$${\lim_{\beta\to\infty}}  {\beta}^{-1}\ln a(\beta) = 0$$
%and $V\dvtx E\times E\to\mathbb R_+\cup\{\infty\}$ is an irreducible cost
%function
%i.e.
%$${
%x_0=y,  x_r=z,\cr
%c(x_0,x_1)\times\cdots\times c(x_{r-1},x_r) > 0.}$$
We consider the
continuous-time homogeneous Markov
process $(X_t)_{t\geq0}$
%$(X, P)$
on $\cX$
whose infinitesimal generator is
\[
\forall f\dvtx \cX\to\R\qquad (Lf) (x) = \sum_{y \in\cX}
c(x,y) \bigl(f(y)-f(x)\bigr).
\]
For $C$ an arbitrary subset of $\cX$,
we define the time $\tau(C)$
of exit from $C$
\[
\tau(C) = \inf\{ t\geq0\dvtx X_t\notin C \}.
\]
The next lemmas provide useful formulas
%for the invariant measure and
for
the laws of the exit time and exit point for an arbitrary subset
of $\cX$.
These formulas are rational fractions of products of the coefficients
of the
matrix of the transition rates,
whose numerators and denominators are most conveniently
written as sums over particular types of graphs.

%de2.1 #&#
\begin{definition}[{[The graphs $G(W)$]}]
Let $W$ be an arbitrary nonempty subset of $\cX$.

An oriented graph on $\cX$ is called a $W$-graph if and only
if:

$\bullet$  there is no arrow starting from a point
of $W$;

$\bullet$  each point of $W^c$ is the initial point of
exactly one arrow;

$\bullet$  for each point $x$ in $W^c$, there exists a
path in the
graph leading from $x$ to $W$.

The set of all $W$-graphs is denoted by $G(W)$.
\end{definition}
If the first two conditions are fulfilled, then
the third condition above is equivalent to:

$\bullet$  there is no cycle in the graph.

%de2.2 #&#
\begin{definition}[{[The graphs $G_{x,y}(W)$]}]
Let $W$ be an arbitrary nonempty subset of $\cX$,
let $x$ belong to $\cX$ and $y$ to $W$.
If $x$ belongs to $W^c$, then the set $G_{x,y}(W)$ is the set of
all oriented graphs on $\cX$ such that:

$\bullet$  there is no arrow starting from a point
of $W$;

$\bullet$  each point of $W^c$ is the initial point of
exactly one arrow;

$\bullet$  for each point $z$ in $W^c$, there exists a
path in the
graph leading from $z$ to $W$;

$\bullet$  there exists a path in the
graph leading from $x$ to $y$.

More concisely, they are the graphs of $G(W)$ which
contain a path leading from $x$ to $y$.

If $x$ belongs to $W$, then
the set $G_{x,y}(W)$ is empty if $x\neq y$
and is equal to $G(W)$ if $x=y$.
\end{definition}

The graphs in $G_{x,y}(W)$ have no cycles.
For any $x$ in $\cX$ and $y$ in $W$, the set $G_{x,y}(W)$ is included
in $G(W)$.

%de2.3 #&#
\begin{definition}[{[The graphs $G(x\nrightarrow W)$]}]
Let $W$ be an arbitrary, nonempty subset of $\cX$, and let $x$ be a point
of $\cX$.

If $x$ belongs to $W$, then the set $G(x\nrightarrow W)$ is
empty.

If $x$ belongs to $W^c$, then
the set $G(x\nrightarrow W)$ is the set of
all oriented graphs on $\cX$ such that:

$\bullet$  there is no arrow starting from a point
of $W$;

$\bullet$  each point of $W^c$ except one, say $y$,
is the initial point of exactly one arrow;

$\bullet$  there is no cycle in the graph;

$\bullet$  there is no path in the graph leading from $x$
to $W$.

The third condition (no cycle) is equivalent to:

$\bullet$  for each $z$ in $W^c\setminus\{y\}$,
there is a path in the graph leading from $z$ to $W\cup\{y\}$.
\end{definition}

%le2.4 #&#
\begin{lemma}\label{aaa}
Let $W$ be an arbitrary, nonempty subset of $\cX$, and let $x$ be a point
of $\cX$.
The set $G(x\nrightarrow W)$ is the union of all the sets
$G_{x,y}(W\cup\{y\}),\break y\in W^c$.
\end{lemma}

In the case $x\in W^c, y\in W$, the definitions
of $G_{x,y}(W)$ and $G(x\nrightarrow W)$ are those
given by Wentzell and Freidlin (1984). We have extended these definitions
to cover all possible values of $x$. With our choice for the
definition of the time of exit $\tau(W^c)$ (the first time greater than
or equal to zero when the chain is outside~$W^c$),
the formulas for the law of $X_{\tau(W^c)}$
and for the expectation of $\tau(W^c)$ will remain valid in all cases.

Let $g$ be a graph on $\cX$. We define
\[
c(g) = \prod_{(x\rightarrow y)\in g} c(x,y).
\]

%le2.5 #&#
\begin{lemma}[(Exit point)]\label{dfe}
For any nonempty subset $W$ of $\cX$, any $y$ in $W$ and $x$ in $\cX$,
\[
P(X_{\tau(W^c)}=y/ X_0=x) = \frac{ \sum_{g\in G_{x,y}(W)}c(g)}{
\sum_{g\in G(W)} c(g)}.
\]
\end{lemma}

%
%le2.6 #&#
\begin{lemma}[(Exit time)]\label{ert}
For any subset $W$ of $\cX$ and $x$ in $\cX$,
\[
E\bigl(\tau\bigl(W^c\bigr)/ X_0=x\bigr) =
\frac{ \sum_{y\in W^c}
\sum_{g\in G_{x,y}(W\cup\{y\})}c(g)} {
\sum_{g\in G(W)} c(g)} = \frac{ \sum_{g\in G(x\nrightarrow W)}c(g)} {
\sum_{g\in G(W)} c(g)}.
\]
\end{lemma}

For instance, if we apply Lemma \ref{ert} to the case where
$W=\cX\setminus\{x\}$, and the process starts from $x\in\cX$,
then we get
\[
E\bigl(\tau\bigl(\{ x \}\bigr)/ X_0=x\bigr) = \frac{1}{\sum_{y\neq x}c(x,y)} = -
\frac{1}{c(x,x)}.
\]
To prove these formulas in continuous time, we
study the involved quantities as functions of the starting point and
we
derive a system of linear equations
with the help of the Markov property.
For instance, let
\[
m(x,y) = P(X_{\tau(W^c)}=y/ X_0=x).
\]
Let $T=\tau(\{ x \})$. We have then
\begin{eqnarray*}
m(x,y) &=& \sum_{z\in W^c} P(X_{\tau(W^c)}=y,
X_T=z/ X_0=x) + P(X_{T}=y/
X_0=x)
\\
&=& \sum_{z\in W^c} P(X_{\tau(W^c)}=y/
X_0=z) P(X_{T}=z/ X_0=x)\\
&&{} +
P(X_{T}=y/ X_0=x).
\end{eqnarray*}
Let
\[
p(x,z) = P(X_{T}=z/ X_0=x) = \frac{c(x,z)}{\sum_{u\neq x}c(x,u)} = -
\frac{c(x,z)}{c(x,x)}.
\]
Then $p(\cdot,\cdot)$ is a matrix of transition probabilities,
and
\[
m(x,y) = \sum_{z\in W^c} p(x,z) m(z,y) +p(x,y).
\]
This is exactly the same equation as in the case of a discrete-time
Markov chain with transition matrix
$p(\cdot,\cdot)$.
This way
the continuous-time formula can be deduced from its discrete-time counterpart.

%s2.2 #&#
\subsection{The Metropolis dynamics}
\label{metrosec}
We suppose from now onward that
we are dealing
with a family of continuous-time homogeneous Markov
processes $(X_t)_{t\geq0}$ indexed by a positive
parameter $\beta$ (the inverse temperature).
In particular,
the state space and the transition rates change with $\beta$.
We suppose that these processes evolve under a Metropolis dynamics.
More precisely, let
$\alpha\dvtx \cX\times\cX\to[0,1]$ be a symmetric
irreducible transition kernel on $\cX$,
that is,
$\alpha(x,y)=\alpha(y,x)$ for $x,y\in\cX$ and
%$${\lim_{\beta\to\infty}}  {\beta}^{-1}\ln a(\beta) = 0$$
%and $V\dvtx E\times E\to\mathbb R_+\cup\{\infty\}$ is an irreducible cost
%function
%i.e.
%
\begin{eqnarray*}
&&\forall y,z\in\cX\times\cX\qquad \exists x_0,x_1,\ldots,x_r, x_0=y, x_r=z,
\\
&&\qquad\alpha(x_0,x_1)\times\cdots\times
\alpha(x_{r-1},x_r) > 0.
\end{eqnarray*}
Let $H\dvtx \cX\to\R$ be an energy defined on $\cX$.
%More precisely,
%we consider the coordinate process $X=(X_n)_{n\in\Bbb N}$ on the
%space ${E^{\Bbb N}}$ defined by
%$X_n\dvtx (\omega_0,\ldots)\mapsto\omega_n$ together with a family of
%probabilities $(P_\beta)$ indexed by $\beta$; under each of them
%the coordinate process is a Markov chain.
%We suppose that these Markov chains follow the Me\-tro\-po\-lis
%dynamics
%associated to the kernel $\alpha$ and the energy $H$,
%namely that their infinitesimal generator is
%$$
%(Lf)(x) =  \sum_{y \in\cX} c(x,y) (f(y)-f(x)),
%$$
%chains are in the Freidlin--Wentzell regime, namely that their
%transition mechanisms satisfy
%$$a(\beta) \exp-\beta V(x,y)  \leq  P_\beta(X_{n+1}=y/X_n=x)
We suppose that
the transition rates $c(x,y)$ are given by
\[
\forall x,y\in\cX\qquad c(x,y) = \alpha(x,y) \exp \bigl( -\b\max \bigl(0,H(y)-H(x)
\bigr) \bigr).
\]
%
%for all $x,y$ in $\cX$, where $\beta\mapsto a(\beta)$
%is a positive function such that
%For $C$ an arbitrary subset of $\cX$ and $t\geq0$,
%we define the time $\tau(C,t)$
%of exit from $C$ after time $t$
%$$\tau(C,t) = \inf\{ s\geq t\dvtx X_s\notin C \}$$
%(we make the convention that $\tau(C)=\tau(C,0)$).

%We define also the time $\theta(C,t)$ of the last visit to
%the set $C$ before time $t$
%$$\theta(C,t) = \sup\{ s\leq t\dvtx X_s\in C \}$$
%(if the chain has not visited $C$ before $t$, we take $\theta(C,t)=0$).
%Remark that $\tau$ is a stopping time when $m$ is deterministic
%whereas $\theta$ is not.\par
%describe the probability
% of a given configuration when the
%nonnucleating process starts from the metastable state.
The irreducibility hypothesis ensures the existence of a
unique
invariant probability measure $\nu$ for the
Markov process $(X_t)_{t\geq0}$.
We have then, for any $x,y\in\cX$ and $t\geq0$,
\[
\nu(x)P(X_t=y/X_0=x) \leq \sum
_{z\in\cX}\nu(z)P(X_t=y/X_0=z) = \nu(y).
\]
In the case where $\alpha(x,y)\in\{ 0,1 \}$ for $x,y\in\cX$,
the invariant measure $\nu$ is the Gibbs distribution
associated to the Hamiltonian $H$ at inverse temperature $\beta$,
and we have
\[
\forall x,y\in\cX, \forall t\geq0\qquad P(X_t=y/X_0=x) \leq
\exp \bigl( -\b\bigl(H(y)-H(x)\bigr) \bigr).
\]
We will send $\beta$ to $\infty$, and we seek asymptotic estimates on
the law of exit from a subset of $\cX$. The exact formulas given in the
previous section can be exploited when the
cardinality of the space $\cX$
and the degree of the communication graph are not too large,
so that the number of terms in the sums is negligible on
the exponential scale. More precisely, let
$\deg(\alpha)$ be the degree of the communication kernel $\alpha$,
that is,
\[
\deg(\alpha) = \max_{x\in\cX} \bigl| \bigl\{ y\in\cX\dvtx \alpha(x,y)>0
\bigr\} \bigr|.
\]
We suppose that $\alpha(x,y)\in\{ 0,1 \}$ for $x,y\in\cX$ and that
\[
\lim_{\beta\to\infty} \frac{1}{\beta} {|\cX| \ln\deg(\alpha)} = 0.
\]
Under this hypothesis, for any subset $W$ of $\cX$, the number of
graphs in $G(W)$ is bounded by
\[
\bigl| G(W) \bigr| \leq {\deg(\alpha)}^{|\cX|} = \exp o(\beta).
\]
From Lemma \ref{dfe},
we have then
for a subset $W$ of $\cX$, $y$ in $W$ and $x$ in $\cX$,
\[
{\deg(\alpha)}^{-|\cX|} \frac{c(g^*_{x,y})}{
c(g^*_{W})} \leq P(X_{\tau(W^c)}=y/
X_0=x) \leq {\deg(\alpha)}^{|\cX|} \frac{c(g^*_{x,y})}{c(g^*_{W})},
\]
where the graphs
$g^*_{x,y}$ and
$g^*_{W}$ are chosen so that
\begin{eqnarray*}
{c\bigl(g^*_{x,y}\bigr)} &=& \max \bigl\{ c(g)\dvtx g\in
G_{x,y}(W) \bigr\},
\\
{c\bigl(g^*_{W}\bigr)} &=& \max \bigl\{ c(g)\dvtx g\in G(W) \bigr\}.
\end{eqnarray*}
For $g$ a graph over $\cX$ we set
\[
V(g) = \sum_{(x\rightarrow y)\in g} \max \bigl(0,H(y)-H(x) \bigr)
\]
so that $c(g)=\exp(-\beta V(g))$.
The previous inequalities yield then
\begin{eqnarray*}
&&\lim_{\beta\to\infty} \frac{1}{\beta} \ln P(X_{\tau(W^c)}=y/
X_0=x)
\\
&&\qquad=\min \bigl\{ V(g)\dvtx g\in G_{x,y}(W) \bigr\} - \min \bigl\{ V(g)\dvtx g\in
G(W) \bigr\}.
\end{eqnarray*}
Similarly, from Lemma \ref{ert}, we obtain that
\begin{eqnarray*}
&&\lim_{\beta\to\infty} \frac{1}{\beta} \ln E\bigl({\tau
\bigl(W^c\bigr)}/ X_0=x\bigr)
\\
&&\qquad=\min \bigl\{ V(g)\dvtx {g\in G( x\nrightarrow W )} \bigr\} - \min \bigl\{ V(g)\dvtx g
\in G(W ) \bigr\}.
\end{eqnarray*}

%s2.3 #&#
\subsection{Cycles and cycle compounds}
\label{cycles}

We say that two states $x,y$ communicate
if either $x=y$
or $\alpha(x,y)>0$.
%(in a more general setting, two configurations are communicating if the
%corresponding transition probability is positive).
A path $\o$ is a sequence $\o=(\o_1,\ldots,\o_n)$
of states such that each state of the sequence communicates
with its successor.
%We denote the set of paths starting from $x$ and
%ending in $y$ by $\tonda{x \to y}$.
A set $\cA$ is said to be connected if any states in $\cA$ can be
joined by a path
in $\cA$, that is,
\begin{eqnarray*}
&&\forall x,y\in\cA\qquad \exists \o_1,\ldots,\o_n\in\cA,
\o_1=x, \o_n=y,
\\
&&\qquad\alpha(\o_1,\o_2)\cdots\alpha(
\o_{n-1},\o_n) %\alpha(\o_{n},y)
>0.
\end{eqnarray*}
We define the communication energy between two states
$x, y$ by
\[
E(x,y)=\min \Bigl\{ \max_{z \in\o} H(z)\dvtx  \o\mbox{ path from $x$
to $y$} \Bigr\}.
\]
%
%Notice that for any configuration $\h$, we have
%$E(\h,\h)=-\infty$, because the set of paths
The
communication energy between two
sets of states $\cA,\cB$ is
\[
E(\cA,\cB)=\min \bigl\{ E(x,y)\dvtx  {x\in\cA,y\in\cB} \bigr\}.
\]
The height of a set of states $\cA$ is
\[
\height(\cA) = \max \bigl\{ %H(x)\dvtx  x\in\cA \big\}.$$
E(x,y)\dvtx x,y\in\cA, x\neq y \bigr\}.
\]

%de2.7 #&#
\begin{definition}\label{cyclesb}
A cycle is a connected set of states $\cA$ such that
\[
%H(x)
\height(\cA) < E(\cA,\cX\setminus\cA).
\]
%
%$$\max \big\{
%E(\sigma,\eta)\dvtx
%
%E(\sigma,\cA^c)\dvtx
%$$
%A cycle is a set of configurations $\cA$ such that
%for some $\lambda\in{\mathbb R}$,
A cycle compound is a connected set of states
$\ocA$ such that
\[
%H(x)
\height(\ocA) \leq E(\ocA,\cX\setminus\ocA).
\]
%
%$$\max \big\{
%E(\sigma,\eta)\dvtx
%
%E(\sigma,\cA^c)\dvtx
%$$
%A cycle compound is a set of configurations
%$\ocA$ such that
%for some $\lambda\in{\mathbb R}$,
\end{definition}
Let us rewrite these definitions directly in terms of the energy $H$.
For any set $\cA$, we have
\[
E(\cA,\cX\setminus\cA)=\min \bigl\{ \max\bigl(H(x),H(y)\bigr)\dvtx  {x\in\cA, y\notin\cA}, \alpha(x,y)>0 \bigr\}.
\]
Notice that the height of a singleton is $-\infty$. Moreover, if $\cA$
is a connected set having at least two elements, then
\[
\height(\cA) = \max \bigl\{ H(x)\dvtx  x\in\cA \bigr\}.
\]
%
%$$
%H(x)
%%\height(\cA)
% <
%E(\cA,\cX\setminus\cA).$$
Thus a cycle is either a singleton or a connected set of states $\cA$ such
that
\[
\forall x,y\in\cA, \forall z\notin\cA\qquad \alpha(y,z)>0 \quad\Longrightarrow\quad H(x)
< % \min \big\{
\max\bigl(H(y),H(z)\bigr).
%E(\cA,\cX\setminus\cA)
\]
%
%$$
%H(x)
%%\height(\cA)
% <
%E(\cA,\cX\setminus\cA).$$
A cycle compound is either a singleton or a connected set of states
$\ocA$ such
that
\[
\forall x,y\in\ocA, \forall z\notin\ocA\qquad \alpha(y,z)>0 \quad\Longrightarrow\quad H(x)
\leq % \min \big\{
\max\bigl(H(y),H(z)\bigr). %\dvtx  {x\in\cA, y\notin\cA},
%E(\cA,\cX\setminus\cA)
\]
Although a cycle and a cycle
compound have almost the same definitions, the structure of these sets
is quite different.
Indeed, the communication under a fixed height $\lambda$
is an equivalence relation, and the cycles are equivalence classes
under this relation.
In particular, two cycles are either disjoint or included one into the other.
With our definition, any singleton is also a cycle of
height $-\infty$.
The next proposition shows that a cycle compound
can have a more complicated structure.
%
%pr2.8 #&#
\begin{proposition}\label{union}
Let $n\geq2$, and
let $\cA_1,\ldots, \cA_n$ be $n$ cycles such that
\[
E(\cA_1,\cX\setminus\cA_1) = \cdots = E(
\cA_n,\cX\setminus\cA_n).
\]
If their union
\[
\overline{\cA} = \bigcup_{i=1}^n
\cA_i %\big\{ x\in\partial\cA_i\dvtx H(x)=E(\cA_i,\cX\setminus\cA_i)
%%\big\}\Big)$$
\]
is connected, then it is a cycle compound.
\end{proposition}
\begin{pf}
If $\overline{\cA}$ is a singleton, then there is nothing to prove.
Let us
suppose that
$\overline{\cA}$
has at least two elements.
Since
$\overline{\cA}$ is connected, then
\[
\height(\overline{\cA}) = \max \bigl\{ H(x)\dvtx  x\in\overline{\cA} \bigr\}.
\]
Moreover,
\[
E(\overline{\cA},\cX\setminus\overline{\cA}) \geq \min_{1\leq i\leq n}
E(\cA_i,\cX\setminus\cA_i) = \max_{1\leq i\leq n}
E(\cA_i,\cX\setminus\cA_i).
\]
For $i\in\{ 1,\ldots,n \}$, since $\cA_i$ is a cycle,
we have
\[
E({\cA_i},\cX\setminus{\cA_i}) \geq \max \bigl\{ H(x)\dvtx  x
\in{\cA_i} \bigr\},
\]
whence
\[
E(\overline{\cA},\cX\setminus\overline{\cA}) \geq \max_{1\leq i\leq n}
\max \bigl\{ H(x)\dvtx  x\in{\cA_i} \bigr\} = \height(\overline{\cA}),
\]
so that $\overline{\cA}$ is a cycle compound.
\end{pf}

Thus two distinct cycle compounds might have a nonempty intersection.
%A set $\overline{\cA}$ is a cycle compound if it is a connected
%set for the communication kernel $\alpha$ and if it is of the form
%%where $\cA_1,\ldots, \cA_n$ are cycles.
Let us introduce a few more definitions.
The bottom of a set $\cG$ of states is
\[
\bottom(\cG) = \Bigl\{ x\in\cG\dvtx  H(x) = \min_{y\in\cG}H(y) \Bigr\}.
\]
It is the set of the minimizers of the energy in $\cG$.
We denote the energy of the states in
$\bottom(\cG)$ by $H(\bottom(\cG))$.
The depth of a set $\cG$ is
\[
\depth(\cG) = E(\cG,\cX\setminus\cG)-H\bigl(\bottom(\cG)\bigr).
\]
The exterior boundary of a subset
$\cG$ of $\cX$ is the set
\[
\partial\cG= \bigl\{x \notin\cG\tc\exists y \in\cG, \alpha(y,x)>0 \bigr\}.
\]
Let us set, for $g$ a graph over $\cX$,
\[
V(g) = \sum_{(x\rightarrow y)\in g} \max \bigl(0,H(y)-H(x) \bigr).
\]
The following results are far from obvious; they are consequences
of the formulas of Section~\ref{lawexit} and the analysis of the cycle
decomposition
\cite{S4,OS3,OS2,CaCe,OS1}.
%
%th2.9 #&#
\begin{theorem}\label{exitcost}
Let $\ocA$ be a cycle compound, let $x\in\ocA$ and
let $y\in\partial\ocA$. We have the identity
\begin{eqnarray*}
&&\min \bigl\{ V(g)\dvtx g\in G_{x,y}(\cX\setminus\ocA) \bigr\} - \min \bigl
\{ V(g)\dvtx g\in G(\cX\setminus\ocA) \bigr\}
\\
&&\qquad= \max \bigl(0,H(y)- E(\ocA,\cX\setminus\ocA) \bigr),
\\
%$${
&&\min \bigl\{ V(g)\dvtx {g\in G(x\nrightarrow \cX\setminus\ocA)} \bigr
\} - \min \bigl\{ V(g)\dvtx g\in G(\cX\setminus\ocA) \bigr\}
\\
&&\qquad= E(\ocA,\cX\setminus\ocA) %\height(\ocA)
-H\bigl(\bottom(\ocA)\bigr).
\end{eqnarray*}
\end{theorem}

Substituting the above identities into the formulas of
Lemmas \ref{dfe} and \ref{ert}, we obtain the following estimates.
%
%co2.10 #&#
\begin{corollary}\label{exitcom}
%For any cycle compound $\ocA$ of $\cX$,
%we have for any
%$x\in\ocA$,
%$y\in\partial\ocA$,
%$$\forall x\in\ocA
%
Let $\ocA$ be a cycle compound, let $x\in\ocA$ and
let $y\in\partial\ocA$. We have
\begin{eqnarray*}
{\deg(\alpha)}^{-|\cX|} &\leq& \frac{P (
X_{\tau(\ocA)}=y/X_0=x )} {
\exp (-\beta
\max (0,H(y)-
E(\ocA,\cX\setminus\ocA)
) )} \leq {\deg(
\alpha)}^{|\cX|},
\\
%For any cycle compound $\ocA$ of $\cX$,
%we have for any
%$x\in\ocA$,
%$y\in\partial\ocA$,
%$$\forall x\in\ocA
%
%Let $\ocA$ be a cycle compound and let $x\in\ocA$.
%We have
{\deg(
\alpha)}^{-|\cX|} &\leq& \frac{
E(\tau(W^c)/ X_0=x)}{
\exp (\beta
\depth(\ocA)
%E(\ocA,\cX\setminus\ocA)
%-\beta H(\bottom(\ocA))
)
} \leq {\deg(\alpha)}^{|\cX|}.
\end{eqnarray*}
\end{corollary}
%
%A path $\o\in\tonda{\xi\to\h}$ such that
%$\max_{\z\in\o} H(\z)= E(\xi,\h)$ is called
%{\textbf{optimal}}; the set of optimal paths
%is denoted by $\opt\tonda{\xi\to\h}$.
%The maximizers in the optimal paths in
%$\opt\tonda{\xi\to\h}$ are called {\textbf{saddles}}
%and their set is denoted by $\cS(\xi,\h)$.
%subcycles of
%$\cA_{\minus}$
%$\cC_d$
%such that
%$$H(\widetilde B(\cA_1))=\cdots=
%H(\widetilde B(\cA_n))=\lambda.$$
%%$H(\widetilde B(\cA_i))= \lambda$ for $i=1,\ldots,n$.
Let $\cY$ be a subset of $\cX$.
A cycle $\cA$ (resp., a cycle compound $\ocA$) included in $\cY$
is said
to be maximal if there is no cycle $\cA'$
(resp., no cycle compound $\ocA{}'$) included in $\cY$
such that $\cA\subsetneq\cA'$ (resp., $\ocA\subsetneq\ocA{}'$).

%le2.11 #&#
\begin{lemma}\label{disjoint}
Two maximal cycle compounds in $\cY$ are either equal or disjoint.
\end{lemma}
\begin{pf}
Let $\ocA_1,\ocA_2$ be two
maximal cycle compounds in $\cY$ which are not disjoint.
Suppose that
\[
E(\ocA_1,\cX\setminus\ocA_1) = E(\ocA_2,
\cX\setminus\ocA_2).
\]
Then
$\ocA_1\cup\ocA_2$ is still a cycle compound included in $\cY$.
By maximality, we must have
$\ocA_1=\ocA_2$.
Suppose that
\[
E(\ocA_1,\cX\setminus\ocA_1) < E(\ocA_2,
\cX\setminus\ocA_2).
\]
Let $x$ be a point of
$\ocA_1\cap\ocA_2$. If
$\ocA_1\setminus\ocA_2\neq\varnothing$,
then
\[
E(x,\cX\setminus\ocA_2) \leq \height(\ocA_1) \leq E(
\ocA_1,\cX\setminus\ocA_1),
\]
which is absurd. Thus
$\ocA_1\subset\ocA_2$,
and by maximality,
$\ocA_1=\ocA_2$.
\end{pf}

We denote by $\cM(\cY)$ the partition of $\cY$ into maximal
cycles, that is,
\[
\cM(\cY) = \{ {\cA}\dvtx  {\cA} \mbox{ is a maximal cycle included in $\cY$} \},
\]
and by $\overline{\cM}(\cY)$ the partition of $\cY$ into maximal
cycle compounds, that is,
\[
\overline{\cM}(\cY) = \{ \overline{\cA}\dvtx  \overline{\cA} \mbox{ is a maximal
cycle compound included in $\cY$} \}.
\]
%
%In fact,
%the elements of $\cM(\cB)$ are
%%It is easy to see that, for $\h\in\cB$, the sets
%$$\cA_\h= \sgraffa{\h} \cup\sgraffa{\xi\tc E(\h,\xi)<E(\h,\cB^c)},
%
%A set $\cA$ such that
%$$\max  \big\{ H(\h)\dvtx {\h\in\cA} \big\} <
% \min \big\{ H(\z)\dvtx \z\in\partial\cA \big\}$$
%is called a {{nontrivial}} \textbf{cycle}.
%A singleton consisting of a
%single configuration is called a trivial cycle.
%peut-etre definir "principal boundary"? dans le cas d'un cycle trivial,
%prendre $\{x\in\partial\cA\dvtx  H(x)\leq H(\cA)\}$?
%The principal boundary of a
%nontrivial cycle $\cA$ is
%$$\Bt(\cA)=\bottom(\partial\cA).$$
%For a trivial cycle $\cA$ we set
%$$\Bt(\cA)=\{ x\in\partial\cA\dvtx H(x)\leq H(\cA) \}.$$
%The depth of a cycle $\cA$ is the quantity
%$$\G(\cA)=H(\Bt(\cA))-H(\bottom(\cA)).$$
%We call $H(\Bt(\cA))$ the height of $\cA$.

%le2.12 #&#
\begin{lemma}\label{exitval}
Let
$\ocA$ be a maximal cycle compound included in a subset $\cD$ of $\cX$,
and let $x$ belong to
$\partial\ocA\cap\cD$.
%The energy
Then
$H(x)$ is not equal to
$E(\ocA,\cX\setminus\ocA)$.
If
$H(x)<
E(\ocA,\cX\setminus\ocA)$,
then we have
$E(x,\cX\setminus\cD)<
E(\ocA,\cX\setminus\ocA)$.
%Let
%$\ocA$ be a maximal cycle compound included in a subset $\cD$ of $\cX$
%and let $x$ belong to
%$\partial\ocA\cap\cD$.
%The energy
%$H(x)$ is not equal to $\height(\ocA)$.
%If
%$H(x)<\height(\ocA)$, then we have
%$E(x,\cX\setminus\cD)<\height(\ocA)$.
%$$\forall x\in\partial\ocA\cap\cD
%H(x)\neq\height(\ocA)$$
%and
%$$H(x)<\height(\ocA) \Rightarrow
\end{lemma}
\begin{pf}
If there was a state $x\in\partial\ocA\cap\cD$
such that
$H(x)=E(\ocA,\cX\setminus\ocA)$, then the set
$\ocA\cup\{ x \}$ would be a cycle compound
included in $\cD$, which would be strictly larger than $\ocA$, and this
would contradict the maximality of $\ocA$.
Similarly, for the second assertion, suppose that
$H(x)< E(\ocA,\cX\setminus\ocA)$,
and let
\[
\cA' = \bigl\{ y\in\cX\dvtx E(x,y)< E(\ocA,\cX\setminus\ocA) \bigr\}.
\]
The set $\cA'$ is a cycle of height strictly less than
$E(\ocA,\cX\setminus\ocA)$ and such that
%satisfying
%, in particular
$E(\cA',\cX\setminus\cA')\geq E(\ocA,\cX\setminus\ocA)$.
% Since
%$H(x)<
%E(\ocA,\cX\setminus\ocA)$,
%we see that
Moreover,
\[
\height \bigl(\ocA\cup\cA'\bigr) \leq E(\ocA,\cX\setminus\ocA)
\leq E\bigl(\ocA\cup\cA',\cX\setminus\bigl(\ocA\cup
\cA'\bigr)\bigr).
\]
%
%.$$
%Moreover
%$$E(\ocA',\cX\setminus\ocA') \geq
%E(\ocA,\cX\setminus\ocA),$$
%whence
%$$E(\ocA\cup\ocA',\cX\setminus(\ocA\cup\ocA'))
% \geq
%E(\ocA,\cX\setminus\ocA).$$
Thus $\ocA\cup\cA'$ is still a cycle compound.
%included in $\cD$.
Because of the maximality of $\ocA$, this cycle compound is
not included in $\cD$.
%we can't have
%${\ocA'}\subset\cD$.
Therefore $\cA'$ intersects $\cX\setminus\cD$ and
$E(x,\cX\setminus\cD)< E(\ocA,\cX\setminus\ocA) $.
\end{pf}

%s3 #&#
\section{The stochastic Ising model}
\label{feat}
The material presented in this section is standard and classical.
In Section~\ref{isisec}, we define the Hamiltonian of the Ising model
with various boundary conditions, and we show the benefit of working
with an irrational magnetic field.
In Section~\ref{sigm}, we define the
stochastic Ising model, and we recall the graphical construction, which provides
a coupling between the various dynamics associated to different
boundary conditions and parameters.
%We first recall some basic definitions of the Ising model.
%In this section we introduce the language to describe geometrically
%the energy landscape of the Ising model and to find the relationships
%between the interesting quantities in different dimensions.

%s3.1 #&#
\subsection{The Hamiltonian of the Ising model}
\label{isisec}
%With each configuration $\sigma\in\cX^d= \sgraffa{-1,+1}^{\Z^d}$,
%With each configuration $\h\in\sgraffa{-1,+1}^{\Z^d}$,
%With each configuration $\h\dvtx \Z^d\to\sgraffa{-1,+1}$,
With each configuration $\sigma\in\{-1,+1\}^{\Z^d}$,
we associate a formal Hamiltonian $H$ defined by
\[
H(\sigma) = - \frac{1}{2} \mathop{\sum
_{\{x,y\} \subset\Z^d}}_{|x-y|=1} \sigma(x) \sigma(y) - \frac h 2 \sum
_{x \in\Z^d} \sigma(x).
\]
%
% The elements of $\Z^d$ will be called the lattice sites and
The
value $\sigma(x)$ is the spin at site $x\in\Z^d$ in the
configuration $\sigma$.
Notice that the first sum runs over the unordered pairs $x,y$ of
nearest neighbors
sites of $\Z^d$.
We denote by $\sigma^x$ the configuration obtained from $\sigma$ by flipping
the spin at site $x$.
The variation of energy caused by flipping the spin at site $x$ is
\[
H\bigl(\sigma^x\bigr)-H(\sigma) = \sigma(x) \biggl(\mathop{\sum
_{y \in\Z^d}}_{|x-y|=1} \sigma(y) + h \biggr).
\]
Given a box $\L$ included in $\Z^d$ and a boundary condition
%$\z\dvtx  \Z^d \setminus\L\longrightarrow\sgraffa{-1,+1}$, we define
${\zeta}\in\{ -1,+1 \}^{\Z^d \setminus\L}$, we define
a function $H_{\L}^{{\zeta}}\dvtx  \{ -1,+1 \}^\L\longrightarrow\R$ by
\[
H_{\L}^{{\zeta}} %(\sigma)  =  - \frac1 2 \sum_{ {x,y \in\L}{|x-y|=1}} \sigma(x)
(\sigma) = -
\frac{1}{2} \mathop{\sum_{\{x,y\} \subset\L}}_{|x-y|=1}
\sigma(x) \sigma(y) - \frac h 2 \sum_{x \in\L} \sigma(x)
- \frac{1}{2} \mathop{\sum_{x \in\L,y\notin\L}}_{|x-y|=1}
\sigma(x) {\zeta}(y)+c_\L^{\zeta},
\]
where $c_\L^{\zeta}$ is a constant depending on $\L$ and ${\zeta}$.
%The Gibbs measure on the infinite lattice is defined in a standard way
%by taking the limit of the finite-volume Gibbs measures.
Since $h$ is positive,
for sufficiently large boxes, the configuration
with all pluses, denoted by $\plus$, is the absolute minimum of the
energy for any boundary condition, and it has the maximal
Gibbs probability.
The configuration with all minuses, denoted by $\minus$,
will play the role of the deepest local minimum in our system,
representing the metastable state.
%We will work with the following shifted energy\dvtx
%H_{\L}^{\z}
%(\sigma)  =
%E_{\L}^{\z}
%(\sigma) -
%E_{\L}^{\z}
%(\minus)
%.
We choose the constant
$c_\L^{\zeta}$ so that
\[
H_{\L}^{{\zeta}} (\minus) = 0.
\]
%
%We drop $\z$ from the notation if $\z=\minus$.
%In this section, we
%consider the $d$-dimensional Ising model in a cube
%$Q$ with a sufficiently large side-length and
%minus boundary conditions.
Sometimes we remove $\L$ and ${\zeta}$ from the notation to alleviate
the text,
writing simply $H$ instead of
$H_{\L}^{{\zeta}}$.
The communication kernel $\alpha$ on
$\{ -1,+1 \}^\L$
is defined by
\[
\forall\sigma\in \{ -1,+1 \}^\L, \forall x\in\L\qquad \alpha\bigl(\sigma,
\sigma^x\bigr)=1
\]
and
$\alpha(\sigma,\eta)=0$ if $\sigma$ and $\eta$ have different spins
in two sites or more.
The space
$\{ -1,+1 \}^\L$ is now endowed with a communication
kernel $\alpha$
and an energy
$H_{\L}^{{\zeta}}$,
we define an associated Metropolis dynamics on it
as in Section~\ref{metrosec}.

We shall identify a configuration of spins
with the support of the pluses in it; this way,
we think of a configuration as a set, and we can perform
the usual set operations on configurations.
For instance,
we denote by $\h\cup\xi$ the configuration in which
the set of pluses is the union of the sets of pluses in
$\h$ and in $\xi$.
We call volume
of a configuration $\h$
the number of pluses in $\h$
and we denote it by $|\h|$.\vadjust{\goodbreak}
We call perimeter
of a configuration $\h$
the number of the interfaces between the pluses and the minuses in $\h$
%the cardinality of the set of nearest neighbors of $\h$
and we denote it by $p(\h)$,
\[
p(\eta) = \bigl| \bigl\{ \{x,y\}\dvtx \eta(x)=+1, \eta(y)=-1, |x-y|=1 \bigr\} \bigr|.
\]
The Hamiltonian of the Ising model can then be rewritten
conveniently as
\[
H(\h) = p(\h)-h|\h|.
\]
Our analysis of the energy landscape will be based on the assumption
that $h$ is an irrational number.
%$h \in\R\setminus\Q$, that will be justified in Remark \ref{r2}.
This hypothesis simplifies in a radical way our study because of the
following lemma.

%le3.1 #&#
\begin{lemma}\label{irrazionale}
%If $h \in\R\setminus\Q$, then a
%Any configuration included in
%$\h$ with the same energy coincides with $\h$.
Let $h$ be an irrational number.
Suppose $\sigma,\eta$ are two configurations such that
$\sigma\subset\eta$ and $H(\sigma)=H(\eta)$. Then $\sigma=\eta$.
\end{lemma}
\begin{pf}
Since $h$ is irrational, the knowledge of the energy
of a configuration determines in a unique way its perimeter
and its volume.
Since $\sigma$ is included in $\eta$ and they have the same volume,
%then $\sigma$ and $\eta$ are equal.
then they are equal.
\end{pf}

In the next section, we build a monotone coupling of the dynamics associated
to different magnetic fields $h$. With the help of this coupling,
we will show in Section~\ref{rer}
that it is sufficient to prove Theorem \ref{mainfv} for irrational
values of the magnetic field. The main point is that the critical
constant $\kappa_d$ depends continuously on $h$ (this is proved in
Lemma \ref{contga}).
% From this coupling and the continuity of
%the critical constants as function

We believe that the main features of the cycle structure should persist
for rational values of $h$.
%, but we do not see the need to investigate
%in this direction.
The assumption that $h$ is irrational
(or at least that it does not belong to some countable set)
is present in most papers to simplify the structure of
the energy landscape, with the only exception of
\cite{MNOS}.
%, where this analysis is made.
In dimension 2, for $2/h$ integer,
there exists a very complicated
cycle compound, consisting of cycles with the same depth
that communicate at the same energy level; see \cite{MNOS}.
This compound is not contained in the metastable cycle
and is compatible with our results.

Our analysis is based on the following attractive inequality.

%le3.2 #&#
\begin{lemma}\label{attrineq}
For any configurations $\h$, $\xi$, we have
%$$\forall\h,\xi
%
\[
H(\h\cap\xi)+H(\h\cup\xi) \le H(\h) + H(\xi).
\]
\end{lemma}
\begin{pf}
This inequality
can be proved with a direct computation; see
Theorem~5.1 of \cite{BAC}.
\end{pf}
%

%s3.2 #&#
\subsection{Graphical construction}
\label{sigm}\label{graphi}
The time evolution of the model
is given by the Metropolis dynamics:
when the system is in the configuration $\h$,
the spin at a site $x \in\L\subset\Z^d$ flips at rate
\[
c_{\L,\b}^{{\zeta}}(x,\h) = \exp \bigl( -\b\max \bigl(0,
H_{\L}^{{\zeta}} \bigl(\h^x\bigr)-
H_{\L}^{{\zeta}} (\h) \bigr) \bigr),\vadjust{\goodbreak}
\]
where
%$[a]_+=\max\sgraffa{a,0}$ is the positive part of $a$ and
the parameter $\b$ is the inverse temperature.
A standard construction yields a continuous-time Markov process
%(reversible with respect to the Gibbs measure)
whose generator is defined by
\[
\forall f\dvtx \{-1,+1\}^{\L} \to\R\qquad (Lf) (\h) = \sum
_{x \in\L} c_{\L,\b}^{{\zeta}} (x,\h)
\bigl(f\bigl(\h^x\bigr)-f(\h)\bigr).
\]
The process in a $d$-dimensional
box $\Lambda$, under magnetic field $h$, with initial condition $\a$
and boundary condition ${\zeta}$ is denoted by
\[
\bigl(\s^{\a,{\zeta}}_{\Lambda,t}, t\geq0\bigr).
\]
To define the process in infinite volume, we consider the weak limit of the
previous process as $\L$ grows to $\Z^d$. This weak limit does not depend
on
the sequence of the boundary conditions; see \cite{Sch} for the details.
Sometimes
we omit $\Lambda$, $\a$ or ${\zeta}$ from
the notation if $\Lambda=\Z^d$, $\a=\minus$,
or ${\zeta}=\minus$, respectively.
%whenever the dimension or the magnetic field are clear from the
%context,
%we also drop $d$ and $h$ from the notation.

In order to compare different processes,
we use a standard construction,
known as the graphical construction,
that allows us to define on the same probability space all the
processes at a given inverse temperature
$\b$, in $\Z^d$ and in any of its finite subsets,
with any initial and boundary conditions and any magnetic field $h$.
We refer to \cite{Sch} for details.
We consider two families of i.i.d.
Poisson processes with rate one,
associated with the sites in $\Z^d$.
Let $x \in\Z^d$. We denote by
$(\t^-_{x,n})_{n\geq1}$
and by
$(\t^+_{x,n})_{n\geq1}$
the arrival times of the two Poisson processes associated to $x$.
Notice that, almost surely, these random times are all distinct.
With each of these arrival times, we associate uniform random
variables
$(u^-_{x,n})_{n\geq1}$,
$(u^+_{x,n})_{n\geq1}$,
and we assume that
these variables are independent of
each other and of the Poisson processes.
We introduce next
an updating procedure in order to define simultaneously
all the processes
on this probability space.
Let $\Lambda$ be a finite subset of $\Z^d$, and let $x\in\Lambda$.
Let
$\varepsilon=-1$ or $\varepsilon=+1$, let $\a$ be an initial
configuration and
let ${\zeta}$ be a boundary condition.
Let $\s$ denote the configuration just before time $\t^\varepsilon_{x,n}$.
The updating rule at time
$\t^\varepsilon_{x,n}$ is the following:

$\bullet$ the spins not at $x$ do not change;

$\bullet$ if
$\s(x)=-\varepsilon$
and
$u^\varepsilon_{x,n}<c_{\L,\b}^{{\zeta}}(x,\s)$, then
the spin at $x$ is reversed.
% the spin at $x$ is turned upside down.

If the set $\L$ is finite, then the above rules define
a Markov process $(\sigma^{\a,{\zeta}}_{\L,t})_{t\geq0}$.
Whenever $\L$ is infinite, one has to be more careful, because
there is an infinite number of arrival times in any finite
time interval, and
it is not possible to order them in an increasing sequence.
However, because the rates are bounded, changes in the system
propagate at a finite speed, and
a Markov process
%$(\sigma_{\L,t})_{t\geq0}$
can still
be defined by taking the limit of finite volume processes; see \cite
{Sch,L} for more details.
In any case
our proofs will
involve mainly boxes whose side length is finite,
although they might grow with $\b$.
From now on, we denote by $P$ and $E$ the probability and
expectation with respect to the family of the Poisson processes
and the uniform random variables.
The graphical construction allows us to take advantage of
the monotonicity properties of the rates
$c_{\L,\b}^{{\zeta}}(x,\s)$.
For
any box $\L$,
any configurations
%$\a\subset\a'$, $\z\subset\z'$,
$\a\leq\a'$, ${\zeta}\leq{\zeta}'$,
%any boxes $\L\subset\L' \subset\Z^d$ and any fields
%$h \le h'$,
we have %for all $t\ge0$
\[
\forall t\geq0\qquad \s^{\alpha,{\zeta}}_{\L,t} \leq
\s^{\alpha',{\zeta}'}_{\L,t}.
\]
The process is also nondecreasing as a function of the magnetic field $h$.

%s3.3 #&#
\subsection{Reduction to irrational fields}
\label{rer}

We show here how the monotonicity of the process as a function of the
magnetic field, together with the continuity of $\Ga_d$ and~$\kappa_d$,
allow us to reduce the study to irrational values of the magnetic field.
Suppose that Theorem \ref{mainfv} has been proved for irrational
values of the magnetic field.
Let $h<h_0$ be a positive rational number,
and let $\kappa<\max(\Gamma_d-dL,\kappa_d)$.
As we will see in Lemma \ref{contga}, the constants $\Ga_d$ and
$\kappa_d$
depend continuously on $h$, and therefore there
exists
an irrational number
$h'$
such that $h<h'<h_0$ and
\[
\kappa<\max\bigl(\Gamma_d'-dL,\kappa_d'
\bigr),
\]
where $\Ga_d'$ and $\kappa_d'$ are the constants associated to the
field $h'$.
Theorem \ref{mainfv} applied to the process
$(\s_{\L_\beta,t}^{-,\minus,h'})_{t\geq0}$ associated
to the field $h'$ yields
\[
\lim_{\beta\to\infty} \P \bigl(\s_{\L_\beta,\tau_\beta}^{-,\minus,h'}(0)=1
\bigr) = 0.
\]
From the graphical construction, we have
\[
{ \s_{\L_\beta,\tau_\beta}^{-,\minus,h}(0) \leq
\s_{\L_\beta,\tau_\beta}^{-,\minus,h'}(0)},
\]
whence
\[
\lim_{\beta\to\infty} \P \bigl(\s_{\L_\beta,\tau_\beta}^{-,\minus,h}(0)=1
\bigr) = 0
\]
as desired.
The second part of Theorem \ref{mainfv} for
rational values of $h$ is proved similarly.
Therefore,
it is sufficient to prove Theorem \ref{main}
for $h$ irrational.
For the remainder of the paper, we will assume that it is the case.
This will allow us to use the result of Lemma \ref{irrazionale} which
implies the other results on the energy landscape proven
in Section~\ref{stc},
in particular Lemma \ref{fondo}.

%s4 #&#
\section{Isoperimetric results}
\label{energyestmates}
In this section we report some specific results on the
energy landscape
of the $d$-dimensional Ising model.
In the two-dimensional case, a very detailed description
can be found in \cite{NS2,NS1}.
In three dimensions, the cycle structure is known only
near the typical transition paths; see \mbox{\cite{N2,N1,AC,BAC}}.
In higher dimensions, we can compute the communication energy
between $\minus$ and $\plus$ by using the results of Neves \cite{N1},
but finer details are still unknown.
In Section~\ref{sie},
we state a discrete isoperimetric inequality which will be used
in the proof of Lemma \ref{control}.
In Section~\ref{sir}, we define the so-called reference path.
Thanks to the isoperimetric results of Neves, we can compute the critical
energy $\Ga_d$ with the help of the reference path.
This is done in Section~\ref{sitm}. As a by-product, we prove that the
energy $\Ga_d$ depends continuously on $h$.
In the inductive proof of Theorem \ref{T2}, we work with mixed boundary
conditions, called
$n\pm$ boundary conditions.
In Section~\ref{six}, we define
the $n\pm$ boundary conditions, and we prove the
required isoperimetric results in boxes with these boundary conditions.
%are the counterpart of those presented in Section \ref{}
%and they

%s4.1 #&#
\subsection{An isoperimetric inequality}
\label{sie}
A $d$-dimensional polyomino is a set which is the finite union of unit
$d$-dimensional cubes.
There is a natural correspondence between configurations and polyominoes.
To a configuration
we associate the polyomino which is the union of the unit cubes
centered at the
sites having a positive spin. The main difference between
configurations and
polyominoes
is that the polyominoes are defined up to translations.
Neves \cite{N1} has obtained a discrete isoperimetric inequality
in dimension $d$, which yields the exact value of
\[
\min \bigl\{ \per(c)\dvtx  c\mbox{ is a $d$-dimensional polyomino of volume $v$}
\bigr\},
\]
where $v\in\N$.
This value is a quite complicated function of the volume $v$, which is
larger than
\[
2d \bigl\lfloor v^{1/d} \bigr\rfloor^{d-1}.
\]
We derive from this the following simplified
isoperimetric inequality.\vspace*{9pt}

\textit{Simplified isoperimetric inequality}.
For a
$d$-dimensional polyomino $c$,
\[
\per(c) \geq 2d \bigl(\mbox{volume}(c) \bigr)^{(d-1)/d}.
\]
\begin{pf}%{Proof of ???}
We rely on the inequality stated above, and we
perform a simple scaling with an integer factor $N$,
\begin{eqnarray*}
&&
\min \bigl\{ \per(c)\dvtx  c\mbox{ $d$-dimensional polyomino of volume $v$} \bigr\}
\\
&&\qquad
\geq \min \bigl\{ \per \bigl( N^{-1/d}c \bigr)\dvtx  c\mbox{ polyomino of
volume $Nv$} \bigr\}
\\
&&\qquad
= N^{({1-d})/{d}} \min \bigl\{ \per( c)\dvtx  c\mbox{ polyomino of volume $Nv$}
\bigr\}
\\
&&\qquad
\geq N^{({1-d})/{d}} 2d \bigl\lfloor(Nv)^{1/d} \bigr
\rfloor^{d-1}.
\end{eqnarray*}
Sending $N$ to $\infty$, we obtain the desired inequality.
\end{pf}

If we had applied the classical isoperimetric inequality in
$\mathbb{R}^d$, then we would have obtained an inequality with
a different constant, namely the perimeter of the
unit ball instead of $2d$. The constant $2d$ is sharp, indeed
there is equality
when
$c$ is a $d$-dimensional cube whose side length is an integer.
We believe that, for polyominoes of volume
equal to $l^d$ where $l$
is an integer,
it is
the only shape realizing the equality,
%is the $d$--dimensional
%cube of side length equal to $l$,
yet we were unable to locate a proof
of this statement in the literature (apart for the three-dimensional case
\cite{BAC}).
We will need the simplified isoperimetric inequality with the correct
constant in the main inductive proof.

%s4.2 #&#
\subsection{The reference path}
\label{sir}
Let $R$ be
a parallelepiped in $\Z^d$ whose vertices belong to
$\Z^d+(1/2,\ldots,1/2)$ and whose sides are parallel to the axis.
A face of $R$ consists of the set of the sites of $\Z^d$
%outside the parallelepiped
which are at distance $1/2$ from the parallelepiped
and which are contained in a given single hyperplane.
%We will often use induction on the dimension.
With a slight abuse of terminology,
we say that a
configuration $\eta$ is obtained by attaching a $(d-1)$-dimensional
configuration $\xi$
to a face of a $d$-dimensional parallelepiped ${\zeta}$
if $\h={\zeta}\cup\xi$ and $\xi$ is contained in a face of ${\zeta}$.
It is immediate to see that in this case
\[
H_{\Z^d}({\zeta}\cup\xi) = H_{\Z^d} ({\zeta}) +
H_{\Z^{d-1}} (\xi).
\]
%
%This identity will be the key for our inductive proof.
We call quasicube a %$d$-dimensional
parallelepiped in $\Z^d$
such that the shortest and the longest side lengths differ
at most by one length unit.
Notice that the faces of a quasicube are
$(d-1)$-dimensional quasicubes.
%Quasi-cubes can be ordered according to their volume.
%A quasi-cube with maximal volume among the quasi-cubes with
%volume smaller than or equal to $i$ is
%said maximal up to volume $\mathbf i$.
From the results of Neves \cite{N1} we see that
there exists an optimal path
% from the energy point of view, the best path
from $\minus$
to $\plus$ made of configurations
which are as close as possible to a cube.
% The following path will be a useful tool in our analysis.
%Let $\L$ be a box whose sides are larger than $2d/h$.
We call reference path in a box $\Lambda$
a path
$\rho=(\rho_0,\ldots,\rho_{|\L|})$
%$\rho\in(\minus\to\plus)$
going from $\minus$ to $\plus$
built with the following algorithm.
In one dimension, $\rho_i$ has exactly $i$ pluses
which form an interval of length $i$.
%In dimension $d$,
%$\rho_i$ is a maximal quasicube up to volume $i$
%with a configuration of a $(d-1)$-dimensional reference path
%attached on one of its largest faces.
%An algorithm to construct a reference path is the following:
%%\bd{refpath}
%In one dimension, the path flips consecutively the spins
%in the box;
In higher dimension, we proceed as follows:
\begin{longlist}[(3)]
\item[(1)] Put a plus somewhere in the box.
\item[(2)] Fill one of the largest faces of the parallelepiped of pluses
(among that contained in the box),
following a $(d-1)$-dimensional reference path.
\item[(3)] Go to step 2 until the entire box is full of pluses.
\end{longlist}
With a reference path $\rho=(\rho_0,\ldots,\rho_{|\L|})$, we
associate a
re\-fe\-ren\-ce cy\-cle path
consisting of the
sequence of cycles
$(\pi_0,\ldots,\pi_{|\L|})$, where for $i=0,\ldots,|\L|$,
the cycle
$\pi_i$ is the maximal cycle of
$\{ -1,+1 \}^\Lambda\setminus\{ \minus,\plus \}$
%$\cC_d\setminus\{ \minus \}$
containing $\rho_i$.
%$$\cC_i=\sgraffa{\o_i}\cup
A~reference path enjoys the following remarkable property:
\[
\forall i<j\qquad E(\rho_i,\rho_j) = \max \bigl\{ H(
\rho_k)\dvtx i\leq k\leq j \bigr\},
\]
that is, it realizes the solution of the minimax problem associated with
the communication energy between any two of its configurations.
%s4.3 #&#
\subsection{The metastable cycle}
\label{sitm}
Let $\L$ be a box whose sides are larger than $2d/h$.
We endow $\L$ with minus boundary conditions.
The metastable cycle $\cC_d$ in the box $\Lambda$ is the
maximal cycle of
\[
\{ -1,+1 \}^\Lambda\setminus\{ \plus \}
\]
containing $\minus$ in the energy landscape associated to $H_\L^-$,
the Hamiltonian in $\L$ with minus boundary conditions. We define
\[
\G_d = \depth(\cC_d) = E(\minus,\plus).
\]
Recall that, by convention, $H(\minus)=0$.
Obviously, a path $\o=(\o_0,\ldots,\o_l)$
going from $\minus$ to $\plus$ satisfies
\begin{eqnarray*}
\max_{0\leq i\leq l}H(\o_l) &\geq& \max
_{0\leq k\leq|\Lambda|} \min \bigl\{ H(\sigma)\dvtx  \sigma\in\{ -1,+1
\}^\Lambda,|\sigma|=k \bigr\}
\\
&=& \max_{0\leq k\leq|\Lambda|} \bigl( \min \bigl\{ p(\sigma)\dvtx  \sigma\in\{ -1,+1
\}^\Lambda, |\sigma|=k \bigr\}-hk \bigr),
\end{eqnarray*}
and the reference path $\rho$ realizes the equality in this inequality.
We conclude therefore that
\[
\G_d = \max_{0\leq k\leq|\Lambda|} H(\rho_k).
\]
When $h$ is irrational,
% and smaller than $1$,
%we find from a direct computation that
there exists a unique value $m_d$ such that
\[
\G_d = H(\rho_{m_d}),
\]
that is, the value $\G_d$ is reached for the configuration of a reference
path having volume $m_d$. We call such a configuration a critical droplet.

%f1
\begin{figure}

\includegraphics{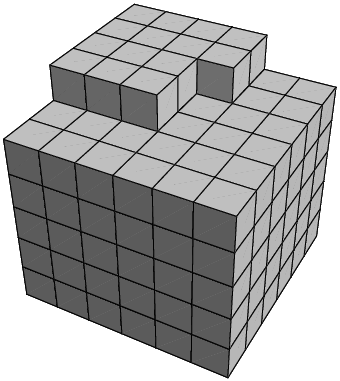}

\caption{A configuration of the reference path.}
\end{figure}

From the results of Neves \cite{N2,N1}
and a direct computation, we derive the following facts.
Let
\[
l_c(d) = \biggl\lfloor\frac{2(d-1)}{h} \biggr\rfloor.
\]
The configuration of volume $m_d$ is a quasicube with sides of length
$l_c(d)$ or \mbox{$l_c(d)+1$}, with a $(d-1)$-dimensional critical droplet
attached on one
of its largest sides. The precise shape of the critical droplet depends
on the value of $h$ (see, e.g., \cite{BAC} for $d=3$);
by the precise shape, we mean the number of sides of the quasicube
which are equal to
$l_c(d)$ and $l_c(d)+1$.
It is possible to derive exact formulas for $m_d$ and $\G_d$, but
they are complicated, and it is necessary to consider various cases
according to the value of $h$.
However, we have $m_1=1$, $\G_1=2-h$ and the
following inequalities
\begin{eqnarray*}
\bigl(l_c(d) \bigr)^d &\leq& m_d \leq
\bigl(l_c(d)+1 \bigr)^d,
\\
2d \bigl(l_c(d) \bigr)^{d-1} -h \bigl(l_c(d)+1
\bigr)^d &\leq& \Ga_d \leq 2d \bigl(l_c(d)+1
\bigr)^{d-1} -h \bigl(l_c(d) \bigr)^d.
\end{eqnarray*}
This yields the following expansions as $h$ goes to $0$:
\[
m_d \sim \biggl( \frac{2 (d-1)}{h} \biggr)^d,\qquad
\G_d \sim2 \biggl( \frac{2 (d-1)}{h} \biggr)^{d-1}.
\]
\begin{lemma}\label{contga}
The energy $\Ga_d$ of the critical droplet in dimension $d$ is a continuous
function of the magnetic field $h$.
\end{lemma}

\begin{pf}
Let $h_0>0$. Let $\L$ be a box of side length larger than $4d/h_0$.
%Let $h\geq h_0$.
From the previous results,
for any $h\geq h_0$, we have the equality
\[
\Ga_d %\max_{0\leq i\leq l}H(\o_l)
= \max_{0\leq k\leq|\Lambda|} \min \bigl\{ H(
\sigma)\dvtx  \sigma\in\{ -1,+1 \}^\Lambda, |\sigma|=k \bigr\}. % =
\]
Given a configuration $\sigma$ of spins in $\Lambda$,
the Hamiltonian $H(\sigma)$ is a continuous function of the
magnetic field $h$.
For $k\leq|\Lambda|$, the number of configurations $\sigma$
such that
$|\sigma|=k$ is finite, thus
the minimum
\[
\min \bigl\{ H(\sigma)\dvtx  \sigma\in\{ -1,+1 \}^\Lambda, |\sigma|=k \bigr\}
\]
is also a continuous function of $h$.
Thus $\Ga_d$ is also a continuous function of $h$ on $[h_0,+\infty[$.
This holds for any $h_0>0$, thus
$\Ga_d$ is a continuous function of $h$ on $]0,+\infty[$.
\end{pf}

Our next goal is to prove that
the maximal depth of the cycles in a reference cycle path is
%not larger
smaller than $\G_{d-1}$.
Let $\rho=(\rho_0,\ldots,\rho_{|\L|})$ be a reference path,
and
let $(\pi_0,\ldots,\pi_{|\L|})$ be the corresponding reference cycle path.
We set
\[
\Delta_d = \max_{0\leq i<m_d}\depth(\pi_i) =
\max_{0\leq i<m_d} \bigl(E(\pi_i,\minus)-E\bigl(\bottom(
\pi_i)\bigr) \bigr).
\]

%pr4.2 #&#
\begin{proposition}\label{dee}
The maximal depth $\Delta_d$ of the cycles in a reference cycle path
is strictly less than $\G_{d-1}$.
\end{proposition}
\begin{pf}
For $i<m_d$ the configuration
$\rho_i$ belongs to $\cC_d$, and we have
\[
E(\pi_i,\minus) = \max_{0\leq j\leq i}H(
\rho_j).
\]
Let us define, for $0\leq i\leq r$,
\begin{eqnarray*}
\uv_i &=& \min \bigl\{ |\sigma|\dvtx \sigma\in\pi_i \bigr\},
\\
\ov_i &=& \max \bigl\{ |\sigma|\dvtx \sigma\in\pi_i \bigr\}.
\end{eqnarray*}
Whenever $i<m_d$, the value $\uv_i$ is the unique integer $v$
such that
\[
H(\rho_{v-1}) = E(\pi_i,\minus).
\]
Thanks to the minimax property of the reference path, we have also
that
$\rho_k\in\pi_i$ for $\uv_i\leq k\leq\ov_i$ whence
\[
E\bigl(\bottom(\pi_i)\bigr) = \min \bigl\{ H(\rho_k)\dvtx
\uv_i\leq k\leq\ov_i \bigr\}.
\]
From the previous identities, we infer that
\begin{eqnarray*}
\Delta_d &=& \max_{0\leq i<m_d} \max \bigl\{ H(
\rho_{\uv_i-1})- H(\rho_{k})\dvtx  \uv_i\leq k\leq
\ov_i \bigr\}
\\
&\leq& \max_{0\leq j\leq i<m_d} \bigl(H(\rho_{j})- H(
\rho_{i}) \bigr).
\end{eqnarray*}
The maximum of the energy
along a $(d-1)$-dimensional reference path
is reached at the value $m_{d-1}$, while the minimum
of the energy is reached at one of the two ends of the path.
Therefore the indices $i^*,j^*$ realizing the maximum of the
right-hand side
correspond, respectively, to a quasicube
$\rho_{i^*}$
and the union
$\rho_{j^*}$
of a quasicube $c^*$
and a $(d-1)$-dimensional critical droplet.
Since
%we must have also
$j^*\leq i^*$, we have
\mbox{$c^*\subset\rho_{j^*}\subset
\rho_{i^*}$}. The quasicubes
$c^*$ and $\rho_{i^*}$ being
subcritical, we have
$H(c^*)<H(\rho_{i^*})$ and therefore
\[
\Delta_d \leq H(\rho_{j^*})- H(\rho_{i^*}) < H(
\rho_{j^*})- H\bigl({c^*}\bigr) \leq \G_{d-1}.
\]
The last inequality holds also when
$c^*$ is too small so that a $(d-1)$-dimensional
critical droplet cannot be attached to one of its faces.
\end{pf}
%
%s4.4 #&#
\subsection{\texorpdfstring{Boxes with $n\pm$ boundary conditions}
{Boxes with n+- boundary conditions}}
\label{six}

Unlike in the simplified model studied in \cite{CM2}, we cannot use here
a direct induction on the dimension $d$. Instead, we introduce special
boundary conditions that make a $d$-dimensional system behave
like a $n$-dimensional system.
For $E$ a subset of $\Z^d$, we define its outer vertex boundary
$\partial^{\mathrm{out}}E$ as
\[
\partial^{\mathrm{out}}E = \bigl\{ x\in\Z^d\setminus E\dvtx \exists y\in
E, |y-x|=1 \bigr\}.
\]
Let $n\in\{ 0,\ldots,d \}$. We define next
mixed boundary conditions for parallelepipeds
with minus on $2n$ faces and plus on
$2d-2n$ faces.\vspace*{9pt}

\textit{$n\pm$ Boundary condition}. Let $R$ be a parallelepiped. We
write $R$ as the product $R=\Lambda_1\times\Lambda_2$, where
$\Lambda_1,\Lambda_2$ are parallelepipeds of dimensions $n,d-n$,
respectively. We consider the boundary conditions on $R$ defined as:

$\bullet$ minus on
$ (\partial^{\mathrm{out}}
\L_1 ) \times\L_2$;

$\bullet$ plus on
$\L_1\times
\partial^{\mathrm{out}}
\L_2$.

We denote by $n\pm$ this boundary condition, and by $H^{n\pm}$
the corresponding Hamiltonian in $R$.
%More generally, if $R$ is a $d$--dimensional parallelepiped, the
The $n\pm$ boundary condition on $R$ is obtained by
putting minuses on the exterior
faces of $R$ orthogonal to the first $n$ axis and pluses on the
remaining faces.

%f2 #&#
\begin{figure}

\includegraphics{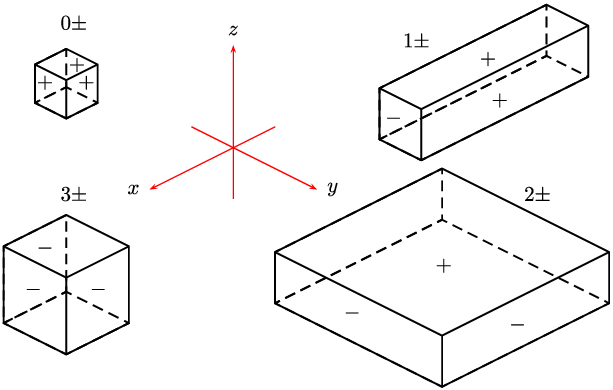}

\caption{$n\pm$ boundary conditions.}
\end{figure}

We will now transfer the isoperimetric results in $\Z^d$
to parallelepipeds with $n\pm$
boundary condition.
%We denote by $(\s^{=}_{\Sigma;t})_{t\geq0}$ the process in $\Sigma$
%evolving with the sandwich boundary conditions $\rho$.

%le4.3 #&#
\begin{lemma}\label{prok}
Let $n\in\{ 1,\ldots,d \}$. Let $R$
be a $d$-dimensional parallelepiped, and let $l$ be the length
of its smallest side.
For any configuration $\sigma$
in $R$ such that $|\sigma|<l$,
there exists an $n$-dimensional configuration $\rho$
such that
\[
|\rho|=|\s|,\qquad H_{\Z^n}(\rho) \leq H^{n\pm}_R(
\sigma).
\]
\end{lemma}
\begin{pf}
The constraint on the cardinality of $\sigma$ ensures that there
is no cluster of pluses connecting two opposite faces of $R$.
We endow $\N^d$ with $n\pm$ boundary conditions by putting
minuses on
\[
\bigl(\{ -1 \}\times\N^{d-1} \bigr)\cup\cdots\cup \bigl(
\N^{n-1}\times\{ -1 \}\times\N^{d-n} \bigr)
\]
and pluses on
\[
\bigl(\N^n\times\{ -1 \}\times\N^{d-n-1} \bigr)\cup\cdots\cup
\bigl(\N^{d-1}\times\{ -1 \} \bigr).
\]
We shall prove the following assertion, which implies the claim of
the lemma. Suppose $n<d$.
For any finite configuration $\sigma$ in $\N^d$, there exists
a configuration $\rho$ in $\N^{d-1}$ such that
\[
|\rho|=|\s|,\qquad H_{\N^{d-1}}^{n\pm}(\rho) \leq H_{\N^d}^{n\pm}(
\sigma).
\]
If we start with a
configuration $\sigma$
in $R$ such that $|\sigma|<l$, then we
apply iteratively this result to the connected components of $\sigma$
(since no connected component of $\sigma$ intersects two opposite faces
of $R$, up to a rotation,
their energies can be computed as if they were in $\N^d$
with $n\pm$ boundary conditions). We end up with a configuration $\eta$
in $\N^n$ with $n\pm$ boundary conditions which satisfies the conclusion
of the lemma. We prove next the assertion.
Let
$\sigma$ be a finite configuration in $\N^d$, and let
$c$ be the polyomino associated to $\sigma$.
We let $c$ fall by gravity along the $(n+1)$th axis
on $\N^{n}\times\{ -1 \}\times\N^{d-n-1}$.\looseness=1

%f3
\begin{figure}

\includegraphics{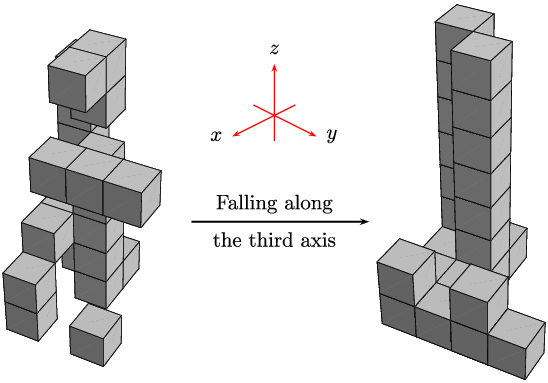}

\caption{The falling operation.}
\end{figure}

The resulting polyomino $\tic$ has the same volume as $c$
and moreover
\[
\per(\tic) \leq \per(c),
\]
because the number of contacts between the unit cubes or
with the boundary condition cannot decrease
through the ``falling'' operation.
We can think of $\tic$ as a stack of $(d-1)$-dimensional polyominoes
$c_0,\ldots,c_k$, which are obtained by intersecting $\tic$ with the layers
\[
L_i = \bigl\{ x=(x_1,\ldots,x_d)\in
\N^d\dvtx  i-\tfrac{1}{2}\leq x_{n+1}<i+\tfrac{1}{2}
\bigr\},\qquad i\in\N.
\]
Since we have let $c$ fall by gravity to obtain $\tic$,
this stack is nonincreasing in the following sense:
for $i$ in $\N$,
the $(d-1)$-dimensional polyomino $c_i$
associated with the layer $L_i$ contains
the $(d-1)$-dimensional polyomino $c_{i+1}$
associated with the layer $L_{i+1}$.
As a consequence,
\[
H_{\N^d}^{n\pm}(\tic) \geq \sum_{i\geq0}
H_{\N^{d-1}}^{n\pm}(c_i) + \operatorname{area}\bigl(
\mathrm{proj}_{n+1}(\tic)\bigr),
\]
where\vspace*{1pt} $\mathrm{proj}_{n+1}(\tic)$ is the orthogonal
projection of $\tic$ on $\N^{n}\times\{ -1 \}\times\N^{d-n-1} $. Let
$\hic$ be a $(d-1)$-dimensional polyomino obtained as the union of
disjoint translates of $c_0,\ldots,c_k$. The polyomino $\hic$ answers
the problem.
\end{pf}

Let $\L$ be a box whose sides are larger than $m_n$.
We construct next a
reference path
$(\rho^{n\pm}_i,0\leq i\leq|\L|)$
in the box $\L$
endowed with
$n\pm$ boundary conditions with the following algorithm:
\begin{longlist}[(4)]
\item[(1)] Compute the maximum number $m$ of plus neighbors for a
minus site in the box (taking into
account the boundary conditions).
\item[(2)] If there is only one site realizing this maximum, put a
plus at this site and
go to step 1.
\item[(3)] Otherwise, compute the maximal length of a segment of minus
sites having all
$m$ plus neighbors.
\item[(4)] Put a plus at a site of a segment realizing the previous maximum
and go to step 1.
%at a site having $d-n$ plus neighbors.
%have $d-n+1$ plus neighbors and select one such face among
%the largest.
%following a $(d-1)$-dimensional reference path
%with $n+1\pm$ boundary conditions.
%$d-n+1$ plus neighbors.
\end{longlist}
As before, the reference path
$(\rho^{n\pm}_i,0\leq i\leq|\L|)$
realizes the solution of the minimax problem associated with
the communication energy between any two of its configurations.
The metastable cycle $\cC_d^{n\pm}$ in the box $\Lambda$ with
$n\pm$ boundary conditions is the
maximal cycle of
\[
\{ -1,+1 \}^\Lambda\setminus\{ \minus,\plus \}
\]
containing $\minus$ in the energy landscape associated to the
Hamiltonian $H^{n\pm}_{\L}$.

%co4.4 #&#
\begin{corollary}\label{dep}
The depth of
the metastable cycle $\cC_d^{n\pm}$ is equal to
$\G_n$.
%$$\G_n = \depth(\cC_d^{n\pm}).$$
\end{corollary}
\begin{pf}
With the help of Lemma \ref{prok}, we can compare the energy along
a path in $\L$ with $n\pm$ boundary conditions with the energy along a
path in $\Z^n$, in such a way that
at each index the configurations in each
path have the same cardinality. This construction implies immediately
that
\[
\depth\bigl(\cC_d^{n\pm}\bigr) \geq \G_n.
\]
To get the converse inequality we simply consider the reference
path
in $\L$ with $n\pm$ boundary conditions.
\end{pf}
%
%co4.5 #&#
\begin{corollary}\label{codep}
The maximal depth $\Delta_d^{n\pm}$ of the cycles in a reference
cycle path
with $n\pm$ boundary conditions
is strictly less than $\G_{n-1}$.
\end{corollary}
\begin{pf}
%By Lemma \ref{prok},
We check that, until the index $m_n$, the energy along the
reference path $(\rho^{n\pm}_i,i\geq0)$ is equal to the
energy along the reference path in $\Z^n$ computed with~$H_{\Z^n}$.
The result follows then from Proposition \ref{dee}.
\end{pf}

%s5 #&#
\section{The space--time clusters}
\label{stc}
%$\t$, $\P_\h$

The goal of this section is to prove Theorem \ref{totcontrole}, which
provides a control on the diameter of the space--time clusters. Theorem
\ref{totcontrole} is used in an essential way in the proof of the lower
bound on the relaxation time, under the following weaker form: for the
dynamics restricted to a small box, the probability of creating a large
$\STC$ before nucleation is SES. We first recall the basic definitions
and properties of the space--time clusters in Section~\ref{badp}. We
next proceed to show that it is very unlikely that large space--time
clusters are formed before nucleation. The main theorem of this
section, Theorem \ref{totcontrole}, is the analog of Lemma 4 in
\cite{DS1}. The proof in \cite{DS1} relies on the fact that in two
dimensions the energy needed to grow, that is, the energy of a
protuberance, is larger than the energy needed to shrink a subcritical
droplet. In higher dimension, we are not able to prove a corresponding
result.
% and the counterparts of the protuberances have a geometry.
Let us give a quick sketch of the proof of
Theorem \ref{totcontrole}.
We consider a set $\cD$ satisfying a technical hypothesis, and we want
to control
the probability of creating a large space--time cluster before exiting
$\cD$.
Typically, the set $\cD$ is a cycle or a cycle compound included in the
metastable cycle.
We use several ideas coming from the theory of simulated annealing
\cite{CaCe}. We decompose $\cD$ into its maximal cycle compounds, and
we show that, before exiting $\cD$, the process is unlikely to make
a large number of jumps between these maximal cycle compounds. Thus, if
a large
space--time cluster is created, then it must be created during a visit to
a maximal cycle compound.
The problem is therefore reduced to control the size
of the space--time cluster created
inside a cycle compound $\ocA$ included in $\cD$.
%More precisely,
%we need to prove the statement of Theorem \ref{totcontrole}
%for a cycle compound.
The key estimate is proved by induction over the depth of the cycle compound.
Suppose we want to prove the estimate for a
cycle compound $\ocA$.
A first key fact, proved
%in Lemma 5.4
in Lemma \ref{fondo}
with the help of the ferromagnetic inequality,
is that in the Ising model under irrational magnetic field
the bottom of every cycle compound is a singleton.
Let $\eta$ be the bottom of
the cycle compound $\ocA$.
We consider now the trajectory of the process starting from a point
of $\ocA$ until it exits from $\ocA$.
In Section~\ref{sectriangle},
in order to control the size of the space--time clusters,
we define a quantity
%denoted by
$\diam \STC(s,t)$ depending on a time interval $[s,t]$.
This quantity is larger than the increase of the maximum of the diameters
of the
space--time clusters created between the times $s$ and $t$. Moreover
this quantity is subadditive with respect to the time; see Lemma \ref
{triangle}.
Our strategy is to look
at the successive visits to $\eta$ and the excursions outside
of $\eta$. Suppose that $\eta$ has only one connected component.
The creation of
a large space--time cluster in a fixed direction has to be achieved during
an excursion outside of $\eta$. Indeed, each time the process comes back
to $\eta$, the growth of the
space--time clusters restarts almost from scratch.\looseness=1

\begin{figure}

\includegraphics{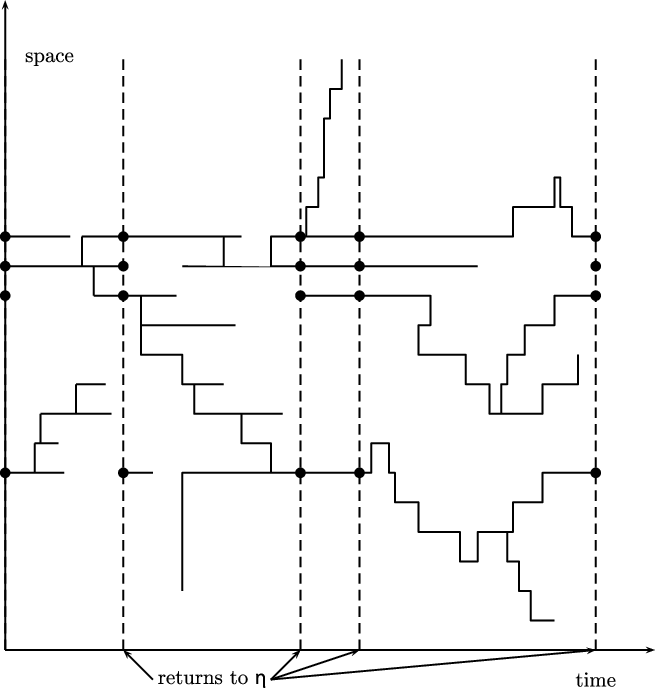}

\caption{Evolution of a STC in dimension 1.}\label{fig4}
\end{figure}

Thus if a large
space--time cluster is created before the exit of $\ocA$, then
it has to be created
during an excursion outside of $\eta$.
The situation is more complicated when the bottom
$\eta$ has several connected components. Indeed, the
space--time clusters associated to one connected component might
change between two consecutive visits to $\eta$. We prove in
Section~\ref{anad} that this does not happen: at each
visit to $\eta$, a given
connected
component of $\eta$ always belong to the same
space--time cluster. This is a consequence of Lemma \ref{stccc}.
Figure \ref{fig4} shows an example of the
space--time clusters associated to a configuration $\eta$
having two connected components. On the evolution depicted in the figure,
the space--time clusters containing
the lower component of $\eta$ at the times of the first two returns
are distinct. We will prove that this cannot occur as long as the process
stays in the cycle compound $\ocA$ (this is the purpose of
Lemma \ref{stccc}).

%f4 #&#

We rely then on a technique going back to the theory of simulated annealing,
which consists of removing the bottom $\eta$ from $\ocA$,
decomposing
$\ocA\setminus\{ \eta \}$ into its maximal cycle compounds and studying
the jumps of the process between these maximal cycle compounds until
the exit of
$\ocA\setminus\{ \eta \}$.
As before,
we show that, before exiting
$\ocA\setminus\{ \eta \}$,
the process is unlikely to make
a large number of jumps between these maximal cycle compounds.
This step is very similar to the initial step, when we considered a general
set $\cD$.
For the clarity of the exposition, we prefer to
repeat the argument rather than to introduce additional notation and make
a general statement.
Using the subadditivity of
$\diam \STC(s,t)$, we conclude that a large space--time cluster has
to be created during a visit to a maximal cycle compound of
$\ocA\setminus\{ \eta \}$.
Now each cycle compound included in
$\ocA\setminus\{ \eta \}$
has a depth strictly smaller than the depth of $\ocA$.
Using the induction hypothesis, we have a control on the
space--time clusters created during each visit to these cycle compounds.
Combining the estimate provided by the induction hypothesis and
the estimate on the number of cycle compounds of
$\ocA\setminus\{ \eta \}$
visited by the process, we obtain a control on the size of the
space--time clusters created during an excursion in
$\ocA\setminus\{ \eta \}$.
Using the estimates presented in Section~\ref{cycles}, we can also
control the
number of visits to $\eta$ before the exit of $\ocA$.
The induction step is completed by combining all the previous
estimates.\looseness=1
%In our setting the cardinality of the state space
%grows with $\b$.
%Nevertheless,
%and we use renewal ideas
%to show that a large space--time cluster has to be formed during an
%excursion outside the bottom $\eta$.
%We then decompose $\ocA\setminus\{ \eta \}$ into its maximal
%cycle compounds.
%This allows us to obtain probability estimates through an
%inductive procedure.
%To realize this program, we need to control the size of the
%space--time clusters with a subadditive quantity.
%Let us consider the Ising model in a volume $Q=\L(\log\b)$.
%Let $\overline{\cA}$ be a cycle compound in $\cC$ and $\h$ its
%unique minimizer for the energy.
%Let us set $\t=\t_{\overline{\cA}^c}$.

%s5.1 #&#
\subsection{Basic definitions and properties}
\label{badp}
Let $\Lambda$ be a subset of $\Z^d$ and let $(\sigma_{\L,t})_{t\geq0}$
be a continuous-time trajectory in
$\{ -1,+1 \}^\Lambda$.
%$\sigma_t(x): \Lambda\times\R^+ \rightarrow\{ -1,+1 \}$.
We endow the set of the space--time points $\Lambda\times\R^+$
with the following connectivity relation:
% We denote by $(x,t)\longleftrightarrow(y,s)$ the event that
the two space--time points $(x,t)$ and $(y,s)$
are connected if $\sigma_{\L,t}(x) = \sigma_{\L,s}(y) = +1$
and:

$\bullet$
either $s=t$ and $|x-y| \le1$;

$\bullet$ or
$x=y$ and $\sigma_{\L,u}(x)=+1$ for $u\in[\min(s,t),\max(s,t)]$.

A space--time cluster of the trajectory $(\sigma_{\L,t})_{t\geq0}$
is a maximal connected component of space--time points.
For $u\leq s\in\R^+$, we denote by $\STC(u,s)$ the space--time clusters
of the trajectory restricted to the time interval $[u,s]$.
Sometimes we deal with a specific initial condition $\alpha$
and boundary conditions ${\zeta}$.
We denote by
$\STC(\s^{\alpha,{\zeta}}_{\L,t},s\leq t\leq u)$
%For $u\leq s\in\R^+$, we denote by $\STC(u,s)$
the space--time clusters
of the trajectory of the process
$(\s^{\alpha,{\zeta}}_{\L,t})_{t\geq0}$
restricted to the time interval $[u,s]$.\vspace*{1pt}

%
%The space--time cluster ($\STC$) at time $t$,
%for the process $\h$,
%that contains the point $(x,s)$
%is the set of points that are connected to $(x,s)$ by a sequence
%of space--time points, i.e.,
%(y_0,s_0),\ldots,(y_k,s_k)\cr
%%\sgraffa{y_i,s_i}_{i=0,\ldots,k}  \mathrm{ with }
%%y_i\in\Phi, s_i\in[0,t], \Acapo
%(y_0,s_0)= (x,s), (y_k,s_k)=(y,s'), \cr
%(y_i,s_i)\longleftrightarrow(y_{i+1},s_{i+1})
% \Big\}.
%Notice that the sites in a $\STC$ must be
%in $\Phi$ so that the boundary condition cannot link different $\STC$.
%%\color{black}
The graphical construction updates the configuration
in two different places independently until a space--time cluster
connects the two places.
We state next a refinement of Lemma 2 of \cite{DS1}, which allows us
to compare processes defined in different volumes or with different
boundary conditions
via the graphical construction described in Section~\ref{graphi}.
%Let $\Lambda_1 \subset\Lambda_2$ be two subsets of $\Z^d$.
%We consider the coupled processes
%$\sigma_{\Lambda_1,t}$ and
%$\sigma_{\Lambda_2,t}$
%in $\Lambda_1$ and $\Lambda_2$ with minus boundary conditions
%on a time interval $[s,u]$.
%%Then, $\STC_{\Phi_1}(s) [\s_{\Phi_2}]=\STC^{s,\Phi_1} [\s_{\Phi_1}]$.
%%In other words, all
%The space--time clusters of the process in $\Lambda_2$ that are not
%a space--time cluster
%of the process in $\Lambda_1$ contain a site in
%$\Lambda_2 \setminus\Lambda_1$.
%Let $A$ be a subset of $\Sigma$. If $\cC$ is a STC
%for the dynamics in $\Sigma$ with $n\pm$ boundary conditions
%such that
%then $\cC$ is also a STC for the dynamics restricted to $A$ with
%$n\pm$ boundary conditions.
%
%For $A$ a subset of $\Z^d$, we define its outer vertex boundary
%$\dout A$ as
%$$\dout A =
%
%le5.1 #&#
\begin{lemma}
\label{cresSTC}
Let $\L$ be a subset of $\Z^d$, and let ${\zeta}$ be a boundary
condition on~$\L$.
Let $x$ be a site of the exterior boundary of $\L$ such that ${\zeta}(x)=+1$.
If $\cC$ is a $\STC$
for the dynamics in $\L$ with ${\zeta}$ as boundary conditions,
and $\cC$ is such that $x$ is not the neighbor of a point of $\cC$, then
$\cC$ is also a $\STC$
for the dynamics in $\L$ with ${\zeta}^x$ as boundary conditions.
\end{lemma}
\begin{pf}
We denote by $\alpha$ the initial configuration.
From the coupling, we have
\[
\forall t\geq0, \forall y\in\L\qquad \s^{\alpha,{\zeta}^x}_{\L,t}(y) \leq
\s^{\alpha,{\zeta}}_{\L,t}(y).
\]
Let $\cC$ be a $\STC$ in
$\STC(\s^{\alpha,{\zeta}}_{\L,t},s\leq t\leq u)$ and suppose that
$\cC$ does not belong to
$\STC(\s^{\alpha,{\zeta}^x}_{\L,t},s\leq t\leq u)$.
Necessarily, there exists a space--time point
$(y,t)$
such that
\[
(y,t)\in\cC,\qquad \s^{\alpha,{\zeta}^x}_{\L,t}(y)=-1,\qquad \s^{\alpha,{\zeta}}_{\L,t}(y)=+1.
\]
We consider the set of the space--time points satisfying the above
condition, and we denote by $(y^*,t^*)$ the space--time point such that
$t^*$ is minimum. This is possible since the number of spin flips in a
finite box is finite in a finite time interval, and moreover the
trajectories are right continuous. At time $t^*$, the spin at site
$y^*$ becomes $+1$ in the process
$(\s^{\alpha,{\zeta}}_{\L,t})_{t\geq0}$, and it remains equal to $-1$
in $(\s^{\alpha,{\zeta}^x}_{\L,t})_{t\geq0}$. We examine next the
neighbors of $y^*$. Let $z$ be a neighbor of $y^*$ in $\L$. If
$\s^{\alpha,{\zeta}}_{\L,t^*}(z)=-1$, then
$\s^{\alpha,{\zeta}^x}_{\L,t^*}(z)=-1$ as well. Suppose that
$\s^{\alpha,{\zeta}}_{\L,t^*}(z)=+1$. The spin at $z$ does not change
at time $t^*$, thus for $s<t^*$ close enough to $t^*$, we have also
$\s^{\alpha,{\zeta}}_{\L,s}(z)=+1$. This implies that $\{ z
\}\times[s,t^*]$ is included in $\cC$. From the definition of
$(y^*,t^*)$, we have that
\[
\forall u\in\bigl[s,t^*\bigr]\qquad \s^{\alpha,{\zeta}^x}_{\L,u}(z)=+1.
\]
We conclude that the neighbors of $y^*$ in $\L$ have the same spins in
$\s^{\alpha,{\zeta}^x}_{\L,t^*}$ and in $\s^{\alpha,{\zeta}}_{\L,t^*}$.
Therefore $y^*$ must have a neighbor in $\Z^d\setminus\L$ whose spin is
different in $\s^{\alpha,{\zeta}^x}_{\L,t^*}$ and in
$\s^{\alpha,{\zeta}}_{\L,t^*}$. The only possible candidate is $x$.
\end{pf}

The next corollary is very close to Lemma 2 of \cite{DS1}.
%
%co5.2 #&#
\begin{corollary}\label{deghstc}
Let $\L_1\subset\L_2$ be two subsets of $\Z^d$, let $\alpha$ be an
initial configuration in $\L_2$
and let ${\zeta}$ be a boundary condition on $\L_2$.
If no $\STC$ of the process
$(\s^{\alpha,{\zeta}}_{\L_2,t},s\leq t\leq u)$
intersects both $\L_1$ and the inner boundary of $\L_2$, then
\[
\forall t\in[s,u]\qquad \s^{\alpha,{\zeta}}_{\L_2,t}|_{\L_1} =
\s^{\alpha,-}_{\L_2,t}|_{\L_1}.
\]
\end{corollary}
We define the diameter
$\diam \cC$ of a space--time cluster $\cC$
by
\[
\diam \cC = \sup \bigl\{ |x-y|_\infty\dvtx  (x,s), (y,t) \in\cC
\bigr\},
\]
where $|\cdot|_\infty$ is the supremum norm given by
\[
\forall x=(x_1,\ldots,x_d)\in\Z^d\qquad
|x|_\infty = \max_{1\leq i\leq d} |x_i|.
\]
Thus $\diam \cC$ is
the
%$\infty$-norm
diameter of the spatial projection of $\cC$.
%Given a unit vector $\underline e$, parallel to one of the coordinate
%axes,
%we denote by
%$\diam_{\underline e}  \STC= \sup_{\{x,s\},\{x',s'\}} (x-x')\cdot
%the diameter of the projection of a $\STC$ in the direction $
% We denote by $\STC_{\L}(t)$ the set of all
%the $\STC$ for the process $\h$
%that are contained in $\L\times[0,t]$.
%s5.2 #&#
\subsection{The bottom of a cycle compound}
We prove here that, when $h$ is irrational, the bottom of a cycle
compound of the Ising model
contains a unique configuration.
Throughout the section, we consider a finite box $Q$ endowed with a
boundary condition $\xi$.
To alleviate the formulas, we write simply $H$ instead of $H_Q^\xi$.

%le5.3 #&#
\begin{lemma}\label{nuovo}
%Let $h \in\R\setminus\Q$.
% Let $\cA_1,\ldots, \cA_n$ be a collection of cycles included into
% $\cA_{\minus}$ such that
% $H(\widetilde B(\cA_i))= \lambda \forall i=1,\ldots,n$
% and $\overline{\cA}=\bigcup_{i=1}^n \tonda{\cA_i \cup\widetilde B(
% is a connected set.
Suppose that $h$ is irrational.
Let $\h$ be a minimizer of the energy in a cycle compound $\overline
{\cA}$.
Then for any ${\zeta}\in\overline{\cA}$,
${\zeta}\cup\h\in\overline{\cA}$ and ${\zeta}\cap\h\in
\overline{\cA}$.
\end{lemma}
\begin{pf}
Let $\eta$ belong to the bottom of $\ocA$.
We assume that
$\overline{\cA}$ is not a singleton; otherwise there is nothing to prove.
Let
$\o=(\o_1,\ldots,\o_n)$
be a path in $\overline{\cA}$ that goes from $\h$ to ${\zeta}$.
%We can assume that for each $k$, $\o_k \not= \o_{k+1}$ (automatique)
We associate with $\o$
a \emph{slim} path
\[
\o\cap\eta = (\o_1\cap\h,\ldots,\o_n\cap\h)
\]
and a \emph{fat} path
\[
\o\cup\eta = (\o_1\cup\h,\ldots,\o_n\cup\h).
\]
Suppose that the thesis is false, and let us set
\[
\k^* = \min \{ k\geq1\dvtx  \o_k\cap\h\notin\overline{\cA}
\mbox{ or } \o_k\cup\h\notin\overline{\cA} \}. %\mbox{ or } \Acapo&& H(\o_k)<\max\tonda{H(\o_k\cap\h),H(\o_k\cup\h)}}.
\]
Notice that $\k^*$ is larger than or equal to $2$.
We will use the attractive inequality
\[
H(\o_k\cap\h)+H(\o_k\cup\h) \le H(
\o_k) + H(\h)
\]
and the fact that $\h$ is a minimizer of the energy in $\overline{\cA}$.
Let us set
\[
\lambda = E \bigl(\overline{\cA}, \{-1,+1\}^{\L} \setminus\overline{
\cA} \bigr).
\]
First, for any $k<\k^*$, the above inequality yields that
\[
\max \bigl( H(\o_{k}\cap\h),H(\o_{k}\cup\h) \bigr) \leq H(
\o_{k})\leq \lambda.
\]
%
%We prove next that this inequality is still valid for $k=\k^*$.
%For any $k\geq2$, the configurations
The configurations
$\o_{\k^*}$ and $\o_{\k^*-1}$ differ for the spin in a single site.
We say that the $\k^*$th spin flip is inside
(resp., outside) $\h$
if this site has a plus spin
(resp., a minus spin)
in $\h$,
that is, if $\o_{\k^*} \vartriangle\o_{{\k^*}-1} \subset\h$
(resp.,
$\o_{\k^*} \vartriangle\o_{{\k^*}-1} \not\subset\h$).
We distinguish two cases, according to the position of
the $\k^*$th spin flip with respect to $\h$:

(i) if %$\o_{\k^*} \vartriangle\o_{k-1} \subset\h$,
the $\k^*$th spin flip is inside $\h$,
then $\o_{\k^*} \cup\h=\o_{\k^*-1} \cup\h$,
so that only the slim path moves and exits $\ocA$ at index $\k^*$.
%but this is larger than or equal to zero.
%This means in particular that $\o_{\k^*-1}\cup\h\in\overline{\cA}$
%and, again by \eqref{attr}, we get
%$$H(\o_{\k^*}) -H(\o_{\k^*}\cup\h)\ge
%H(\o_{\k^*-1}\cap\h)-H(\h)\ge0.$$
Thus
\[
\o_{\k^*-1}\cap\h\in\ocA,\qquad \o_{\k^*}\cap\h\notin\ocA
\]
and these two configurations communicate, therefore
\[
\max\bigl( H( \o_{\k^*-1}\cap\h), H( \o_{\k^*}\cap\h) \bigr)
\geq \lambda.
\]
We distinguish again two cases:

$\bullet$  $H( \o_{\k^*-1}\cap\h)\geq\lambda$. Since
$ H( \o_{\k^*-1}\cap\h)\leq
H( \o_{\k^*-1})\leq\lambda$,
then
$\o_{\k^*-1}\cap\h$ and
$\o_{\k^*-1}$ have both an energy equal to $\lambda$, and by
Lemma \ref{irrazionale},
we conclude that
$\o_{\k^*-1}\cap\h=
\o_{\k^*-1}$
and
$\o_{\k^*-1}$ is included in $\h$.
Since we are assuming that the slim path moves at step $\k^*$,
the original path and the slim path undergo the same spin flip so that they
must coincide also at step $\k^*$, contradicting the assumption that
$\o_{\k^*}\cap\h\notin\overline{\cA}$.
%The case where the fat path exits $\overline{\cA}$ is analogous.

$\bullet$  $H( \o_{\k^*}\cap\h)\geq\lambda$.
By the attractive inequality
\[
H(\o_{\k^*}) -H(\o_{\k^*}\cap\h) \ge H(\o_{\k^*-1}\cup
\h)-H(\h) \geq 0,
\]
whence
\[
H(\o_{\k^*}\cap\h) \leq H(\o_{\k^*}) \leq \lambda.
\]
Thus
$\o_{\k^*}\cap\h$ and
$\o_{\k^*}$ have both an energy equal to $\lambda$. By Lemma \ref
{irrazionale},
we conclude that
$\o_{\k^*}\cap\h=
\o_{\k^*}$,
contradicting the assumption that
$\o_{\k^*}\cap\h\notin\overline{\cA}$.

We consider next the second case. The argument is very similar in the
two dual cases (i) and (ii), yet it seems necessary to handle them
separately.

%E(\ocA,\cX\setminus\ocA)
(ii) if %$\o_{\k^*} \vartriangle\o_{k-1} \subset\h$,
the $\k^*$th spin flip is outside $\h$,
then $\o_{\k^*} \cap\h=\o_{\k^*-1} \cap\h$,
so that only the fat path moves and exits $\ocA$ at index $\k^*$.
%but this is larger than or equal to zero.
%This means in particular that $\o_{\k^*-1}\cup\h\in\overline{\cA}$
%and, again by \eqref{attr}, we get
%$$H(\o_{\k^*}) -H(\o_{\k^*}\cup\h)\ge
%H(\o_{\k^*-1}\cap\h)-H(\h)\ge0.$$
Thus
\[
\o_{\k^*-1}\cup\h\in\ocA,\qquad \o_{\k^*}\cup\h\notin\ocA
\]
and these two configurations communicates, therefore
\[
\max\bigl( H( \o_{\k^*-1}\cup\h), H( \o_{\k^*}\cup\h) \bigr)
\geq \lambda.
\]
We distinguish again two cases:

$\bullet$  $H( \o_{\k^*-1}\cup\h)\geq\lambda$. Since
$ H( \o_{\k^*-1}\cup\h)\leq
H( \o_{\k^*-1})\leq\lambda$,
then
$\o_{\k^*-1}\cup\h$ and
$\o_{\k^*-1}$ have both an energy equal to $\lambda$, and by
Lemma \ref{irrazionale},
we conclude that
$\o_{\k^*-1}\cup\h=
\o_{\k^*-1}$
and
$\o_{\k^*-1}$ contains $\h$.
Since we are assuming that the fat path moves at step $\k^*$,
the original path and the fat path undergo the same spin flip so that they
must coincide also at step $\k^*$, contradicting the assumption that
\mbox{$\o_{\k^*}\cup\h\notin\overline{\cA}$}.
%The case where the fat path exits $\overline{\cA}$ is analogous.

$\bullet$  $H( \o_{\k^*}\cup\h)\geq\lambda$.
By the attractive inequality
\[
H(\o_{\k^*}) -H(\o_{\k^*}\cup\h) \ge H(\o_{\k^*-1}\cap
\h)-H(\h) \geq 0,
\]
whence
\[
H(\o_{\k^*}\cup\h) \leq H(\o_{\k^*}) \leq \lambda.
\]
Thus
$\o_{\k^*}\cup\h$ and
$\o_{\k^*}$ have both an energy equal to $\lambda$. By Lemma \ref
{irrazionale},
we conclude that
$\o_{\k^*}\cup\h=
\o_{\k^*}$,
contradicting the assumption that
$\o_{\k^*}\cup\h\notin\overline{\cA}$.
\end{pf}
\begin{lemma}\label{fondo}
Suppose that $h$ is irrational.
The bottom
$\bottom(\overline{\cA})$ of any cycle compound $\overline{\cA}$
contains a single configuration.
%$|\cF(\overline{\cA})|=1$ for any cycle compound $\overline{\cA}$.
\end{lemma}
\begin{pf}
If $\h_1,\h_2 \in\bottom(\overline{\cA})$, then by Lemma \ref{nuovo}
we have also
$\h_1 \cup\h_2 \in\overline{\cA}$ and $\h_1 \cap\h_2 \in
\overline{\cA}$,
so that $H(\h_1) + H(\h_2) \le H(\h_1 \cup\h_2)+H(\h_1 \cap\h
_2) $.
But by the attractive inequality,
\[
H(\h_1 \cup\h_2)+H(\h_1 \cap
\h_2) \le H(\h_1) + H(\h_2),
\]
so that $\h_1\cup\h_2$ and $\h_1\cap\h_2$ are also in
$\bottom(\overline{\cA})$. Lemma \ref{irrazionale} implies that
$\h_1\cup\h_2=\h_1\cap\h_2$, showing that $\h_1=\h_2$.
\end{pf}

%s5.3 #&#
\subsection{The space--time clusters in a cycle compound}
\label{anad}
In this section, we study
%the energy landscape in
some properties of the paths contained in suitable cycle compounds.
In order to avoid unnecessary notation,
with a slight abuse of terms, we consider space--time clusters
associated to a discrete time trajectory.
%the connected components in a path.
In other words, in this section the word ``time''
means ``index of the configuration in the trajectory,''
and the space--time clusters considered here are
pure geometrical objects.
We will use these geometrical results in order to control
the diameter of the space--time clusters of our processes.

As in the previous
section, we consider a finite box $Q$
endowed with a boundary condition $\xi$.
To alleviate the formulas, we write simply $H$ instead of $H_Q^\xi$.
A~connected component of a configuration $\sigma$ is a maximal
connected subset of the plus sites of $\sigma$
\[
\bigl\{ x\in\Z^d\dvtx \sigma(x)=+1 \bigr\},
\]
two sites being connected if they
are nearest neighbors on the lattice.
We denote by $\cC(\sigma)$ the connected components of $\sigma$.
If $C\in\cC(\sigma)$, then we define its energy as
\[
H(C) = \bigl|\bigl\{ \{ x,y \}\dvtx x\notin C, y\in C, |x-y|=1 \bigr\} \bigr|-h|C|.
\]
In particular, we have
% for any configuration $\sigma$
%
\[
H(\sigma) = \sum_{C\in\cC(\sigma)}H(C).
\]
Let $\o=(\o_0,\ldots,\o_r)$ be a path of configurations in the box $Q$.
We endow the set of the space--time points $Q\times\N$
with the following connectivity relation associated to $\o$:
% We denote by $(x,t)\longleftrightarrow(y,s)$ the event that
the two space--time points $(x,i)$ and $(y,j)$
are connected if $\o_i(x) = \o_j(y) = +1$
and:

$\bullet$
either $i=j$ and $|x-y| \le1$;

$\bullet$ or
$x=y$ and $|i-j|=1$.

A space--time cluster of the path $\o$
is a maximal connected component of space--time points in $\o$.
We consider a domain $\cD$, which is a set of configurations
satisfying the following hypothesis.\vspace*{9pt}

\textit{Hypothesis on $\cD$.} The configurations in $\cD$
are such that:

$\bullet $ There exists $v_\cD$ (independent of $\beta$)
such that
$|\sigma|\leq v_\cD$ for any $\sigma\in\cD$.\eject

%$\bullet $For any connected
%component $C$ of a configuration of $\cD$,
%we have $H(C)>H(\minus)$.

$\bullet $ If $\sigma\in\cD$
and $C$ is a connected component of $\sigma$,
then we have $H(C)>H(\minus)$.
%$\bullet $For any $\sigma\in\cD$,
%any connected component $C$ of $\sigma$,
%we have $H(C)>H(\minus)$.

$\bullet $ If $\sigma\in\cD$
and $\eta$ is
%a configuration
such that
$\eta\subset\sigma$ and
$H(\eta)\leq H(\sigma)$,
then
$\eta\in\cD$.

%le5.5 #&#
\begin{lemma}\label{stccc}
Let $\overline{\cA}$ be
a cycle compound included in $\cD$,
and let $\eta$ be the unique configuration of
$\bottom(\overline{\cA})$.
Let
$\o = (\o_0,\ldots,\o_r)$
be a path in $\overline{\cA}$ starting
at $\eta$ and ending at $\eta$.
Let $C$ be a connected component of $\eta$. Then
the space--time sets
$C\times\{0\}$ and
$C\times\{r\}$ belong to the same space--time cluster of $\o$.
\end{lemma}
\begin{pf}
From Lemma \ref{fondo}, we know that
$\bottom(\overline{\cA})$ is reduced to a single configuration $\h$.
By Lemma \ref{nuovo}, the path
\[
\o\cap\eta = (\o_0\cap\h,\ldots,\o_r\cap\h)
\]
is still a path
in $\overline{\cA}$ that goes from $\eta$ to $\eta$.
Moreover, the space--time clusters of
$\o\cap\eta$ are included in those of $\o$, therefore
it is enough to prove the result for the path
\mbox{$\o\cap\eta$}.
Let $\ot$ be the path obtained from
$\o\cap\eta$ by removing all the space--time clusters of
$\o\cap\eta$ which do not intersect $\eta\times\{0\}$.
The path $\ot$ is still admissible, that is, it is a sequence
of configurations such that each configuration communicates
with its successor.
Let $i\in\{ 0,\ldots,r \}$. We have
$\ot_i \subset \o_i\cap\eta$.
%Thanks to the hypothesis on the domain $\cD$,
Since
$\ot_i$ is obtained from $\o_i\cap\eta$ by removing some connected
components of
$\o_i\cap\eta$, the second hypothesis on
the domain $\cD$ yields that
$H(\ot_i) \leq H(\o_i\cap\eta)$. With the help of the third
hypothesis on $\cD$,
we conclude that
$\ot_i$ is in $\cD$.
%$$\forall i\in\{ 0,\ldots,r \}
%H(\ot_i) \leq H(\o_i\cap\eta),$$
In particular the whole path $\ot$ stays in $\cD$.
Suppose that the path $\ot$ leaves $\ocA$ at some index $i$,
so that
$\ot_{i}\neq\o_{i}\cap\eta$.
We consider two cases:

$\bullet$ $\ot_{i-1}=\o_{i-1}\cap\eta$. In this case,
the spin flip between
$\o_{i-1}\cap\eta$
and
$\o_{i}\cap\eta$ creates a new STC which does not intersect
$\eta\times\{0\}$, hence
$\ot_{i}=\o_{i-1}\cap\eta$.
This contradicts the fact that
$\ot$ leaves $\ocA$ at index $i$.

$\bullet$ $\ot_{i-1}\neq\o_{i-1}\cap\eta$. Since we
have also
$\ot_i\neq\o_{i}\cap\eta$, then
by Lemma \ref{irrazionale}
we have the strict inequality
\[
\max \bigl(H(\ot_{i-1}),H(\ot_i) \bigr) < \max \bigl(H(
\o_{i-1}\cap\eta),H(\o_i\cap\eta) \bigr) \leq E(\ocA,\cX
\setminus\ocA).
\]
However, since
$\ot$ leaves $\ocA$ at index $i$, we have also
\[
\max \bigl(H(\ot_{i-1}),H(\ot_i) \bigr) \geq E(\ocA,\cX
\setminus\ocA),
\]
which is absurd. Thus
the path $\ot$ stays also in $\ocA$.
Since
\[
H(\ot_r) \leq H(\o_r\cap\eta),\qquad \ot_r\subset
\o_r\cap\eta=\eta,
\]
we have
$\ot_r=\eta$
by Lemma \ref{irrazionale}.
The path $\ot$ is included in
$\eta\times\{ 0,\ldots,r \}$,
hence,
for any connected component $C$ of $\eta$, the
space--time cluster of $\ot$ containing $C\times\{ r \}$ is
included in
$C\times\{ 0,\ldots,r \}$, so that its intersection
with $\eta\times\{ 0 \}$,
which is not empty by construction,
must be equal to
%is
$C\times\{ 0 \}$.
\end{pf}
%
%pb, ils ne sont pas forcement connectes?
%s5.4 #&#
\subsection{Triangle inequality for the diameters of the STCs}
\label{sectriangle}
In the sequel, we consider a trajectory of the process
$(\s_{Q,t},t\geq0)$ in a finite box $Q$,
and we study its space--time clusters.
For $s<t$, we define
\[
\diam \STC(s,t) = \max \biggl(\mathop{\sum_{\cC\in\STC
(s,t)}}_{\cC\cap(Q\times\{ s,t \})\neq\varnothing}
\diam \cC, \mathop{\max_{\cC\in\STC(s,t)}}_{ \cC\cap(Q\times\{
s,t \})=\varnothing} \diam \cC
\biggr).
\]
%
%$${
%,
The main point of this awkward definition is the following triangle
inequality.

%le5.6 #&#
\begin{lemma}\label{triangle}
For any $s<u<t$, we have
\[
\diam \STC(s,t) \leq \diam \STC(s,u) + \diam \STC(u,t).
\]
\end{lemma}
\begin{pf}
When we look at the restriction
to the time intervals
$(s,u)$ and $(u,t)$
of a $\STC$ in $\STC(s,t)$ which is alive at time $u$,
this $\STC$ splits into several $\STC$ belonging to
$\STC(s,u)\cup\STC(u,t)$.
Yet the diameter of
the initial $\STC$ is certainly less than the sum of all the diameters
of the $\STC$
in
$\STC(s,u)\cup\STC(u,t)$ which are alive at time $u$.
The proof is quite tedious; however, since this inequality is fundamental
for our argument we provide a detailed verification.
%We write
%$${
%{
%{\cC\in\STC(s,t)}{
%{\cC\cap(Q\times\{ u \})\neq\varnothing}
%}}
% +
%{
%{\cC\in\STC(s,t)}
%{\cC\cap(Q\times\{ s,t \})\neq\varnothing}}
%{\cC\cap(Q\times\{ u \})=\varnothing}}
% +
First, we have
%that
%
\[
\mathop{\mathop{\sum_{\cC\in\STC(s,t)}}_{
{\cC\cap(Q\times\{ s,t \})\neq\varnothing}}}
_{\cC\cap(Q\times\{ u \})\neq\varnothing} \diam \cC \leq \mathop{\sum_{\cC\in\STC(s,u)}}_
{\cC\cap(Q\times\{ u \})\neq\varnothing}
\diam \cC + \mathop{\sum_{\cC\in\STC(u,t)}}_
{\cC\cap(Q\times\{ u \})\neq\varnothing}
\diam \cC.
\]
Next, if $\cC\in\STC(s,t)$ and
$\cC\cap(Q\times\{ u \})=\varnothing$,
then
$\cC\in\STC(s,u)\cup\STC(u,t)$.
Thus
\[
\mathop{\mathop{\sum_{\cC\in\STC(s,t)}}_{
\cC\cap(Q\times\{ s,t \})\neq\varnothing}}_{
{\cC\cap(Q\times\{ u \})=\varnothing
}}
\diam \cC \leq \mathop{\mathop{\sum_{\cC\in\STC(s,u)}}_
{\cC\cap(Q\times\{ s \})\neq\varnothing}}_
{\cC\cap(Q\times\{ u \})=\varnothing}
\diam \cC + \mathop{\mathop{\sum_{\cC\in\STC(u,t)}}_
{\cC\cap(Q\times\{ u \})=\varnothing}}_
{\cC\cap(Q\times\{ t \})\neq\varnothing}
\diam \cC.
\]
Summing the two previous inequalities, we get
\begin{eqnarray*}
\mathop{\sum_{\cC\in\STC(s,t)}}_
{\cC\cap(Q\times\{ s,t \})\neq\varnothing} \diam \cC &\leq&
\mathop{\sum_{\cC\in\STC(s,u)}}_{
\cC\cap(Q\times\{ s,u \})\neq\varnothing} \diam \cC +
\mathop{\sum_{\cC\in\STC(u,t)}}_
{\cC\cap(Q\times\{ u,t \})\neq\varnothing} %{\cC\cap(Q\times\{ t \})\neq\varnothing}
\diam \cC
\\
&\leq& \diam \STC(s,u) + \diam \STC(u,t).
\end{eqnarray*}
%
%either
%$\cC\in\STC(s,u)$ or
%$\cC\in\STC(u,t)$.
Moreover,
if
$\cC\in\STC(s,t)$,
$\cC\cap(Q\times\{ s,t \})=\varnothing$ and
$\cC\cap(Q\times\{ u \})\neq\varnothing$,
then
\[
\diam \cC \leq %\diam \STC(s,t) =
\mathop{\mathop{\sum_{\cC\in\STC(s,u)}}_
{\cC\cap\bigl(Q\times\{ u \}\bigr)\neq\varnothing}}_
{\cC\cap(Q\times\{ s \})=\varnothing}
\diam \cC + \mathop{\mathop{\sum_{\cC\in\STC(u,t)}}_
{\cC\cap(Q\times\{ u \})\neq\varnothing}}_
{\cC\cap(Q\times\{ t \})=\varnothing}
\diam \cC.
\]
Finally
if
$\cC\in\STC(s,t)$,
$\cC\cap(Q\times\{ s,u,t \})=\varnothing$,
then
$\cC\in\STC(s,u)\cup\STC(u,t)$ and
\[
\diam \cC \leq %\diam \STC(s,t) =
\mathop{\max_{\cC\in\STC(s,u)}}_{
\cC\cap(Q\times\{ s,u \})=\varnothing}
\diam \cC + \mathop{\max_{\cC\in\STC(u,t)}}_{
\cC\cap(Q\times\{ u,t \})=\varnothing} \diam \cC.
\]
The two previous inequalities yield
\[
\mathop{\max_{\cC\in\STC(s,t)}}_{
\cC\cap(Q\times\{ s,t \})=\varnothing} \diam \cC \leq \diam
\STC(s,u) + \diam \STC(u,t),
\]
and the proof is complete.
\end{pf}

%s5.5 #&#
\subsection{The diameter of the space--time clusters}
\label{diamostc}
We consider boxes that grow slowly with $\beta$.
This creates a major
complication in the description of the energy landscape, but it allows
us to obtain very strong estimates
that will be used to control entropy effects in the dynamics of growing
droplets.
We make the following hypothesis on the volume of the box
$Q$.\vspace*{9pt}

\textit{Hypothesis on $Q$.} The box $Q$ is such that
$|Q| = \exp o(\ln\beta)$, which means that
\[
\lim_{\beta\to\infty} \frac{\ln|Q|}{\ln\beta} = 0.
\]
Let $n\in\{ 0,\ldots,d \}$. As in Section~\ref{anad},
we consider
%We consider a boundary condition $\z$ on $Q$ and
%domain $\cD$, which is a set of configurations
a set of configurations
$\cD$
in the box $Q$ satisfying the following hypothesis.\vspace*{9pt}

\textit{Hypothesis on $\cD$.} The configurations in $\cD$
are such that:

$\bullet$ There exists $v_\cD$ (independent of $\beta$)
such that
$|\sigma|\leq v_\cD$ for any $\sigma\in\cD$.

%$\bullet $For any connected
%component $C$ of a configuration of $\cD$,
%we have $H(C)>H(\minus)$.

$\bullet$ If $\sigma\in\cD$
and $C$ is a connected component of $\sigma$,
then we have
\[
H_Q^{n\pm}(C) > H_Q^{n\pm}(\minus).
\]

$\bullet$ If $\sigma\in\cD$
and $\eta$ is
%a configuration
such that
$\eta\subset\sigma$ and
$H_Q^{n\pm}(\eta)\leq H_Q^{n\pm}(\sigma)$,
then
$\eta\in\cD$.

The hypothesis on $\cD$ ensures that
the number of the energy values
of the configurations in $\cD$ with $n\pm$ boundary conditions
is bounded by a value independent
of $\beta$.
Indeed, for any $\sigma\in\cD$,
\[
H_Q^{n\pm}(\sigma) = \sum_{C\in\cC(\sigma)}H_Q^{n\pm}(C),
\]
where
$\cC(\sigma)$ is the set of the connected components of $\sigma$.
Yet
there are at most $v_\cD$ elements in $\cC(\sigma)$,
and any element of $\cC(\sigma)$ has volume at most $v_\cD$;
hence the number of possible values for $H$
is at most
$c(d)^{(v_\cD)^2}$ where $c(d)$ is a constant depending on the
dimension $d$ only.
%In particular, we have
%$\lambda_1=H(\minus)$ and
%$\lambda_p=\height(\cC_d)$.
Let next
\[
\delta_0<\delta_1<\cdots<\delta_p
\]
be the possible
values for the difference of the energies of two configurations of~$\cD$,
%depth of a cycle compound included in $\cD$,
that is,
\[
\{ \delta_0,\ldots,\delta_p \} = \bigl\{
\bigl|H_Q^{n\pm}(\sigma)-H_Q^{n\pm}(\eta) \bigr|_+\dvtx \sigma,\eta\in\cD \bigr\}.
\]

\vspace*{9pt}

\textit{Notation}.
We will study the space--time clusters associated to different processes.
For $\alpha$ an initial configuration
and $\zeta$ a boundary condition, we denote by
\[
\STC\bigl( \s^{\alpha,{\zeta}}_{Q,t}, s\leq t\leq u\bigr)
\]
the $\STC$ associated to the trajectory of the process
$(\s^{\alpha,{\zeta}}_{Q,t})_{t\geq0}$ during the time interval $[s,u]$.
Accordingly,
\[
\diam \STC\bigl( \s^{\alpha,{\zeta}}_{Q,t}, s\leq t\leq u\bigr)
\]
is equal to
$\diam \STC(
s,u)$
computed for the $\STC$ of the process
$(\s^{\alpha,{\zeta}}_{Q,t})_{t\geq0}$
on the time interval $[s,u]$.
%
%th5.7 #&#
\begin{theorem}\label{totcontrole}
%Let $\z$ be a boundary condition on $Q$ and
Let $n\in\{ 1,\ldots,d \}$.
For any $K>0$, there exists a value $D$ which depends
only on $v_\cD$ and $K$ such that,
for $\beta$ large enough,
we have
%There exists a positive constant $c(d,h)$ depending
%on the dimension $d$ and the magnetic field $h$ only,
%and a constant $k_0(\ocA)$ depending on $\ocA$ such that
%$$\forall\alpha\in\cD
%
%
\[
\forall\alpha\in\cD\qquad P \bigl( \diam \STC\bigl( \s^{\alpha,n\pm}_{Q,t},
0\leq t\leq\tau(\cD)\bigr)\geq D \bigr) %\left(
%for the process
% $(\s^{\alpha,\z}_{Q,t}, t\geq0)$
%}\\
\leq \exp(-\beta K).
\]
\end{theorem}

To alleviate the formulas, we drop the superscripts which do not vary, like
the boundary conditions $n\pm$ and sometimes the initial configuration
$\alpha$.
Throughout the proof we fix an integer $n\in\{ 1,\ldots,d \}$, and
$\s_{Q,t}$
stands for
$\s^{\alpha,n\pm}_{Q,t}$.
For\vspace*{1pt} $\cA$ an arbitrary set and $t\geq0$,
we define the time $\tau(\cA,t)$
of exit from $\cA$ after time $t$
\[
\tau(\cA,t) = \inf \{ s\geq t\dvtx \sigma_{Q,s}\notin\cA \}.
\]
Let $\cE$ be a subset of $\cD$.
We consider the decomposition of
$\cE$ into its
maximal cycle compounds $\ocM(\cE)$, and
we look at the successive jumps between the elements of
$\ocM(\cE)$.
%(we make the convention that $\tau(C)=\tau(C,0)$).
For $\gamma\in\cE$, we denote by
\[
\bpi(\gamma,\cE)
\]
the maximal cycle compound of
$\cE$ containing $\gamma$.
Let $\alpha\in\cE$ be the initial configuration.
We define recursively a sequence of random times and maximal cycle
compounds included in
$\cE$,
\[
\matrix{ %format\lambda &\lambda &\lambda\\
\btau_0=0,&&
\bpi_0=\bpi(\alpha,\cE),
\cr
\btau_1=\tau(
\bpi_0,\btau_0),&&\bpi_1=\bpi(
\sigma_{Q,\btau_1}, \cE),
\cr
\vdots&& \vdots
\cr
\btau_k=\tau(
\bpi_{k-1},\btau_{k-1}),&& \bpi_k=\bpi(
\sigma_{Q,\btau_k},\cE),
\cr
\vdots&& \vdots
\cr
\btau_R=\tau(
\bpi_{R-1},\btau_{R-1}),&& \bpi_R=\bpi(
\sigma_{Q,\btau_R},\cE),
\vspace*{1pt}\cr
\btau_{R+1}=\tau(\cE).&&}
\]
The sequence $(\bpi_0,\ldots,\bpi_{R-1},\bpi_R)$ is
the path of the maximal cycle compounds
in $\cE$
visited by $(\sigma_{Q,t})_{t\geq0}$,
and it is denoted by
$\bpi(\cE)$.
%It is a random variable with values in the pruned cycle path space
%$${\overline\Psi(D,G) =
%1\leq k<r \big\}\times\big\{ \{z\}\dvtx z\in D^c \big\}}$$
%of finite sequences of sets of equivalent cycles starting in $G$,
%traveling through $\overline{\Cal M}(D\setminus G)$ and ending
%in $D^c$.
We first obtain a control on the random length $R(\cE)$ of
$\bpi(\cE)$.
%
%pr5.8 #&#
\begin{proposition}\label{long}
There exists a constant $c>0$ depending only
on $v_\cD$ such that,
for any subset $\cE$ of $\cD$,
for $\beta$ large enough,
\[
\forall\alpha\in\cE, \forall r\geq1\qquad \P \bigl(R(\cE)\geq r \bigr) \leq
\frac{1}{c}\exp(-\beta cr).
\]
\end{proposition}
\begin{pf}
Let us set $\ocA_0=\bpi(\alpha,\cE)$.
We write
\[
\P\bigl(R(\cE)=r\bigr) = \sum_{\ocA_1,\ldots,\ocA_r
\in\ocM(\cE)} \P \bigl( \bpi(
\cE)= (\ocA_0,\ocA_1,\ldots,\ocA_r) \bigr).
\]
Let
${\ocA_1,\ldots,\ocA_r}$ be a fixed path in
$\ocM(\cE)$.
With the help of the Markov property, we have
\begin{eqnarray*}
&&
\P \bigl(\bpi(\cE)= (\ocA_0,\ocA_1,\ldots,
\ocA_r) \bigr)
\\
&&\qquad
= \sum_{\alpha_1\in\ocA_1\cap\partial\ocA_0,\ldots,
\alpha_r\in\ocA_r\cap\partial\ocA_{r-1}} \P \pmatrix{  \mbox{$
\bpi(\cE)= (\ocA_0,\ocA_1,\ldots,\ocA_r)$}\cr
\mbox{$\sigma_{Q,\btau_1}=\alpha_1, \ldots, \sigma_{Q,\btau_{r}}=
\alpha_r$} }
\\
&&\qquad
= \sum_{\alpha_1\in\ocA_1\cap\partial\ocA_0,\ldots,
\alpha_r\in\ocA_r\cap\partial\ocA_{r-1}} \P \bigl( \sigma^\alpha_{Q,\btau_1}=
\alpha_1 \bigr) \cdots \P \bigl( \sigma^{\alpha_{r-1}}_{Q,\btau_{r}}=
\alpha_r \bigr).
\end{eqnarray*}
Using the hypothesis on $Q$ and $\cD$, for $\varepsilon>0$
and for $\beta$ large enough,
we can bound the prefactor appearing in Corollary \ref{exitcom} by
\[
{\deg(\alpha)}^{|\cX|} \leq \exp(\beta{\varepsilon}).
\]
For $i\in\{ 1,\ldots,r \}$, let $a_i$ in $\cE$ be such that
$H(a_i)=
E(\ocA_{i-1},\cX\setminus\ocA_{i-1})$.
Applying next Corollary \ref{exitcom}, we obtain
\begin{eqnarray*}
&&
\P \bigl( \bpi(\cE)= (\ocA_0,\ocA_1,\ldots,
\ocA_r) \bigr)
\\
&&\qquad
\leq %\sum_{\alpha_1\in\ocA_1,\ldots,
\sum_{\alpha_1\in\ocA_1\cap\partial\ocA_0,\ldots,
\alpha_r\in\ocA_r\cap\partial\ocA_{r-1}} \exp(r
\beta{\varepsilon}) \prod_{i=1}^r \exp
\bigl(-\beta\max \bigl(0,H(\alpha_i)- H(a_i) \bigr) \bigr)
\\
% \leq
%|\cE|^r
%-
&&\qquad
\leq %\sum_{\alpha_1\in\ocA_1,\ldots,
\sum_{\alpha_1\in\ocA_1\cap\partial\ocA_0,\ldots,
\alpha_r\in\ocA_r\cap\partial\ocA_{r-1}} \exp(r\beta{\varepsilon})
\exp \bigl(-\beta \delta_1 \bigl| \bigl\{ i\leq r\dvtx  H(\alpha_i)>H(a_i)
\bigr\} \bigr| \bigr).
\end{eqnarray*}
For $1\leq i\leq r$, the point $\alpha_i$
belongs to $\partial\ocA_{i-1}$.
By Lemma \ref{exitval},
this implies
that
$H(\alpha_i)\neq H(a_i)$.
%$\height(\ocA_{i})\neq\height(\ocA_{i-1})$
Moreover there is no strictly decreasing sequence
of energy values of length larger than $p+2$
(recall that
$\delta_0<\delta_1<\cdots<\delta_p$ are the possible
values for the difference of the energies of two configurations of~$\cD$).
Therefore
\[
\bigl| \bigl\{ i\leq r\dvtx  H(\alpha_i)>H(a_i) %\height(\ocA_{i})> \height(\ocA_{i-1})
\bigr
\} \bigr| \geq \biggl\lfloor\frac{r}{p+2} \biggr\rfloor. % \geq \frac{r}{2p}-1
\]
We conclude that
\[
\P \bigl( \bpi(\cE)= (\ocA_0,\ocA_1,\ldots,
\ocA_r) \bigr) \leq |\cE|^r %\exp(r\beta{\varepsilon})
\exp \biggl( r
\beta{\varepsilon} %\frac{\beta\delta_1 r}{2p}+\beta\delta_1\Big)
-{\beta\delta_1} \biggl\lfloor
\frac{r}{p+2} \biggr\rfloor \biggr)
\]
and
\[
\P\bigl(R(\cE)=r\bigr) \leq \bigl|\ocM(\cE) \bigr|^r |\cE|^r
%e^{\beta\delta_1}
\exp \biggl( r\beta{\varepsilon} %\frac{\beta\delta_1 r}{2p}+\beta\delta_1\Big)
-{\beta
\delta_1} \biggl\lfloor\frac{r}{p+2} \biggr\rfloor \biggr).
\]
By Lemmas \ref{disjoint} and \ref{fondo}, the map which associates to
each maximal cycle compound
its bottom is
one to one, hence
$ |\ocM(\cE) |\leq
|\cE|$.
The hypothesis on $\cD$ yields that,
for $\varepsilon>0$ and for $\beta$ large enough,
\[
%|\cE|+
|\cE| \leq {v_\cD} |Q |^{v_\cD}
\leq \exp(\beta{\varepsilon}),
\]
whence
%Thus for $\beta$ large enough, we have
%
\[
\P\bigl(R(\cE)=r\bigr) \leq \exp \biggl( 3r\beta{\varepsilon}
-{\beta\delta_1} \biggl\lfloor\frac{r}{p+2} \biggr
\rfloor \biggr).
\]
Choosing $\varepsilon$ small enough and resumming this inequality, we
obtain the desired estimate.
\end{pf}
%
%Let $K>0$. There exists a value $D$ which depends
%only on $v_\cD$ and $K$ such that,
%for $\beta$ large enough,
%for any
%cycle compound
%$\ocA$
%included in $\cD$, we have
%$$\forall\alpha\in\ocA
%
%Let $\lambda_1<\cdots<\lambda_p$ be the energy values
%of the configurations in $\cC_d$, i.e.,
%$$\{
% =

We start now the proof of Theorem \ref{totcontrole}. We consider the
decomposition of $\cD$ into its maximal cycle compounds $\ocM(\cD)$ in
order to reduce the problem to the case where $\cD$ is a cycle
compound.
%control the size
%of the
%space--time clusters
%created inside such a cycle compound.
We decompose
\begin{eqnarray*}
&&
\P \bigl( \diam \STC\bigl(0, \tau(\cD) \bigr)\geq D \bigr) \\
&&\qquad\leq \P\bigl(R ( \cD)
\geq r\bigr)
\\
&&\qquad\quad{} + \sum_{0\leq k<r} \P \bigl( \diam \STC\bigl(0, \tau (
\cD) \bigr)\geq D, R ( \cD) =k \bigr).
\end{eqnarray*}
Let us fix $k<r$.
We write, using the notation defined before
Proposition \ref{long}, and
setting $\ocA_0=\bpi(\alpha,\cD)$,
\begin{eqnarray*}
&&\P \bigl( \diam \STC\bigl(0, \tau(\cD) \bigr)\geq D, R ( \cD) =k \bigr)
% \leq
%P_\alpha\bigg(
%$\diam \STC(0,
%)\geq D'$}
%(\ocA_1,\ldots,\ocA_k)$}
%}\bigg)
\\
&&\qquad
\leq \sum_{\ocA_1,\ldots,\ocA_k
\in\ocM(\cD)} %\in\ocM(\ocA\setminus\{ \eta \})}
\P \pmatrix{\mbox{$\displaystyle \sum_{0\leq j\leq k}\diam \STC(
\tau_j, \tau_{j+1})\geq D$} \vspace*{2pt}\cr \mbox{$\displaystyle \bpi(\cD)= (
\ocA_0,\ocA_1,\ldots,\ocA_k)$} }
% \leq
%P_\alpha\bigg(
%$\diam \STC(\tau_i,
%(\ocA_1,\ldots,\ocA_k)$}
%}\bigg)
\\
&&\qquad\leq \sum_{\ocA_1,\ldots,\ocA_k
\in\ocM(\cD)} %\in\ocM(\ocA\setminus\{ \eta \})}
\sum_{j=0}^{k} \sum
_{\alpha_j\in\ocA_j} \P \pmatrix{\mbox{$\displaystyle \diam \STC\bigl(
\s^{\alpha,n\pm}_{Q,t}, \tau_j\leq t\leq
\tau_{j+1}\bigr)\geq D/r$}\vspace*{2pt}\cr \mbox{$\displaystyle  \s^{\alpha,n\pm}_{Q,\btau_j}=
\alpha_j, \bpi(\cD)= (\ocA_0,\ocA_1,\ldots,
\ocA_k)$}}
\\
&&\qquad\leq \sum_{\ocA_1,\ldots,\ocA_k
\in\ocM(\cD)} %\in\ocM(\ocA\setminus\{ \eta \})}
\sum_{j=0}^{k}
\sum_{\alpha_j\in\ocA_j} %\P\big(
\P \bigl( \diam
\STC\bigl( \s^{\alpha_j,n\pm}_{Q,t}, 0\leq t\leq \tau(
\ocA_{j})\bigr)\geq D/r \bigr). %\\
% \leq
%0\leq t\leq
\end{eqnarray*}
Given a value $K$, we choose $r$ such that
%$$\delta_1-\frac{\delta_1r}{8p}<-2K',
$cr>2K$, where $c$ is the constant appearing in
Proposition \ref{long}.
%and $D'$ such that
%$D'/r>D_i(2K')$ where
%$D_i(2K')$ is the value given by the induction hypothesis associated
%to $2K'$.
We choose then $\varepsilon>0$ such that
$r\varepsilon< K$.
%$$\varepsilon<
%%and the integer $r$ so that
%$$old
%-\frac{\delta_1r}{4p} < -2K',
%and a value $D'$ such that
By Lemmas \ref{disjoint} and~\ref{fondo}, the map which associates to
each maximal cycle compound
its bottom is
one to one, hence
\[
\bigl| \ocM(\cD) \bigr| \leq |\cD| \leq \exp(\beta\varepsilon).
\]
The last inequality holds for
$\beta$ large,
thanks to the hypothesis on $\cD$.
%$|\cD|\leq\exp(\beta\varepsilon$.
Combining the previous estimates, we obtain, for $\beta$ large enough,
\begin{eqnarray*}
&&
\P \bigl( \diam \STC\bigl(0, \tau(\cD) \bigr)\geq D \bigr) \\
&&\qquad\leq
\frac{1}{c} %\exp(-\beta cr)
\exp(-2\beta K)
\\
%+ 3\varepsilon k
%-\frac{\delta_1k}{2p}\big)\Big)
&&\qquad\quad{} + r^2\exp(\beta r\varepsilon)
\mathop{\max_{\ocA\in\ocM(\cD)}}_{
\alpha\in\ocA} %\P\big(
\P \bigl( \diam \STC\bigl( \s^{\alpha,n\pm}_{Q,t}, 0\leq t\leq \tau(
\ocA)\bigr)\geq D/r \bigr). %\\
% \leq
%+
%.
\end{eqnarray*}
%
%|\cE|+
%maximal cycle compounds $\ocM(\cD)$ and
To conclude, we need to
control the size
of the
space--time clusters
created inside a cycle compound $\ocA$ included in $\cD$. More precisely,
we need to prove the statement of Theorem \ref{totcontrole}
for a cycle compound.
We shall prove the following
result by induction on the depth of the cycle
compound.\vspace*{9pt}

\textit{Induction hypothesis at step $i$}:
For any $K>0$, there exists $D_i$
depending only on $v_\cD$ and $K$ such that,
for $\beta$ large enough,
for any
cycle compound
$\ocA$
included in $\cD$ having depth less than or equal
to $\delta_i$,
\[
\forall\alpha\in\ocA\qquad \P \bigl( \diam \STC\bigl(
\s^{\alpha,n\pm}_{Q,t}, 0\leq t\leq\tau(\ocA)\bigr) \geq
D_i \bigr) \leq \exp(-\beta K).
\]
Once this result is proved,
to complete the proof of Theorem \ref{totcontrole}, we simply
choose $D$ such that
\[
\frac{D}{r} > \max \bigl\{ D_i(2K)\dvtx 0\leq i\leq p \bigr\},
\]
where $D_i(2K)$ is the constant associated to $2K$ in the induction hypothesis.
We proceed next to the inductive proof.
Suppose that $\ocA$ is a cycle compound of depth~$0$.
Then $\ocA=\{ \eta \}$ is a singleton and therefore
\[
\diam \STC\bigl(0,\tau(\ocA)\bigr) \leq \sum_{C\in\cC(\eta)}
\diam C+1 \leq v_\cD+1. %\max \{ \diam  C\dvtx  C\mbox{ connected component of }\eta \}
\]
Let $i\geq0$. Suppose that the result has been proved for all
the cycle compounds included in $\cD$ of depth less than or equal
to $\delta_i$.
%More precisely, we suppose the following result
%(the integer $n\in\{ 1,\ldots,d \}$ is fixed).
%
Let now $\ocA$ be a cycle compound of depth
$\delta_{i+1}$.
By Lemma \ref{fondo}
the bottom of $\ocA$ consists of a
unique\vspace*{1pt} configuration~$\eta$.
Let $\alpha\in\ocA$ be a starting configuration.
We study next the process
$(\s^{\alpha,n\pm}_{Q,t})_{t\geq0}$, and unless stated otherwise, the
$\STC$ and the quantities like
$\diam \STC$ are those associated to this process.
We define the time $\theta$ of the last visit to $\eta$
before the time $\tau(\ocA)$, that is,
\[
\theta = \sup\bigl\{ s\leq\tau(\ocA)\dvtx \sigma_{Q,s}=\eta \bigr\}
\]
[if the process does not visit $\eta$ before $\tau(\ocA)$,
then we take $\theta=0$].
%Remark that $\theta$ is not a stopping time.
Considering the random times
$\tau(\ocA\setminus\{ \eta \})$,
$\theta$ and $\tau(\ocA)$,
we have by Lemma \ref{triangle},
\begin{eqnarray*}
\diam \STC\bigl(0,\tau(\ocA)\bigr) &\leq& \diam \STC\bigl(0,\tau\bigl(\ocA
\setminus\{ \eta \}\bigr)\bigr)
\\
&&{} + \diam \STC\bigl(\tau\bigl(\ocA\setminus\{ \eta \}\bigr), \theta\bigr) + \diam
\STC\bigl(\theta,\tau(\ocA)\bigr).
\end{eqnarray*}
Indeed, if
$\tau(\ocA\setminus\{ \eta \})<\tau(\ocA)$, then
$\tau(\ocA\setminus\{ \eta \})\leq\theta\leq\tau(\ocA)$,
and the above inequality holds. Otherwise,
if $\tau(\ocA\setminus\{ \eta \})=\tau(\ocA)$, then
$\theta=0$ and the second term of the
right-hand side vanishes.
%With Lemma \ref{srcex}, we can also control the diameter of
%the space--time clusters between
%$\theta$ and $\tau(\ocA))$.
%Let $D>0$ and let $\alpha\in\ocA$. We write
Let $D>0$, and let us write
\begin{eqnarray*}
\P \bigl( \diam \STC\bigl(0,\tau(\ocA)\bigr) \geq D \bigr) &\leq& \P \bigl( \diam
\STC\bigl(0,\tau\bigl(\ocA\setminus\{ \eta \}\bigr)\bigr)\geq D/3 \bigr)
\\
&&{} + \P \bigl( \diam \STC\bigl(\tau\bigl(\ocA\setminus\{ \eta \}\bigr), \theta
\bigr)\geq D/3 \bigr)
\\
&&{} + \P \bigl( \diam \STC\bigl(\theta,\tau(\ocA)\bigr)\geq D/3 \bigr).
\end{eqnarray*}
We will now consider different starting points, hence we use
the more explicit notation for the $\STC$.
From the Markov property,
we have
%$$\tau(\theta(\eta,\tau(\ocA)),\ocA)=\tau(\ocA)=
%
\begin{eqnarray*}
&&
\P \bigl( \diam \STC\bigl( \s^{\alpha,n\pm}_{Q,t}, \tau\bigl(\ocA
\setminus\{ \eta \}\bigr)\leq t\leq \theta\bigr) \geq D/3 \bigr) %\\
\\
&&\qquad
\leq\P \bigl( \diam \STC\bigl( \s^{\eta,n\pm}_{Q,t}, 0\leq t\leq\theta
\bigr) \geq D/3 \bigr)
\end{eqnarray*}
and
\begin{eqnarray*}
&&
\P \bigl( \diam \STC\bigl( \s^{\alpha,n\pm}_{Q,t}, \theta\leq t\leq
\tau(\ocA)\bigr) \geq D/3 \bigr)
\\
&&\qquad
\leq \P \bigl( %\diam \STC(0,\tau(\ocA))
\diam \STC\bigl( \s^{\eta,n\pm}_{Q,t},
0\leq t\leq\tau(\ocA)\bigr) \geq D/3, \tau(\ocA)= \tau\bigl(\ocA\setminus\{ \eta
\}\bigr) \bigr)
\\
&&\qquad\leq \P \bigl( \diam \STC\bigl( \s^{\eta,n\pm}_{Q,t}, 0\leq t\leq
\tau\bigl(\ocA\setminus\{ \eta \}\bigr)\bigr) \geq D/3 \bigr),
\end{eqnarray*}
whence
\begin{eqnarray*}
&&\P \bigl( \diam \STC\bigl( \s^{\alpha,n\pm}_{Q,t}, 0\leq t\leq\tau(
\ocA)\bigr) \geq D \bigr)
\\
&&\qquad\leq2\sup_{\gamma\in\ocA}\P \bigl( \diam \STC\bigl( \s^{\gamma,n\pm}_{Q,t},
0\leq t\leq\tau\bigl(\ocA\setminus\{ \eta \}\bigr)\bigr)\geq D/3 \bigr)
\\
&&\qquad\quad{}+ \P \bigl( \diam \STC\bigl( \s^{\eta,n\pm}_{Q,t}, 0\leq t\leq
\theta\bigr) \geq D/3 \bigr).
\end{eqnarray*}
We first control the size of the
space--time clusters
created during an excursion
outside the bottom $\eta$.
%%
%
%le5.9 #&#
\begin{lemma}\label{srcex}
For any $K'>0$, there exists $D'$
depending only on $v_\cD,K'$ such that,
for $\beta$ large enough,
for any
$\alpha\in\ocA$,
\[
\P \bigl( \diam \STC\bigl( \s^{\alpha,n\pm}_{Q,t}, 0\leq
t\leq \tau\bigl(\ocA\setminus\{ \eta \}\bigr) \bigr)\geq D' \bigr)
\leq \exp\bigl(-\beta K'\bigr).
\]
\end{lemma}
\begin{pf}
The argument is very similar to the initial step of the proof
of Theorem \ref{totcontrole}, that is, we reduce the problem to the maximal
cycle compounds included in
$\ocA\setminus\{ \eta \}$. Although it is possible to include
these two steps
in a more general result, for the clarity of the exposition, we prefer to
repeat the argument rather than to introduce additional notations.
We consider the decomposition of
$\ocA\setminus\{ \eta \}$ into its
maximal cycle compounds $\ocM(\ocA\setminus\{ \eta \})$.
Each cycle compound of
$\ocM(\ocA\setminus\{ \eta \})$ has a depth strictly less
than $\delta_{i+1}$; hence we can apply the induction hypothesis
and control the size
of the
space--time clusters
created inside such a cycle compound.
We decompose next
\begin{eqnarray*}
&&
\P \bigl( \diam \STC\bigl(0, \tau\bigl(\ocA\setminus\{ \eta \}\bigr) \bigr)\geq
D' \bigr)\\
&&\qquad\leq \P\bigl(R \bigl( \ocA\setminus\{ \eta \}\bigr) \geq
r\bigr)
\\
&&\qquad\quad{}+ \sum_{0\leq k<r} \P \bigl( \diam \STC\bigl(0, \tau
\bigl( \ocA\setminus\{ \eta \}\bigr) \bigr)\geq D', R \bigl( \ocA
\setminus\{ \eta \}\bigr) =k \bigr).
\end{eqnarray*}
Let us fix $k<r$ and, denoting
simply $\ocM=\ocM(\ocA\setminus\{ \eta \})$,
we write, using the notation defined before
Proposition \ref{long},
and setting $\ocA_0=\bpi(\alpha,
\ocA\setminus\{ \eta \}
)$,
\begin{eqnarray*}
\hspace*{-4pt}&&
\P \bigl( \diam \STC\bigl(0, \tau\bigl(\ocA\setminus\{ \eta \}\bigr) \bigr)\geq
D', R \bigl( \ocA\setminus\{ \eta \}\bigr) =k \bigr) %\\
% \leq
%P_\alpha\bigg(
%$\diam \STC(0,
%)\geq D'$}
%(\ocA_1,\ldots,\ocA_k)$}
%}\bigg)
\\
\hspace*{-4pt}&&\qquad\leq \sum_{\ocA_1,\ldots,\ocA_k
\in\ocM} %\in\ocM(\ocA\setminus\{ \eta \})}
\P \pmatrix{
\mbox{$\displaystyle \sum_{0\leq j\leq k}\diam \STC(
\tau_j, \tau_{j+1})\geq D'$}\vspace*{2pt}\cr \mbox{$\displaystyle \bpi
\bigl(\ocA\setminus\{ \eta \}\bigr)= (\ocA_0,\ocA_1,
\ldots,\ocA_k)$}} %\\
% \leq
%P_\alpha\bigg(
%$\diam \STC(\tau_i,
%(\ocA_1,\ldots,\ocA_k)$}
%}\bigg)
\\
\hspace*{-4pt}&&\qquad\leq \sum_{\ocA_1,\ldots,\ocA_k
\in\ocM} %\in\ocM(\ocA\setminus\{ \eta \})}
\sum_{0\leq j\leq k} \sum
_{\alpha_j\in\ocA_j} \P \pmatrix{\mbox{$\displaystyle  \diam \STC\bigl(
\s^{\alpha,n\pm}_{Q,t}, \tau_j\leq t\leq
\tau_{j+1}\bigr)\geq D'/r$}\vspace*{2pt}\cr
\mbox{$\displaystyle  \s^{\alpha,n\pm}_{Q,\btau_j}=
\alpha_j, %\sigma_{Q,\btau_j}=\alpha_j,
\bpi\bigl(\ocA\setminus\{ \eta \}\bigr)= (
\ocA_0,\ocA_1,\ldots,\ocA_k)$}}
\\
\hspace*{-4pt}&&\qquad\leq \sum_{\ocA_1,\ldots,\ocA_k
\in\ocM} %\in\ocM(\ocA\setminus\{ \eta \})}
\sum_{0\leq j\leq k} \sum
_{\alpha_j\in\ocA_j} %\P\big(
\P \bigl( \diam \STC\bigl(
\s^{\alpha_j,n\pm}_{Q,t}, 0\leq t\leq \tau(\ocA_{j})\bigr)
\geq D'/r \bigr). %\\
% \leq
%0\leq t\leq
\end{eqnarray*}
Given a value $K'$, we choose $r$ such that
%$$\delta_1-\frac{\delta_1r}{8p}<-2K',
$cr>2K'$, where $c$ is the constant appearing in
Proposition \ref{long}
and $D'$ such that
$D'/r>D_i(2K')$ where
$D_i(2K')$ is the value given by the induction hypothesis at step $i$ associated
to $2K'$.
Notice that this value is uniform with respect to
the cycle compound $\ocA\subset\cD$ of depth $\delta_{i+1}$
because all the cycle compounds of $\ocM$ are included in $\cD$
and have a depth
at most equal to $\delta_i$.
We choose then $\varepsilon>0$ such that
$r\varepsilon< K'$.
%$$\varepsilon<
%%and the integer $r$ so that
%$$old
%-\frac{\delta_1r}{4p} < -2K',
%and a value $D'$ such that
By Lemmas \ref{disjoint} and \ref{fondo}, the map which associates to
each maximal cycle compound
its bottom is
one to one, hence
\[
\bigl| \ocM\bigl(\ocA\setminus\{ \eta \}\bigr) \bigr| \leq \bigl| \ocA\setminus\{ \eta \} \bigr|
\leq |\cD| \leq \exp(\beta\varepsilon).
\]
The last inequality holds for
$\beta$ large,
thanks to the hypothesis on $\cD$.
%$|\cD|\leq\exp(\beta\varepsilon$.
Combining the previous estimates, we obtain, for $\beta$ large enough,
\begin{eqnarray*}
&&
\P \bigl( \diam \STC\bigl(0, \tau\bigl(\ocA\setminus\{ \eta \}\bigr) \bigr)\geq
D' \bigr)
\\
&&\qquad
\leq\bigl|\ocM\bigl(\ocA\setminus\{ \eta \}\bigr) \bigr|^{r-1}r^2 \bigl|\ocA
\setminus\{ \eta \} \bigr| \exp\bigl(-2\beta K'\bigr) + %\exp-\frac{\beta\delta_1r}{4p}.
\frac{1}{c}\exp(-\beta cr) %\sum_{k\geq r}
%+ 3\varepsilon k
%-\frac{\delta_1k}{2p}\big)\Big)
\\
&&\qquad
\leq r^2\exp\bigl(\beta\bigl(r\varepsilon-2 K'\bigr)
\bigr) + \frac{1}{c} \exp\bigl(-2\beta K'\bigr).
% \leq
%+
%.
\end{eqnarray*}
%
%|\cE|+
The last quantity is less than
$\exp(-\beta K')$ for $\beta$ large enough.
\end{pf}

The remaining task is to control
the space--time clusters between
$\tau(\ocA\setminus\{ \eta \})$ and~$\theta$, which amounts to control
%$$\diam \STC(\tau(\ocA\setminus\{ \eta \}),
%
\[
\P \bigl( \diam \STC\bigl( \s^{\eta,n\pm}_{Q,t}, 0\leq t\leq\theta
\bigr) \geq D/3 \bigr).
\]
%
%In the remaining of this proof, we write simply
%$\theta$ instead of $\theta(\eta,\tau(\ocA))$.
We suppose that
$\tau(\ocA\setminus\{ \eta \})<\tau(\ocA)$ (otherwise $\theta=0$)
and that the process is in $\eta$ at time $0$.
To the continuous-time trajectory
$ (
\s^{\eta,n\pm}_{Q,t},
0\leq t\leq
\theta )$,
we associate a discrete path $\o$ as follows:
\[
\matrix{ T_0=0,&&\o_0=\sigma_{Q,0}=\eta,
\vspace*{2pt}\cr
T_1=\min \{ t>T_0\dvtx \sigma_{Q,t}\neq
\o_0 \}, && \o_1=\sigma_{Q,T_1},
\vspace*{2pt}\cr
T_2=\min \{ t>T_1\dvtx \sigma_{Q,t}\neq
\o_1 \}, && \o_2=\sigma_{Q,T_2},
\vspace*{2pt}\cr
\vdots&&
\vdots
\vspace*{2pt}\cr
T_{k}=\min \{ t>T_{k-1}\dvtx \sigma_{Q,t}\neq
\o_{k-1} \},&& \o_k=\sigma_{Q,T_{k}},
\vspace*{2pt}\cr
\vdots&&
\vdots
\vspace*{2pt}\cr
T_{S-1}=\min \{ t>T_{S-2}\dvtx \sigma_{Q,t}\neq
\o_{S-2} \},&& \o_{S-1}=\sigma_{Q,T_{S-1}},
\vspace*{2pt}\cr
%T_S=\tau(\ocA\setminus\{ \eta \}),
T_S=\theta, &&\o_{S}=
\sigma_{Q,T_{S}}=\eta.}
\]
Let $R$ be the number of visits of the path $\o$ to $\eta$, that is,
\[
R = \bigl| \{ 1\leq i\leq S\dvtx \o_i=\eta \} \bigr|.
\]
We define then the indices
$\phi(0),\ldots,\phi(R)$
of the successive visits
to $\eta$ by setting $\phi(0)=0$ and for $i\geq1$,
\[
\phi(i) = \min \bigl\{ k\dvtx  k>\phi(i-1), \o_k=\eta \bigr\}.
\]
The times
$\tau_0,\ldots,\tau_R$
corresponding to these indices are
\[
\tau_i=T_{\phi(i)},\qquad 0\leq i\leq R.
\]
Each
subpath
\[
\ot^i=\bigl(\o_k,\phi(i)\leq k\leq\phi(i+1)\bigr)
\]
is an excursion outside $\eta$ inside $\ocA$.
We denote by
$\cC(\eta)$ the connected components of $\eta$.
Let $C$ belong to $\cC(\eta)$.
By Lemma \ref{stccc}, the space--time sets
$C\times\{\phi(i)\}$ and
$C\times\{\phi(i+1)\}$ belong to the same space--time cluster of $\ot^i$;
therefore they are also in the same space--time cluster of
$\STC(\tau_i,\tau_{i+1})$.
Thus the space--time set
\[
C\times\{ \tau_0,\ldots,\tau_R \}
\]
belongs to one space--time cluster of $\STC(0,\theta)$.
The following computations deal with the process
$(\s^{\eta,n\pm}_{Q,t})_{t\geq0}$ starting from
$\eta$ at time $0$. Hence all the $\STC$ and the exit times
are those associated to this process.
Let $\cC$ belong to $\STC(0,\theta)$.
%From the previous construction,
We consider two cases:

$\bullet$
If
$\cC\cap (\eta\times\{ \tau_0,\ldots,\tau_R \} ) =
\varnothing$,
then there exists $i\in\{ 0,\ldots, R-1 \}$ such that
\[
\cC\in\STC(\tau_i,\tau_{i+1}),\qquad \cC\cap\bigl(\eta\times\{
\tau_i,\tau_{i+1} \}\bigr)=\varnothing.
\]
Therefore
\[
\diam \cC \leq %\max_{0\leq i<R}
\mathop{\max_{\cC\in\STC(\tau_i,\tau_{i+1})}}_{
\cC\cap(Q\times\{ \tau_i,\tau_{i+1} \})=\varnothing}
\diam \cC.
\]

$\bullet $
If
$\cC\cap (\eta\times\{ \tau_0,\ldots,\tau_R \}
) \neq
\varnothing$, then
there exists a connected component
$C\in\cC(\eta)$ and $i\in\{ 0,\ldots,R \}$ such that
$\cC\cap (C\times\{ \tau_i \} ) \neq
\varnothing$.
From the previous discussion, we conclude that
$C\times\{ \tau_0,\ldots,\tau_R \}$ is included in $\cC$.
In fact,
for any $C$ in $\cC(\eta)$, we have
\[
\mbox{either}\quad \cC \cap \bigl(C\times\{ \tau_0,\ldots,
\tau_R \} \bigr) = \varnothing \quad\mbox{or}\quad C\times\{
\tau_0,\ldots,\tau_R \} \subset \cC.
\]
For $C$ in $\cC(\eta)$ and
$i\in\{ 0,\ldots,R-1 \}$,
we denote by
$\STC(\tau_i,\tau_{i+1})(C)$
the space--time cluster of
$\STC(\tau_i,\tau_{i+1})$
containing $C\times\{ \tau_i,\tau_{i+1} \}$.
The space--time cluster $\cC$ is thus included in the set
\[
\mathop{\bigcup_{ C\in\cC(\eta)}}_
{ C\times\{ 0, \theta \}\subset\cC} \bigcup
_{0\leq i<R} \STC(\tau_i,\tau_{i+1})
(C).
\]
For any
$C\in\cC(\eta)$, the space--time set
\[
\bigcup_{0\leq i<R} \STC(\tau_i,
\tau_{i+1}) (C)
\]
is connected, and its diameter is bounded by
\[
2\max_{0\leq i <R} \diam \STC(\tau_i,
\tau_{i+1}) (C).
\]
The factor $2$ is due to the fact that the two sites realizing the diameter
might belong to two
different excursions outside $\eta$.
Therefore
\[
\diam \cC \leq \mathop{\sum_{C\in\cC(\eta)}}_
{ C\times\{ 0, \theta \}\subset\cC}
2\max_{0\leq i <R} \diam \STC(\tau_i,
\tau_{i+1}) (C).
\]
From the inequality obtained in the first case,
we conclude that
\[
\mathop{\max_{\cC\in\STC(0,\theta)}}_{
\cC\cap(Q\times\{ 0, \theta \})=\varnothing} \diam \cC \leq \max
_{0\leq i<R} \mathop{\max_{\cC\in\STC(\tau_i,\tau_{i+1})}}_{
\cC\cap(Q\times\{ \tau_i,\tau_{i+1} \})=\varnothing}
\diam \cC.
\]
We sum next the inequality of the second case over all the
elements of $\STC(0,\theta)$ intersecting $Q\times\{ 0, \theta \}$.
Since
two distinct $\STC$ of
$\STC(0,\theta)$
do not intersect at time~$0$, they do not meet the same connected
components of $\eta$,
and
we obtain
\[
\mathop{\sum_{\cC\in\STC(0,\theta)}}_{
\cC\cap(Q\times\{ 0,\theta \})\neq\varnothing} \diam \cC \leq
\sum_{C\in\cC(\eta)} 2\max_{0\leq i <R} \diam \STC(
\tau_i,\tau_{i+1}) (C).
\]
Putting together the two previous inequalities, we conclude that
\[
\diam \STC(0, \theta) \leq 2|\eta| \max_{0\leq i <R} \diam \STC(
\tau_i,\tau_{i+1}).
\]
%
%Moreover, for any
%${C\in\cC(\eta)}$,
%$${\diam \STC(0, \theta)(C) \leq
%2\rad(0, \theta)(C)
%+\diam C\cr
% \leq
%2\max_{0\leq i <R}
%+\diam C\cr
% \leq
%2\max_{0\leq i <R}
%+\diam C,}$$
%therefore
%$${\diam \STC(0, \theta)
% \leq
%2\max_{0\leq i <R}
%.}$$
We write
\begin{eqnarray*}
&&
\P \bigl( \diam \STC( 0,\theta) \geq D/3 \bigr)
\\
&&\qquad
\leq\P (R\geq r ) +\sum_{0\leq k< r} \P \bigl( \diam \STC( 0,
\theta) \geq D/3, R=k \bigr).
\end{eqnarray*}
For a fixed integer $k$, the previous inequalities and the
Markov property yield
\begin{eqnarray*}
&&
\P \bigl( \diam \STC( 0,\theta) \geq D/3, R=k \bigr)
\\
&&\qquad
\leq \P \Bigl( 2|\eta| \max_{0\leq i <k} \diam \STC(
\tau_i,\tau_{i+1}) \geq D/3, R=k \Bigr)
\\
&&\qquad
\leq k \P \bigl( 2|\eta| \diam \STC( % \s^{\eta,n\pm}_{Q,t},
0,\tau_1)
\geq D/3, \tau_1< \tau(\ocA) %\tau(\ocA\setminus\{ \eta \}))
\bigr).
\end{eqnarray*}
Recalling that
\[
T_1=\min \{ t>T_0\dvtx \sigma_{Q,t}\neq\eta \},\qquad
\tau_1=\min \{ t>T_1\dvtx \sigma_{Q,t}=\eta \},
\]
we claim that, on the event $\tau_1
<\tau(\ocA)$, we have
\[
\diam \STC(0,\tau_{1}) \leq \diam \STC(T_1,
\tau_{1})+1.
\]
Indeed, let $\cC$ belong to
$\STC(0,\tau_{1})$. If $\cC$ is in
$\STC(T_1,\tau_{1})$, then obviously
\[
\diam \cC \leq \diam \STC(T_1,\tau_{1}).
\]
Otherwise, the set
$\cC\cap(Q\times[T_1,\tau_1])$ is the union of several elements of
$\STC(T_1,\tau_{1})$, say $\cC_1,\ldots,\cC_r$, which all intersect
$Q\times\{ T_1 \}$.
The spin flip leading from $\eta$ to
$\sigma_{Q,T_1}$ can change only by one the sum of the diameters of
the $\STC$ present
at time $0$. This spin flip occurred in $\cC$ if and only if
\[
\cC\cap\bigl(Q\times\{ 0 \}\bigr)\neq\cC\cap\bigl(Q\times\{ T_1 \}
\bigr),
\]
thus
\[
\diam \cC \leq \sum_{1\leq i\leq r} \diam \cC_i +
1_{
\cC\cap(Q\times\{ 0 \})\neq\cC\cap(Q\times\{ T_1 \})
}.
\]
Summing over all the elements of $\STC(0,\tau_1)$
which intersect $Q\times\{ 0 \}$, we obtain
the desired inequality.
Reporting in\vspace*{1pt} the previous computation and conditioning with respect
to
$\s^{\eta,n\pm}_{Q,T_1}$,
we get
\begin{eqnarray*}
&&
\P \bigl( \diam \STC( 0,\theta) \geq D/3, R=k \bigr)
\\
&&\qquad\leq k \P \biggl( \diam \STC( % \s^{\eta,n\pm}_{Q,t},
T_1,\tau_1)
\geq \frac{D}{6 |\eta|}-1, \tau_1< \tau(\ocA)
\biggr)
\\
&&\qquad\leq \sum_{\gamma\in\ocA\setminus\{ \eta \}} k \P \biggl( \s^{\eta,n\pm}_{Q,T_1}=
\gamma, \diam \STC( % \s^{\eta,n\pm}_{Q,t},
T_1,\tau_1)
\geq \frac{D}{6 |\eta|}-1, \tau_1< \tau(\ocA)
\biggr)
\\
&&\qquad
\leq |\ocA |k \max_{\gamma\in
\ocA} %\ocA\setminus\{ \eta \}}
P \biggl( \diam \STC
\bigl( \s^{\gamma,n\pm}_{Q,t}, 0\leq t\leq \tau\bigl(\ocA\setminus\{
\eta \}\bigr)\bigr) \geq \frac{D}{6 |\eta|}-1 \biggr).
\end{eqnarray*}
Summing over $k$, we arrive at
\begin{eqnarray*}
&&
\P \bigl( \diam \STC( 0,\theta) \geq D/3 \bigr) \\
&&\qquad\leq \P (R\geq r )
\\
&&\qquad\quad{} + r^2 |\ocA | \max_{\gamma\in
\ocA} \P \biggl( \diam \STC
\bigl( \s^{\gamma,n\pm}_{Q,t}, 0\leq t\leq\tau\bigl(\ocA\setminus\{
\eta \}\bigr)\bigr) \geq \frac{D}{6 |\eta|}-1 \biggr).
\end{eqnarray*}
By the Markov property,
the variable $R$ satisfies for any $n,m\geq0$,
\begin{eqnarray*}
\P (R\geq n+m) &=& \P \bigl(\phi(n+m)<\tau(\ocA) \bigr)
\\
&=& \P \bigl(\phi(n)<\tau(\ocA), \phi(n+m)<\tau(\ocA) \bigr)
\\
&=& \P \bigl(\phi(n)<\tau(\ocA) \bigr) \P \bigl(\phi(m)<\tau(\ocA) \bigr)
\\
&=& \P (R\geq n) \P (R\geq m).
\end{eqnarray*}
Therefore the law of $R$ is the discrete geometric distribution and
\[
\forall n\geq0\qquad \P (R\geq n) = \biggl( \frac{E(R)} {
1+E(R)} \biggr)^n
\leq \exp-\frac{n}{1+E(R)}.
\]
By Corollary \ref{exitcom},
or more precisely its discrete-time counterpart,
%for $\varepsilon>0$ and $\beta$ large enough,
for $\beta$ large enough,
\[
{E(R)} \leq %{E( \tau(\ocA))}  \leq
\exp \bigl(
\tfrac{3}{2}\beta\depth(\ocA) \bigr) \leq \exp(2\beta\delta_{i+1})-1.
\]
Choosing
\[
r=\beta^2 \exp(2\beta\delta_{i+1}),
\]
we obtain
from the previous inequalities that
\begin{eqnarray*}
&&
\P \bigl( \diam \STC( 0,\theta) \geq D/3 \bigr)\\
&&\qquad\leq \exp-\beta^2 \\
&&\qquad\quad{}+
\beta^4 \exp(4\beta\delta_{i+1}) % \times
|\ocA |
\\
&&\hspace*{10.5pt}\qquad\quad{}\times\max_{\gamma\in
\ocA} \P \biggl( \diam \STC\bigl( \s^{\gamma,n\pm}_{Q,t},
0\leq t\leq\tau\bigl(\ocA\setminus\{ \eta \}\bigr)\bigr) %\geq(D/6v_\cD)\big)
\geq
\frac{D}{6v_\cD}-1 \biggr).
\end{eqnarray*}
We complete now the induction step at rank $i+1$.
Let $K>0$ be given. Let $K'>0$ be such that
$4\delta_{i+1}-K'<-3K$, and let $D'$ associated to $K'$
as in Lemma \ref{srcex}. Let $D''$ be such that
\[
\frac{D''}{6v_\cD}-1> D',\qquad \frac{D''}{3}> D'.
\]
%
%For any $D>D_{i+1}$,
Thanks to the hypothesis on $\cD$ and $Q$, for $\beta$ large enough,
\[
|\ocA | \leq |\cD | \leq \exp(\beta K).
\]
From the previous computation, we have
\[
\P \bigl( \diam \STC\bigl( \s^{\eta,n\pm}_{Q,t}, 0\leq t\leq\theta
\bigr) \geq D''/3 \bigr) \leq \exp-\beta^2
+\beta^4 \exp(-2\beta K).
\]
Since $D''/3>D'$, we have also
for any
$\gamma\in\ocA$,
\[
\P \bigl( \diam \STC\bigl( \s^{\gamma,n\pm}_{Q,t}, 0\leq t \leq\tau
\bigl(\ocA\setminus\{ \eta \}\bigr)\bigr)\geq D''/3
\bigr) \leq \exp(-3\beta K).
\]
Substituting the previous inequalities into the inequality obtained
before Lem\-ma~\ref{srcex},
we conclude that, for any
$\alpha\in\ocA$,
\[
\P \bigl( \diam \STC\bigl( \s^{\alpha,n\pm}_{Q,t}, 0\leq t\leq\tau(
\ocA)\bigr) \geq D'' \bigr) \leq %\exp-\beta^2
\bigl(
\beta^4 +3\bigr)\exp(-2\beta K)
\]
and the induction is completed.

%s6 #&#
\section{The metastable regime}
\label{metare}
%{\bf Sandwich boundary conditions.}
The goal of this section is to prove Theorem \ref{T2}, which states
roughly the following.
Under an appropriate
hypothesis on the initial law
and on the initial $\STC$,
for any $\kappa<\kappa_d$,
%Let
%$\tau_\beta$ be a time satisfying
%$$\limsup_{\b\to\infty}
the probability that
%for the process
%$(\s^{n\pm,\xi}_{\Sigma,t})_{t\geq0}$,
%there exists
%a $\STC$ in
%$\STC_\xi(0,\tb)$
a space--time cluster
of diameter
larger than
$\exp({\b L_d})$
is created before time $\exp(\beta\kappa)$
is $\SES$.
The hypothesis is satisfied by the law of a typical configuration
in the metastable regime.
This result allows us to control the speed of propagation of large supercritical
droplets.
As already pointed out by
Dehghanpour and Schonmann,
the control of this speed is a crucial point for the study of metastability
in infinite volume.
This estimate is quite delicate, and it is performed by induction over the
dimension. More precisely, we consider a set of the form
\[
\L^{n}\bigl(\exp(\beta L)\bigr) \times\L^{d-n}(\ln\beta)
\]
with $n\pm$ boundary conditions and we do the proof by induction over $n$.
The process in this set and with these boundary conditions behaves roughly
like the process in dimension $n$.
Proposition \ref{PT2} handles the case $n=0$.
A difficult point is that the growth of the supercritical droplet is more
complicated than a simple growth process. Indeed, supercritical
droplets might
be helped when they touch some clusters of pluses,
which were created independently.
Therefore we cannot proceed as in the simpler growth model
handled in \cite{CM2}.
To tackle this problem, we introduce an hypothesis on the initial law and
on the initial space--time clusters. The hypothesis on the initial law
guarantees that regions which are sufficiently far away are decoupled.
The hypothesis
on the initial space--time clusters provides a control on the
space--time clusters initially present in the configuration.
The point is that these two hypotheses are satisfied by the law of the
process in a fixed good
region until the arrival of the first supercritical
droplets.

The key ingredient in this part of the proof
is the lower bound on the time needed to cross parallelepipeds of
the above kind.
Heuristically, we will take into account the effect of the growing supercritical
droplet\vspace*{1pt} by using suitable boundary conditions, that is,
by using the Hamiltonian $H^{n\pm}$ instead of $H^-$.
Moreover, at the time when the configuration in the parallelepiped
starts to feel this effect,
it is rather likely that the parallelepiped is not void,
so that we have to consider more general initial configurations.

In any fixed $n$-small parallelepiped, it is very unlikely that
nucleation occurs before $\tb$, or that a large space--time cluster is
created before nucleation. However, the region under study contains
an exponential number of $n$-small parallelepipeds. Thus the previous
events will occur somewhere.
In Proposition \ref{nuclei}, we show that these events occur in at
most $\ln\ln\beta$
places.
%$n$-small bowes.
The proof uses the hypothesis
on the initial law and a simple counting argument.
The proof of Theorem \ref{T2} relies on a notion already used in bootstrap
percolation, namely boxes crossed by a space--time cluster; see
Definition \ref{cro}.
An $n$-dimensional box $\Phi$ is said to be crossed by a $\STC$
before time $t$
if,
for the dynamics restricted to
$\Phi\times
\L^{d-n}(\ln\beta)$,
%with initial configuration $\xi$ and
%$n\pm$ boundary condition,
there exists a space--time cluster
whose projection on the first
$n$ coordinates intersects two opposite faces of $\Phi$.
%The box is said to be crossed vertically in case the opposite faces
%which are touched
The point is that, if a box is crossed by a space--time cluster in some
time interval, then it is also crossed in the dynamics restricted to the
box with appropriate boundary conditions.
These appropriate boundary conditions are obtained as follows.
We put $n\pm$
boundary conditions on the restricted box exactly as on the large box,
and we
put $+$
boundary conditions on the faces which are normal to the direction
which is crossed.
%Let us now sketch the proof of
%Theorem \ref{T2}.
The induction step is long, and it is decomposed in eleven steps.

We will use the notation defined in Sections~\ref{energyestmates}
and \ref{stc}.
%mainly use the supremum norm
%$$\forall x=(x_1,\ldots,x_d)\in\Z^d
%|x|_\infty =
%%Given a site $x\in{\mathbb Z}^d$ we denote by $(x_1,\ldots,x_d)$ its
%%coordinates and given a set $C\subset{\mathbb Z}^d$
%We denote by
%$\dinf$ the distance associated to the supremum norm
%and we define the
%%$\dinf$
%diameter
%$\diam C$
%of a subset $C$ of ${\mathbb Z}^d$ by
%$$
% \big\}.$$
%Thus $\diam C$ is the sidelength of the minimal cube
%surrounding $C$.
Our main objective is to control the maximal diameter of
the $\STC$ created in a finite volume before the relaxation time.
Let $d\geq1$, let $n\in\{ 0,\ldots, d \}$ and
let us
consider
a parallelepiped $\Sigma$ in $\Z^{d}$
of the form
\[
\Sigma = \L^{n}(L_\beta) \times\L^{d-n}(\ln\beta),
\]
where
$\L^{n}(L_\beta)$ is a $n$-dimensional cubic box of side length
$L_\beta$,
$\L^{d-n}(\ln\beta)$ is a $d-n$-dimensional cubic box
of side length
$\ln\beta$ and the length $L_\beta$ satisfies
%, $h_\beta$ satisfy
%
\[
L_\beta \geq \ln\beta,\qquad %\lim_{\b\to\infty}
%l_\beta = +\infty,
\limsup
_{\b\to\infty} \frac{1}{\beta} \ln L_\beta < +\infty.
\]
%
%and basis $\L_\beta^{d-1}$, which is a translate of
%$$\Sigma = \L_\beta^{d-1} \times\{ 0,1 \}.$$
%Let $R$ be a cylinder with basis is a $(d-1)$-dimensional cubic box $
%$$R = \L^{d-1}\times\{ 0,\ldots,L \}.$$
%We call {\emph{floor}} of $R$
%its bottom face
%$\L^{d-1}\times\{ 0 \}$
%and {\emph{ceiling}} of $R$ its top face
%$\L^{d-1}\times\{ L \}$.
We set $\k_0=L_0=\Ga_0=0$, and for $n\geq1$
\[
\k_n = \frac{1}{n+1} (\G_1+\cdots+\G_n),\qquad
L_n = \frac{\Gamma_n-\kappa_{n}}{n}.
\]
%
%For $n\in\{ 0,\ldots,d \}$,
In the sequel we consider a time $\tb$ satisfying
\[
\limsup_{\b\to\infty} \frac{1}{\beta} \ln \tau_\beta <
\k_n.
\]
We say that a probability $\P(\cdot)$ is super-exponentially small
in $\beta$ (written in short $\SES$) if it satisfies
\[
\lim_{\beta\to\infty} \frac{1}{\beta}\ln \P(\cdot) = -\infty.
\]
%
%$\kappa<\kappa_n$ and we define
%$$\tau_\beta = \exp(\beta\k).$$
%s6.1 #&#
\subsection{Initial law}
We estimate the speed of growth of exponentially large droplets by
bounding from below the time needed by a large droplet to cross some tiles.
In each tile, we use
$n\pm$ boundary conditions in order to take into account the effect of
the droplet.
A major difficulty is to control the configuration until the arrival of the
supercritical droplets.
We introduce
an adequate hypothesis on the initial law
describing the configuration into the tile when the droplet enters.
This is achieved with the help of the following definitions.\vspace*{9pt}

\textit{$n$-small parallelepipeds}.
Let $n\geq1$.
A parallelepiped is $n$-small if all its sides have a length
larger than $\ln\ln\beta$ and smaller than $n\ln\beta$.
A parallelepiped is $0$-small
if all its sides have a length
larger than $\ln\ln\beta$ and smaller than $2\ln\ln\beta$.\vspace*{9pt}

\textit{Restricted ensemble}.
Let $n\geq0$.
We denote by $m_n$ the volume of the $n$-dimensional critical droplet.
%and by $\Ga_n$ its energy.
Let $Q$ be an $n$-small
%translate of the
parallelepiped.
%$$Q_n = \L^n(n\ln\beta)
%included in $\Sigma$.
The restricted ensemble ${\cal R}_n(Q)$ is the set of the
configurations $\sigma$ in $Q$
such that
$|\sigma| \leq m_n$
and $H_Q^{n\pm}(\sigma)\leq\Ga_n$,
that is,
\[
{\cal R}_n(Q) = \bigl\{ \sigma\in \{ -1,+1 \}^Q\dvtx |\sigma|
\leq m_n, H_Q^{n\pm}(\sigma)\leq\Ga_n
\bigr\}.
\]
%
%$$|\sigma| \leq m_n,
%H^{n\pm}_Q(\sigma) \leq \Ga_n.$$
%With this definition, we can check that
We observe that
$\cR_n(Q)$ is a cycle compound and that
\[
E \bigl(\cR_n(Q), \{ -1,+1 \}^Q\setminus
\cR_n(Q) \bigr) = \Ga_n.
\]
%
%{\bf Restricted ensemble}.
%Let $n\geq0$.
%%and by $\Ga_n$ its energy.
%Let $Q$ be an $n$-small
%%translate of the
%parallelepiped.
%%$$Q_n = \L^n(n\ln\beta)
%%\times\L^{d-n}(\ln\beta)$$
%%included in $\Sigma$.
%The restricted ensemble ${\cal R}_n(Q)$ is the set of the
%configurations $\sigma$ in $Q$
%such that
%$|\sigma| \leq m_n$ and
%$H_Q^{n\pm}(\sigma)\leq\Ga_n$, where
%$m_n$ is the volume of the $n$ dimensional critical droplet
%and $\Ga_n$ is its energy.
%H^{n\pm}_Q(\sigma) \leq \Ga_n.$$
%
Notice that the restricted ensemble satisfies the
hypothesis on the domain $\cD$ stated at the beginning of
Section~\ref{diamostc}.
We introduce next the hypothesis on the initial law, which is preserved until
the arrival of the supercritical droplets and which allows us to
perform the induction.\vspace*{9pt}

\textit{Hypothesis on the initial law at rank $n$.}
At rank $n=0$ we simply assume that the initial law $\mu$ is the Dirac mass
on the configuration
equal to $-1$ everywhere on $\Sigma$.
At rank $n\geq1$,
we will work with an initial law
$\mu$ on the configurations in $\Sigma$ satisfying
the following condition.
For any family $(Q_i,i\in I)$ of $n$-small parallelepipeds
included in $\Sigma$ such that
%
%$\bullet$ the sides of any parallelepiped of the family are larger
%than $\ln\ln\beta$ and less
%%than $n\ln\beta$,
%
%$\bullet$
two parallelepipeds of the family are at distance larger than
\[
5(d-n+1)\ln\ln\beta,
\]
we have the following estimates:
for any family of configurations
$(\sigma_i,i\in I)$ in the parallelepipeds $(Q_i,i\in I)$,
\[
\mu ( \forall i\in I, \sigma|_{Q_i}=\sigma_i ) \leq
\prod_{i\in I} \bigl(\phi_n(
\beta) \rho_{Q_i}^{n\pm}(\sigma_i) \bigr),
\]
where
\[
\rho_{Q_i}^{n\pm}(\sigma_i) = %
\cases{ \exp \bigl(-\beta H^{n\pm}_{Q_i}(
\sigma_i) %-H^{n\pm}_{Q_i}(\underline{\minus})
\bigr), &\quad if $\sigma_i\in{\cal
R}_n(Q_i)$,
\vspace*{2pt}\cr
\exp(-\beta\Ga_n), &\quad if $\sigma_i\notin{\cal
R}_n(Q_i)$,} % \mbox{ otherwise }
\]
and
$\phi_n(\beta)$
is a function depending only upon $\beta$
which is $\exp o(\beta)$, meaning that
\[
\lim_{\beta\to\infty} \frac{1} {
\beta} \ln\phi_n(\beta) =
0.
\]

\vspace*{9pt}

\textit{Hypothesis on the initial $\STC$ at rank $n$}.
We take also into account the presence of $\STC$ in the initial
configuration $\xi$. These
$\STC$ are unions of clusters of pluses present in $\xi$, we denote
them by
$\STC(\xi)$.
We suppose that for
any
$n$-small parallelepiped $Q$ included in $\Sigma$,
\[
\mathop{\sum_{\cC\in\STC(\xi)}}_
{\cC\cap Q\neq\varnothing}\diam\cC \leq
(d-n+1)\ln\ln\beta.
\]
%
%s6.2 #&#
\subsection{Lower bound on the nucleation time}
\label{finitebox}
In this section we give a
lower bound on the nucleation time in a finite box.
The proof rests on a coupling with the dynamics conditioned in
the restricted ensemble, which we define next.\vspace*{9pt}

\textit{Dynamics conditioned to stay in $\cR_n(Q)$.}
We denote by
$(\st^{n\pm,\xi}_{Q,t},t\geq0)$
the process
$(\s^{n\pm,\xi}_{Q,t},t\geq0)$
conditioned to stay
in $\cR_n(Q)$.
Its rates
$\widetilde c^{n\pm}_{Q}(x,\s)$
are identical to those of the process\vadjust{\goodbreak}
$(\s^{n\pm,\xi}_{Q,t},t\geq0)$
whenever $\s^x$ belongs to
$\cR_n(Q)$ and they are equal to $0$ whenever
$\s^x\notin\cR_n(Q)$.
As usual, we couple the processes
\[
\bigl(\st^{n\pm,\xi}_{Q,t},t\geq0\bigr),\qquad \bigl(
\s^{n\pm,\xi}_{Q,t},t\geq0\bigr)
\]
so that
\[
\forall\xi\in\cR_n(Q), \forall t< \tau\bigl(\cR_n(Q)
\bigr)\qquad \st^{n\pm,\xi}_{Q,t}= \s^{n\pm,\xi}_{Q,t}.
\]
Finally the measure
$\mt_{Q}^{n\pm}$ defined by
\[
\forall\s\in\cR_n(Q)\qquad \mt_{Q}^{n\pm}(\s) =
\frac{
\mu_{Q}^{n\pm}(\s)} {
\mu_{Q}^{n\pm}(\cR_n(Q))}
\]
is a stationary measure for the process
$(\st^{n\pm,\xi}_{Q,t},t\geq0)$.\vspace*{9pt}

\textit{Local nucleation}.
We say that local nucleation occurs
before $\tb$
%before time $\tb$
%$\exp(\beta\k)$
in the parallelepiped $Q$ starting from
$\xi$ if the process
$(\s^{n\pm,\xi}_{Q,t}, t\geq0)$
%creates a configuration of volume $m_n$ before $\tb$,
exits
$\cR_n(Q)$
before $\tb$.
In words, local nucleation occurs if the process
creates a configuration of energy larger than $\Ga_n$ or of volume
larger than $m_n$
before $\tb$,
that is,
\[
\max \bigl\{ H_Q^{n\pm} \bigl(\s^{n\pm,\xi}_{Q,t}
\bigr)\dvtx t\leq\tb \bigr\} > \Ga_n \quad\mbox{or}\quad \max \bigl\{ \bigl|
\s^{n\pm,\xi}_{Q,t} \bigr|\dvtx t\leq\tb \bigr\} > m_n.
\]
%
%we have
%$$\tau(\cR_n(Q))\leq\tb,$$
%where we recall that the time $\tau(\cR_n(Q))$
%of exit from $\cR_n(Q)$ is defined as
%$$\tau(\cR_n(Q)) = \inf\{ t\geq0\dvtx
%%\exp(\beta\k).$$
%We say that the dynamics creates a $D$--large STC before time $T$
%We consider the process in the box $Q$,
%under minus boundary conditions.
%For sufficiently large $\b$, the shape and energy of the
%saddle is the same as in the case of periodic boundary conditions.
%Although, strictly speaking, $\log\b$ is not a constant,
%The hypothesis on $Q$ ensures that the state space
%grows very slowly with $\b$ (sub-exponentially),
%so that most of the finite-volume results in [CC], [OS],
%and [MO] can be extended.

% \textbf{recurrence:}
% Quote Theorem 3.1 in [MNOS]:
% the probability that the process does not visit the set of states
% with stability level larger than $\g$ for a time $e^{\b(\g+\e)}$ is $

%/or analogous results under other approaches \end{theorem}

%le6.1 #&#
\begin{lemma}\label{fugaup}
Let $n\geq0$, and let $Q$ be a parallelepiped.
We consider the process
$(\s^{n\pm,\mt}_{Q,t},t\geq0)$
in the box $Q$ with $n\pm$ boundary conditions
and initial law the measure
$\mt_{Q}^{n\pm}$.
For any deterministic time $\tb$,
we have for $\b\geq1$,
\begin{eqnarray*}
&&
P \pmatrix{\mbox{local nucleation occurs before $\tau_\beta$}
\cr
\mbox{in the process } \bigl(\s^{n\pm,\mt}_{Q,t}, t\geq0\bigr)}
%(\si^{n\pm,\xi}_{Q,t}, t\geq0)\\
%volume $m_n$}
\\
%4\b
%|Q|^{m_n+1}\tau_\beta
% \exp(-\b\G_{n})
% +
&&\qquad\leq 4
\b(m_n+2)^2 |Q|^{2m_n+2}\tau_\beta \exp(-
\b\G_{n}) + \exp\bigl(-\b|Q|\tb\ln\b\bigr).
\end{eqnarray*}
\end{lemma}
%
%$$(\st^{n\pm,\xi}_{Q,t},t\geq0),
%$$
%
\begin{pf}
To alleviate the text, we drop $\mt$ from the notation, writing
$\s^{n\pm}_{Q,t}$
instead of
$\s^{n\pm,\mt}_{Q,t}$.
To the continuous-time Markov process
$(\st^{n\pm}_{Q,t},t\geq0)$,
we associate in a standard way a discrete-time Markov chain
\[
\bigl(\st^{n\pm}_{Q,k},k\in\N\bigr).
\]
We define first
the time of jumps. We set $\tau_0=0$ and for
$k\geq1$,
\[
\tau_{k} = \inf \bigl\{ t>\tau_{k-1}\dvtx  \st^{n\pm}_{Q,t}
\neq \st^{n\pm}_{Q,\tau_{k-1}} \bigr\}.
\]
We define then
\[
\forall k\in\N\qquad \st^{n\pm}_{Q,k} = \st^{n\pm}_{Q,\tau_{k}}.
\]
Let $X$ be the total number of
arrival times less than $\tb$
of all the
Poisson processes associated to\vadjust{\goodbreak}
the sites of the box $Q$.
The law of $X$ is Poisson with parameter
$\lambda=|Q|\tau_\beta$.
Next, for any $N\geq\lambda$,
\begin{eqnarray*}
P(X\geq N) &=& \sum_{i\geq N}\frac{\lambda^i}{i!}\exp(-
\lambda)
\leq \lambda^N \exp(-\lambda) \sum
_{i\geq N}\frac{N^{i-N}}{i!}
\\
&=& \biggl(\frac{\lambda}{N} \biggr)^N \exp(-\lambda) \sum
_{i\geq N}\frac{N^{i}}{i!} \leq \biggl(\frac{\lambda}{N}
\biggr)^N \exp(N-\lambda).
\end{eqnarray*}
Thus
\[
P(X\geq4\b\lambda) \leq \exp(-\b\lambda\ln\b).
\]
%
%For $A$ a subset of $\Z^d$, we define its outer vertex boundary
%$\dout A$ as
%$$\dout A =
%Let us define
%Any configuration in
%$\din\cR_n(Q)$ has volume $m_n$, therefore
%(\din\cR_n(Q)) \leq
% \leq
%|Q|^{m_n}
The measure
$\mt_{Q}^{n\pm}$
is a stationary measure for the Markov chain
$(\st^{n\pm}_{Q,k})_{k\geq0}$, thus
\begin{eqnarray*}
P \bigl(\tau\bigl(\cR_n(Q)\bigr)\leq\tb \bigr) &\leq& P \bigl(\exists
t\leq\tb, \s^{n\pm}_{Q,t}\notin \cR_n(Q) \bigr)
\\
&\leq& P \bigl(X\leq4\b\lambda, \exists t\leq\tb, \s^{n\pm}_{Q,t}
\notin \cR_n(Q) \bigr) + P (X> 4\b\lambda).
\end{eqnarray*}
The second term is already controlled. Let us estimate the first term,
\begin{eqnarray*}
&&
P \bigl(X\leq4\b\lambda, \exists t\leq\tb, \s^{n\pm}_{Q,t}\notin \cR_n(Q) \bigr)
\\
&&\qquad
\leq P \bigl(X\leq4\b\lambda, \exists k\leq X, \s^{n\pm}_{Q,0},
\ldots, \s^{n\pm}_{Q,k-1}\in %\din
\cR_n(Q),
\s^{n\pm}_{Q,k}\notin \cR_n(Q) \bigr)
\\
&&\qquad
\leq \sum_{1\leq k\leq4\beta\lambda} %\mathop{\sum_{\eta\in\cR_n(Q)}}_{\rho\not\in\cR_n(Q)}
\sum
_{\eta\in\cR_n(Q)} \sum_{\rho\in\partial\cR_n(Q)} P \bigl(
\s^{n\pm}_{Q,k-1}= \st^{n\pm}_{Q,k-1}=
\eta, \s^{n\pm}_{Q,k}=\rho \bigr).
\end{eqnarray*}
Next, for any
$\eta\in\cR_n(Q), \rho\in\partial\cR_n(Q)$,
\begin{eqnarray*}
&&
P \bigl( %\exists k\leq X
\s^{n\pm}_{Q,k-1}= \st^{n\pm}_{Q,k-1}=
\eta, \s^{n\pm}_{Q,k}=\rho \bigr) %|Q|\tau_\beta
\\
&&\qquad
\leq \mt_{Q}^{n\pm}( \eta) \exp \bigl( -\b\max \bigl(0,
H^{n\pm}_Q(\rho)-H^{n\pm}_Q(\eta) \bigr)
\bigr)
\\
&&\qquad
\leq \exp \bigl( -\b\max \bigl( H^{n\pm}_Q(
\rho),H^{n\pm}_Q(\eta) \bigr) \bigr)
\\
&&\qquad
\leq \exp \bigl( -\b E \bigl(\cR_n(Q), \{ -1,+1 \}^Q
\setminus \cR_n(Q) \bigr) \bigr)
\leq \exp(-\b
\Ga_n).
\end{eqnarray*}
Coming back in the previous inequalities, we get
\begin{eqnarray*}
&&
P \bigl(X\leq4\b\lambda, \exists t\leq\tb, \s^{n\pm}_{Q,t}\notin \cR_n(Q) \bigr)
\\
&&\qquad
\leq 4\b\lambda \bigl|\cR_n(Q) \bigr| \bigl|\partial\cR_n(Q) \bigr| \exp(-\b
\Ga_n)
\\
&&\qquad
\leq 4\b\lambda (m_n+1) |Q|^{m_n} (m_n+2)
|Q|^{m_n+1} \exp(-\b\Ga_n),
\end{eqnarray*}
since the number of pluses in
a configuration of
%$\cR_n(Q)$ and
$\partial\cR_n(Q)$ is at most $m_n+1$.
Putting together the previous inequalities, we arrive at
\begin{eqnarray*}
&&P \bigl(\tau\bigl(\cR_n(Q)\bigr)\leq\tb \bigr)
\\
&&\qquad\leq 4\b(m_n+2)^2 |Q|^{2m_n+2}\tau_\beta
\exp(-\b\G_{n}) + \exp\bigl(-\b|Q|\tb\ln\b\bigr)
\end{eqnarray*}
as required.
\end{pf}\eject
%
%For the lower bound, we are not able to compute the inner resistance
%of the critical cycle
%(and follow the standard strategy outlined in Lemma 4.2 in [MO2])
%to prove recurrence in $\minus$.
%Instead, we use the attractivity of the dynamics to prove the
%uniform bound, valid for
%sufficiently large $\b$:
%$$\forall\eta\in\cC_d
%e^{-\b(\G_d-\G_{d-1}+4\e)}.$$
%Indeed, by attractivity,
%$$\PPm{\tp<t} = \inf_{\h\in\cC_d}\PPP{\h}{\tp<t}.$$
%On the other hand, by Lemma \ref{profondita},
%the depth of the deepest cycle in the reference cycle-path is
%at most $\G_{d-1}$ so that by Theorem ref{3.7},
%the probability to follow this cycle-path within
%$e^{\b(\G_{d-1}+\e)}$
%is higher than
%$$2^{-\log^d \b} e^{-\b(\G_d - \G_{d-1} + 3\e)}.$$
%
%We can partition $[0,e^{\b k}]$ into intervals with length
%$e^{\b(\G_{d-1}+\e)}$ getting the bound
%1-\tonda{1 - e^{-\b(\G_d-\G_{d-1}+4\e)}}^{e^{\b(k - \G_{d-1}-\e)}}
%and eqref{fuga1} immediately follows.

%s6.3 #&#
\subsection{\texorpdfstring{Local nucleation or creation of a large $\STC$}
{Local nucleation or creation of a large STC}}
The condition on the initial law and the initial $\STC$ implies that
the process
is initially in a metastable state. We will need to control the $\STC$
created until the
arrival of the supercritical droplets.
Let $Q$ be a parallelepiped included in $\Sigma$.
To build the $\STC$ of the process
$(\s^{n\pm,\xi}_{Q,t}, t\geq0)$
we take into account the $\STC$ initially
present in $\xi$, and we denote by
$\STC_\xi(0,t)$ the resulting $\STC$ on the time interval [0,t$]$.
Hence an element of
$\STC_\xi(0,t)$ is either a $\STC$ of $\STC(0,t)$ which is born
after time $0$, or it is
the union of the $\STC$ of
$\STC(0,t)$ which intersect an initial $\STC$ of
$\STC(\xi)$.
We define then
$\diam\STC_\xi(0,t)$
as
in Section~\ref{sectriangle} by
\[
\diam \STC_\xi(0,t) = \max \biggl( \mathop{\sum
_{\cC\in\STC_\xi(0,t)}}_{
\cC\cap(Q\times\{ 0,t \})\neq\varnothing} \diam \cC,
\mathop{\max
_{\cC\in\STC_\xi(0,t)}}_{
\cC\cap(Q\times\{ 0,t \})=\varnothing} \diam \cC \biggr).
\]
To control this quantity, we will rely on the following inequality:
\[
\diam \STC_\xi(0,t) \leq \mathop{\sum_{\cC\in\STC(\xi)}}_
{\cC\cap Q\neq\varnothing}
\diam\cC + \diam \STC(0,t).
\]
The first term will be controlled with the help of the hypothesis on
the initial $\STC$,
the second term with the help of Theorem \ref{totcontrole}.\vspace*{9pt}

%{\bf Local nucleation}.
%We say that local nucleation occurs
%before time $\tb$
%%$\exp(\beta\k)$
%in the parallelepiped $Q$ starting from
%$\xi$ if for the process
%$(\s^{n\pm,\xi}_{Q,t}, t\geq0)$
%we have
%$$\tau(\cR_n(Q))\leq\tb,$$
\textit{Local nucleation}.
%We say that local nucleation occurs
%before time $\tb$
%%$\exp(\beta\k)$
%in the parallelepiped $Q$ starting from
%$\xi$ if the process
%$(\s^{n\pm,\xi}_{Q,t}, t\geq0)$
%creates a configuration of energy larger than or equal to $\Ga_n$
%before $\tb$,
%i.e.,
%$$\max\Big\{
%H_Q^{n\pm}\big(\s^{n\pm,\xi}_{Q,t}\big)\dvtx t\leq\tb
% \Big\} \geq \Ga_n.$$
%creates a configuration of volume $m_n$ before $\tb$,
%i.e.,
%$$\max\big\{
% \big\} \geq m_n.$$
%we have
%$$\tau(\cR_n(Q))\leq\tb,$$
%where we recall that the time $\tau(\cR_n(Q))$
We say that local nucleation occurs
before $\tb$
%before time $\tb$
%$\exp(\beta\k)$
in the parallelepiped $Q$ starting from
$\xi$ if the process
$(\s^{n\pm,\xi}_{Q,t}, t\geq0)$
%creates a configuration of volume $m_n$ before $\tb$,
exits
$\cR_n(Q)$
before $\tb$.
In words, local nucleation occurs if the process
creates a configuration of energy larger than $\Ga_n$ or of volume
larger than $m_n$
before $\tb$,
that is,
\[
\max \bigl\{ H_Q^{n\pm} \bigl(\s^{n\pm,\xi}_{Q,t}
\bigr)\dvtx t\leq\tb \bigr\} > \Ga_n \quad\mbox{or}\quad \max \bigl\{ \bigl|
\s^{n\pm,\xi}_{Q,t} \bigr|\dvtx t\leq\tb \bigr\} > m_n.
\]

\vspace*{9pt}

\textit{Creation of large $\STC$}.
We say that the dynamics creates a large $\STC$
before time $\tb$
in the parallelepiped $Q$ starting from
$\xi$ if for the process
$(\s^{n\pm,\xi}_{Q,t}, t\geq0)$, we have
%there exists $\cC$ in $STC(0,\tb)$ such that
%
\[
\diam\STC(0,\tb) \geq \ln\ln\beta.
\]
We denote by
$\cR(Q)$ the event
\[
\cR(Q) = \pmatrix{ \mbox{neither local nucleation nor creation}
\cr
\mbox{of a
large $\STC$ occurs before time $\tb$}
\cr
%$\exp(\beta\k)$
\mbox{in the
parallelepiped $Q$ starting from $\xi$}}.
\]

The next proposition gives
a control on the number of these events in a box of
subcritical volume until time $\tb$.
%
%pr6.2 #&#
\begin{proposition}
\label{nuclei}
%Let $(\s^{n\pm,\xi}_{\Sigma,t}, t\geq0)$
%be the process in $\Sigma$ with $n\pm$ boundary condition and
%initial configuration $\xi$.
%$\xi$ satisfies the hypothesis at rank $n$.
%the law $\mu$ of the initial configuration $\xi$
%is the Dirac mass
%on the configuration
%equal to $-1$ everywhere on $\Sigma$.
Let $n\in\{ 1,\ldots, d \}$.
We suppose that
the
hypothesis on the initial law at rank $n$ is satisfied.
%and the hypothesis on the initial $\STC$ are satisfied.
Let $R_\beta$ be a parallelepiped whose volume satisfies
\[
\limsup_{\b\to\infty} \frac{1}{\beta} \ln |R_\beta| \leq
nL_n.
\]
%
%Let $\phi(\beta)$ be a function going to $\infty$ with $\beta$.
The probability that
for the process
$(\s^{n\pm,\xi}_{\Sigma,t},
t\geq0)$
%0\leq t\leq\tb)\]
%$\phi(\beta)$
$\ln\ln\beta$
local nucleations or creations of a large $\STC$
occur
before time $\tau_\beta$
in $n$-small
parallelepipeds
%translates of $Q_n$
included in $R_\beta$ which are pairwise at distance larger than
$5(d-n+1)\ln\ln\beta$ is super-exponentially small in $\beta$.
\end{proposition}
\begin{pf}
Let us rephrase more precisely the event described in the statement
of the proposition:
there exists a family
$(Q_i,i\in I)$ of
$\ln\ln\beta$ $n$-small
parallelepipeds included in $R_\beta$ such that
\[
\forall i,j\in I\qquad i\neq j \quad\Rightarrow\quad d(Q_i,Q_j)
> 5(d-n+1)\ln\ln\beta,
\]
%
%nucleation events
%occurring in subboxes of $S$
%which are pairwise at distance larger
%than $\ln\beta$, and such that nucleation occurs in each box
and
for $i\in I$,
%nucleation occurs before time
the event
$\cR(Q_i)$ does not occur for the process
$(\s^{n\pm,\xi}_{Q_i,t}, t\geq0)$.
Denoting this event by $\cE$, we have
\[
P(\cE) \leq \sum_{(Q_i)_{i\in I}} P \biggl( %\begin{matrix}
%(\s^{n\pm,\xi}_{Q_i,t}, t\geq0)\\
\bigcap_{i\in I} \cR(Q_i )^c
\biggr)
\]
where the sum runs over all the possible choices of boxes
$(Q_i)_{i\in I}$.
We condition next on the initial configurations $(\s_i,i\in I)$ in the boxes
$(Q_i,{i\in I})$,
\begin{eqnarray*}
P(\cE) &\leq& \sum_{(Q_i)_{i\in I}} \sum
_{(\s_i)_{i\in I}} P \biggl( \bigcap_{i\in I}
\cR(Q_i)^c %\begin{matrix}
%(\s^{n\pm,\xi_i}_{Q_i,t}, t\geq0)\\
\big| \forall i\in I,
\xi|_{Q_i}=\s_i \biggr)
\\
&&\hspace*{52pt}{}\times \mu ( \forall i\in I, \xi|_{Q_i}=\s_i ).
\end{eqnarray*}
Once
the initial configurations $(\s_i,i\in I)$
are fixed, the nucleation events in the boxes
$(Q_i,i\in I)$
become independent because
they depend on Poisson processes associated to disjoint boxes.
Thanks to
the geometric condition imposed on the boxes, we can apply the estimates
given by the hypothesis on the initial law $\mu$,
\begin{eqnarray*}
P(\cE) &\leq& \sum_{(Q_i)_{i\in I}} \sum
_{(\s_i)_{i\in I}} \prod_{i\in I} P \bigl(
\cR(Q_i )^c \big| \xi|_{Q_i}=\s_i \bigr)
%%\mbox{for }i\in I,
%(\s^{n\pm,\xi_i}_{Q_i,t}, t\geq0)\\
\phi_n(\beta)
\rho_{Q_i}^{n\pm}(\s_i)
\\
&=& \sum_{(Q_i)_{i\in I}} \prod_{i\in I}
\biggl( \phi_n(\beta) \sum_{\s_i} P
\bigl(\cR(Q_i)^c \big| \xi|_{Q_i}=\s_i
\bigr) \rho_{Q_i}^{n\pm}(\s_i) \biggr).
\end{eqnarray*}
Let us fix $i\in I$, and let us estimate the term inside the big parenthesis.
Let $Q$ be an $n$-small box.
We write
\begin{eqnarray*}
&&
\sum_{\eta} P \bigl(\cR(Q)^c |
\xi|_{Q}=\eta \bigr) \rho_{Q}^{n\pm}(\eta)
\\
&&\qquad
\leq \sum_{\eta} P\pmatrix{ %\mbox{for }i\in I,
\mbox{the
process } \bigl(\s^{n\pm,\eta}_{Q,t}, t\geq0\bigr)
\vspace*{1pt}\cr
\mbox{nucleates before time $\tb$} } \rho_{Q}^{n\pm}(\eta)
\\
&&\qquad\quad{}+ \sum_{\eta} P \pmatrix{ %\mbox{for }i\in I,
\mbox{the
process } \bigl(\s^{n\pm,\eta}_{Q,t}, t\geq0\bigr) \mbox{ creates}
\vspace*{1pt}\cr
\mbox{a large $\STC$ before nucleating} } \rho_{Q}^{n\pm}(
\eta).
\end{eqnarray*}
First,
by Theorem \ref{totcontrole}, the
probability that the process
$(\s^{n\pm,\eta}_{Q,t}, t\geq0)$
creates a large $\STC$ before nucleating is $\SES$.
Second,
\[
\sum_{\eta\notin\cR_n(Q)} P \pmatrix{ %\mbox{for }i\in I,
\mbox{the
process } \bigl(\s^{n\pm,\eta}_{Q,t}, t\geq0\bigr)
\vspace*{1pt}\cr
\mbox{nucleates before time $\tb$} } \rho_{Q}^{n\pm}(\eta)
\leq 2^{|Q|} \exp(-\beta\Ga_n).
\]
Third, for $ \eta\in\cR_n(Q)$,
using the notation of Section~\ref{finitebox},
\[
\rho_{Q}^{n\pm}(\eta) \leq \bigl|\cR_n(Q)\bigr|
\mut_{Q}^{n\pm}(\eta) \leq (m_n+1)|Q|^{m_n}
\mut_{Q}^{n\pm}(\eta),
\]
whence, using Lemma \ref{fugaup},
\begin{eqnarray*}
&&
\sum_{\eta\in\cR_n(Q)} P \pmatrix{ %\mbox{for }i\in I,
\mbox{the
process } \bigl(\s^{n\pm,\eta}_{Q,t}, t\geq0\bigr)
\vspace*{1pt}\cr
\mbox{nucleates before time $\tb$}} \rho_{Q}^{n\pm}(\eta)
\\
&&\qquad\leq %2^{|Q|}
(m_n+1)|Q|^{m_n} %\sum_{\xi\in\cR_n(Q)}
P \pmatrix{
\mbox{the process } \bigl(\s^{n\pm,\mt}_{Q,t}, t\geq0
\bigr)
\vspace*{1pt}\cr
\mbox{nucleates before time $\tb$}} %\left(
%(\st^{n\pm,
%}_{Q,t}, t\geq0)
%volume $m_n$}
%P
%%\mbox{for }i\in I,
%(\st^{n\pm}_{Q,t}, t\geq0)\mbox{ reaches the}\\
%%\rho_{Q}^{n\pm}(\xi)
\\
&&\qquad\leq 4\b(m_n+2)^3 (n\ln\beta)^{d(2m_n+3)}
\tau_\beta \exp(-\b\G_{n}) + \SES. % \exp(-\b|Q|\tb\ln\b).
%{4\beta m_n(n\ln\beta)^{d(2m_n+1)}}
\end{eqnarray*}
Substituting
these estimates into the last inequality on $P(\cE)$, we obtain
\begin{eqnarray*}
P(\cE) &\leq& %\SES +
\bigl( |R_\beta | (n\ln
\beta)^d \bigl(2^{(n\ln\beta)^d}
+4\beta{(m_n+2)^3(n\ln\beta)^{2dm_n+3d}}
\tau_\beta \bigr) \\
&&\hspace*{134pt}{}\times \phi_n(\beta)\exp(-\beta
\Ga_n) +\SES \bigr)^{ |I|}.
\end{eqnarray*}
Since $|I|=
\ln\ln\beta$ and
%$\kappa<\kappa_n=\Ga_n-nL_n$,
%
\[
\limsup_{\b\to\infty} \frac{1}{\beta} \ln \bigl(
|R_\beta| \tau_\beta \phi_n(\beta) \exp(-\beta
\Ga_n) \bigr) < nL_n+\k_n-\G_n = 0,
\]
we conclude that the above quantity is $\SES$.
\end{pf}

%s6.4 #&#
\subsection{Control of the metastable space--time clusters}
\label{fvs}
The key result
is the following control on the size of the space--time clusters
in the configuration.
The next proposition states the result at rank 0, the theorem thereafter
states the result at rank $n\geq1$.
%
%pr6.3 #&#
\begin{proposition}\label{PT2}
%Let $d\geq1$ be an integer and
%let $L$ be positive real number.
%$n\geq0$ be two integers.
%$\L_\beta=\L(\exp(\beta L))$ be a cubic
%box of sidelength $\exp(\beta L)$.
%where
%$\L^{d}$ is a $d$ dimensional cubic box and
%$\L^{n}$ is an $n$ dimensional cubic box.
%metastable
We suppose that the law $\mu$ of the initial configuration
$\xi$ satisfies the hypothesis at rank $0$.
Let
$\tau_\beta$ be a time satisfying
\[
\limsup_{\b\to\infty} \frac{1}{\beta} \ln \tau_\beta <
\k_0=0.
\]
The probability\vspace*{1pt} that
a $\STC$
%a space--time cluster
of diameter
larger than
$\ln\ln\beta$
is created
%before time
in the process
$(\s^{0\pm,\xi}_{\Sigma,t}, 0\leq t\leq
\tau_\beta
)$
is $\SES$.
%super--exponentially small in $\beta$.
\end{proposition}
\begin{pf}
With $n=0$, we have
\[
\Sigma = \L^{d}(\ln\beta),\qquad \Ga_0=\k_0=L_0=m_0=0,
\]
the boundary condition is plus on
$\partial^{\mathrm{out}}\Sigma$ and
$\cR_0(Q)=\{ \minus \}$
for any box $Q$.
By the hypothesis
on $\mu$ at rank $0$,
%amounts to saying that
the initial law $\mu$ is the Dirac mass on the configuration
equal to $-1$ everywhere on $\Sigma$.
Now
\begin{eqnarray*}
&&P \bigl(\exists \cC\in\STC(0, \tau_\beta )\mbox{ with }\diam\cC\geq
\ln\ln\beta \bigr)
\\
&&\qquad\leq P\pmatrix{
\mbox{there are at least $\ln\ln\beta$ arrival times less than $\tb$
}
\cr
\mbox{for the Poisson processes associated to the sites of $\Sigma$}}
\\
&&\qquad=
P(X\geq\ln\ln\beta), %P\left(
%impact of a spatial Poisson process in
%$\Sigma\times[0,\tau_\beta]$
%}\\
%Poisson processes associated to
%the sites of $\Sigma$}
%,=
\end{eqnarray*}
where $X$ is a variable whose law is Poisson with parameter
\[
\lambda = |\Sigma|\tau_\beta = (\ln\beta)^d\tb.
\]
So
\[
P(X\geq\ln\ln\beta) = \sum_{k\geq\ln\ln\beta} \exp(-\lambda)
\frac{\lambda^k}{k!} \leq 3\lambda^{\ln\ln\beta},
\]
which is $\SES$.
\end{pf}

For the case $n=0$ the initial configuration is $-1$ everywhere, and
all the $\STC$
born before the initial configuration are dead. In the case $n\geq1$,
the situation is more
delicate, and we must deal with $\STC$ born in the past.
To build the $\STC$ of the process
$(\s^{n\pm,\xi}_{\Sigma,t}, t\geq0)$
we take into account the $\STC$ initially
present in $\xi$, and we denote by
$\STC_\xi(0,t)$ the resulting $\STC$ on the time interval [0,t$]$.
Hence an element of
$\STC_\xi(0,t)$ is either a $\STC$ of $\STC(0,t)$ which is born
after time $0$, or it is
the union of the $\STC$ of
$\STC(0,t)$ which intersect an initial $\STC$ of
$\STC(\xi)$.
We recall that
$\STC(\xi)$ denotes the initial $\STC$
present in $\xi$, and these $\STC$s are unions of clusters
of pluses of $\xi$.
%
%th6.4 #&#
\begin{theorem}\label{T2}
%Let $d\geq1$ be an integer and
%let $L$ be positive real number.
%$n\geq0$ be two integers.
%$\L_\beta=\L(\exp(\beta L))$ be a cubic
%box of sidelength $\exp(\beta L)$.
%where
%$\L^{d}$ is a $d$ dimensional cubic box and
%$\L^{n}$ is an $n$ dimensional cubic box.
%metastable
Let $n\in\{ 1,\ldots,d \}$.
%We suppose that the law $\mu$ of the initial configuration
%$\xi$ satisfies the hypothesis at rank $n$.
We suppose that
both the
hypothesis on the initial law at rank $n$
and on the initial $\STC$ present in $\xi$ are satisfied.
Let
$\tau_\beta$ be a time satisfying
\[
\limsup_{\b\to\infty} \frac{1}{\beta} \ln \tau_\beta <
\k_n.
\]
The probability that,
for the process
$(\s^{n\pm,\xi}_{\Sigma,t})_{t\geq0}$,
there exists
%a $\STC$ in
a space--time cluster in
%$\STC$ in
$\STC_\xi(0,\tb)$
%a space--time cluster
of diameter
larger than
$\exp({\b L_n})$
%is created
%%before time
%in the process
%$(\s^{n\pm,\xi}_{\Sigma,t}, 0\leq t\leq
%)$
is $\SES$.
%super--exponentially small in $\beta$.
\end{theorem}

Theorem \ref{T2} is proved by induction over $n$.
We suppose that the result at rank $n-1$ has been proved
and that a $\STC$ of diameter
larger than $\exp(\beta L_n)$ is formed before time $\tb$.
The induction step is long, and it is decomposed in the eleven following
steps:\vspace*{9pt}
%Strategy of the proof of
%Theorem \ref{T2}.}
%We explain here in a loose way the strategy of the proof of
%Theorem \ref{T2}, without entering into the
%technical details. We consider
%a parallelepiped $\Sigma$
%of the form
%$$\Sigma = \L^{n}(\exp(\beta L)) \times\L^{d-n}(\ln\beta)  $$
%with $n\pm$ boundary conditions and the time $\tb=\exp(\beta\kappa)$,
%where $\kappa<\kappa_d$.
%We say that a parallelepiped is $n$-small if
%if all its sides have a length
%larger than $\ln\ln\beta$ and smaller than $n\ln\beta$.

\textit{Step} 1: \textit{Reduction to a box $R_{i,j}$ of side length of
order $\exp(\beta L_n)$}. By a trick going back to the work of
Aizenmann and Lebowitz on bootstrap percolation \cite{AL}, there exists
a $\STC$ of diameter between $\exp(\beta L_n)/2$ and $\exp(\beta
L_n)+1$ which is formed before time $\tb$. In particular there exists a
box $R_{i,j}$ of side length of order $\exp(\beta L_n)$ which is
crossed by a $\STC$ before time $\tb$.\vspace*{9pt}

\textit{Step} 2: \textit{Reduction to a box $S_{i}$
of side length of order $\exp(\beta L_n)/\ln\beta$
devoid of bad events}.
Thanks to Proposition \ref{nuclei}, the number of bad events,
like local nucleation or creation of a large $\STC$, is at most
$\ln\ln\beta$, up to a $\SES$ event.
By a simple counting argument, there exists a box $S_i$ of side length
of order
$\exp(\beta L_n)/\ln\beta$ in the $n$th direction which is
crossed vertically before $\tb$ and in which no bad events occur.
We consider next the dynamics in this box $S_i$
with either $n\pm$ or $n-1\pm$ boundary conditions.\vspace*{9pt}

\textit{Step} 3: \textit{Control of the diameters of the
$\STC$ born in $S_i$ with $n\pm$ boundary conditions}.
By construction, for the dynamics in the box $S_i$
with $n\pm$ boundary conditions, no bad events occur before
time $\tb$, and therefore the process stays in the metastable state.
Until time $\tb$, only small droplets are created, and they survive for
a short time.
We quantify this
in Lemma \ref{essx}, where we prove that
any $\STC$
in $\STC_\xi
(\s^{n\pm,\xi}_{S_i,t}, 0
\leq t\leq
\tau_\beta
)$
has a diameter at most
$(d-n+2)\ln\ln\beta$.\vspace*{9pt}

\textit{Step} 4: \textit{Reduction to a flat box $\Delta_{i,j}\subset S_i$
of height $\ln\beta$
crossed vertically in a
time
$\exp (\beta(\kappa-L_n) )/
(\ln\beta)^2$}.
The box $S_{i}$
has height of order $\exp(\beta L_n)/\ln\beta$.
In the dynamics restricted to $S_i$ with
$n-1\pm$ boundary conditions, the box $S_i$
is vertically crossed by a $\STC$ in a time $\tb$.
From the result of step 3,
we conclude that the crossing $\STC$ emanates either from
the bottom or the top of $S_i$ because the vertical crossing
can occur only with the help of the boundary conditions. This $\STC$
has to be born close to the top or the bottom of $S_i$,
and it propagates then toward the middle plane of $S_i$.
We partition $S_i$ in slabs of height $\ln\beta$, the number of
these slabs is
of order $\exp(\beta L_n)/(\ln\beta)^2$.
By summing the crossing times of each of these slabs, we obtain that
one slab,
denoted by $\Delta_{i,j}$,
has to be crossed vertically in a time
$\exp (\beta(\kappa-L_n) )/
(\ln\beta)^2$.
We denote by $\cT_a$ the following event:
%$\bullet$ The event $\cR(S_i)$ occurs.
At time $a$, the set
$\Delta_{i,j}$ has not been touched by an $\STC$ emanating from top
or bottom of $S_i$
in the process
$(\s^{n-1\pm,\xi}_{S_i,t}, t\geq0)$.
We denote by $\cV_b$ the following event:
At time $b$,
the set
$\Delta_{i,j}$ is vertically crossed
in the process
$(\s^{n-1\pm,\xi}_{S_{i},t},t\geq0)$.
We show that there exist two integer values $a<b$ such that
$b-a<
\exp (\beta(\kappa-L_n) )/
(\ln\beta)^2$, and the events
$\cT_a$ and
$\cV_b$ both occur.\vspace*{9pt}

\textit{Step} 5: \textit{Conditioning on the configuration at the time of arrival
of the large $\STC$}.
We want to estimate the probability of the event
$\cT_a\cap
\cV_b$. This event will have a low probability because it requires
that the slab
$\Delta_{i,j}$ is vertically crossed
too quickly, before it had time to relax to equilibrium.
To this end, we condition with respect to the configuration in
$\Delta_{i,j}$
at time $a$,
and we estimate the probability of the vertical
crossing in a time $b-a$.
We first replace the condition that
no bad events occur before
time $\tb$ by the weaker condition that
no bad events occur before
time $a$ (otherwise the conditioned dynamics after time $a$ would be
much more complicated).
We then perform the conditioning
with respect to the configuration in
$\Delta_{i,j}$
at time $a$.
We denote by $\zeta$ this configuration, by $\nu$ its law
and by $\STC(\zeta)$ the $\STC$ present in $\zeta$.
The idea is to apply the induction hypothesis to the process
in $\Delta_{i,j}$ between times $a$ and $b$. To this end, we check
that $\nu$ and $\STC(\zeta)$ satisfy the hypothesis at rank
$n-1$.\vspace*{9pt}

\textit{Step} 6: \textit{Check of the hypothesis on the initial $\STC$ at rank $n-1$}.
We use the initial hypothesis on the $\STC$ at rank $n$ and the fact
that no bad events, like nucleation or creation of a large $\STC$,
occur until time $a$ to obtain the appropriate control on the
$\STC$ at time $a$.
The factor $(d-n+1)\ln\ln\beta$ is tuned adequately to perform
the induction step. The condition is stronger at step $n$
than at step $n-1$. Indeed, the hypothesis is done at step $n$ on
the initial $\STC$, and because of the metastable dynamics, the
diameters of the $\STC$ might increase by $\ln\ln\beta$ until
the arrival of the supercritical droplets. Thus the hypothesis on
the $\STC$ at rank $n-1$ is still fulfilled.\vspace*{9pt}

\textit{Step} 7: \textit{Check of the hypothesis on the initial law at rank $n-1$.}
Similarly, we use the hypothesis on the initial law at rank $n$ and the fact
that no bad events, like nucleation or creation of a large $\STC$,
occur until time $a$ to obtain the appropriate decoupling on the law of
the configuration
at time $a$.
The hypothesis on the law at rank $n$ implies that small boxes at distance
larger than $5(d-n+1)\ln\beta$ are independent. Until time $a$,
no bad events occur, and hence the metastable dynamics inside a small box
$Q$
can
only be influenced by events happening\vadjust{\goodbreak} at distance $\ln\ln\beta$
from $Q$,
that is, inside a slightly larger box $R$.
This way we obtain the appropriate decoupling on boxes which are
at distance larger than
$5(d-n+2)\ln\beta$.\vspace*{9pt}

\textit{Step} 8: \textit{Comparison of
$\mt_{R}^{n\pm}|_{Q}$
and $\rho_{Q}^{n-1\pm}$.}
To obtain the appropriate bounding factor
%for the law in a small box,
we have to
prove that
if $Q,R$ are two parallelepipeds which are $n$-small and such that
$Q\subset R$,
then for any configuration $\eta$ in $Q$,
\[
\mt_{R}^{n\pm} ( \sigma|_{Q} =\eta ) \leq
\phi_{n-1}(\beta) \rho_{Q}^{n-1\pm}(\eta),
\]
where $\phi_{n-1}(\beta)$
is a function depending only upon $\beta$.
This is done with the help of three geometric lemmas.
First we show that a configuration $\sigma$
having at most $m_n$ pluses
and such that
$H^{n\pm}_{R}(\sigma)
%-H^{n\pm}_{Q}(\underline{\minus})
\leq
\Ga_{n-1}
$
can have at most
$m_{n-1}$ pluses.
The next point is that, when the number of pluses in the configuration
$\eta$
is less
than $m_{n-1}$, the Hamiltonian in $R$ with $n\pm$ boundary conditions
will always be larger than
the Hamiltonian in $Q$ with $n-1\pm$ boundary conditions,
up to a polynomial correcting factor.\vspace*{9pt}

\textit{Step} 9:
\textit{Reduction to a box $\Phi$
of side length of order $\exp(\beta L_{n-1})$.}
We are now able to apply the induction hypothesis at rank $n-1$:
Up to a $\SES$ event, there is
no space--time cluster of diameter
larger than
$\exp({\b L_{n-1}})$
for the process in $\Delta$ with
$n-1\pm$ boundary conditions.
Therefore the vertical crossing of $\Delta$ has to occur
in a box $\Phi$ of side length
of order
$\exp({\b L_{n-1}})$.\vspace*{9pt}

\textit{Step} 10:
\textit{Reduction to boxes $\Phi_i\subset\Phi$
of vertical side length of order $\ln\beta/\break\ln\ln\beta$.}
We partition $\Phi$ in slabs $\Phi_i$
of height $\ln\beta/\ln\ln\beta$, the number of
these slabs is
of order $\ln\ln\beta$.
We can choose a subfamily of slabs such that two slabs of
the subfamily are at
distance larger than $5(d-n+2)\ln\ln\beta$.
Since $\Phi$ endowed with $n-1\pm$ boundary conditions is vertically
crossed before time
$\exp (\beta(\kappa-L_n) )/
(\ln\beta)^2$, so are each of these slabs $\Phi_i$.\vspace*{9pt}

\textit{Step} 11: \textit{Conclusion of the induction step.}
Each slab $\Phi_i$ is crossed, and each of these crossings
implies that a large $\STC$ is created.
The dynamics in each slab $\Phi_i$ with
$n-1\pm$ boundary conditions are essentially independent,
thanks to the boundary conditions and the hypothesis on the initial
law. It follows that
the probability of creating simultaneously these $\ln\ln\beta$
large $\STC$ is $\SES$.

%f5 #&#
\begin{figure}

\includegraphics{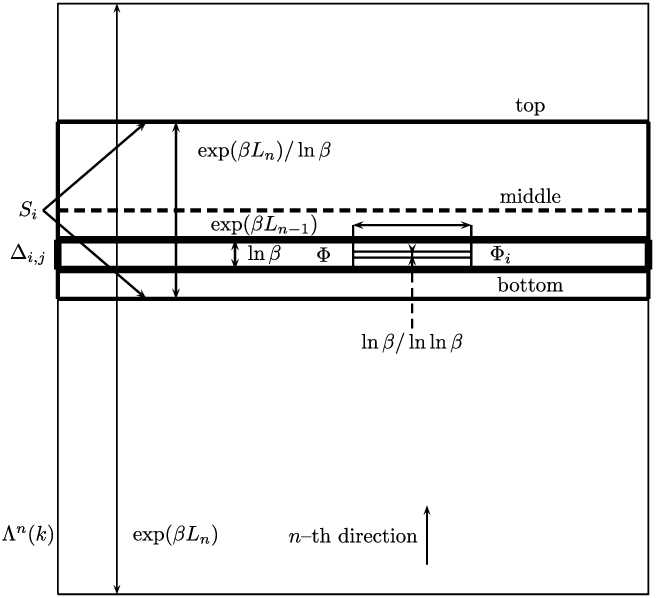}

\caption{Reduction from $\L^n(k)$ to $\Phi_i$.}
\end{figure}

We start now the precise proof, which follows the above strategy.
We suppose that the result at rank $n-1$ has been proved
and that a $\STC$ of diameter
larger than $\exp(\beta L_n)$ is formed before time $\tb$.\vspace*{9pt}

\textit{Step} 1: \textit{Reduction to a box $R_{i,j}$
of side length of order $\exp(\beta L_n)$.}
Let us consider the function
\[
f(t) = \max \bigl\{ \diam\cC\dvtx \cC\in\STC_\xi(0,t) \bigr\}.
\]
This function is nondecreasing; it changes when a spin flip creates
a larger $\STC$ by merging two or more existing $\STC$. Suppose there\vadjust{\goodbreak}
is a
spin flip at time $t$. Just before the spin flip, the largest $\STC$ had
diameter at most
\[
f(t-) = \mathop{\lim_{ s<t}}_{ s\to t}f(s).
\]
Hence after the spin flip, the largest $\STC$ has diameter at most
$2f(t-)+1$. Therefore
\[
\forall t\geq0\qquad f(t) \leq 2f(t-)+1.
\]
With the same reasoning applied to a specific $\STC$, we get
the following result.
%
%le6.5 #&#
\begin{lemma}
\label{reason}
Let $D$ be such that
\[
D \geq \max \bigl\{ \diam\cC\dvtx \cC\in\STC(\xi) \bigr\}.
\]
Let $\cC$ be a $\STC$ in $\STC_\xi(0,t)$ having diameter larger
than $D$.
There exists $s\leq t$ and $\cC'$ a
$\STC$ in $\STC_\xi(0,s)$ such that
\[
\cC'\subset\cC,\qquad D\leq\diam\cC'\leq2D.
\]
%
% $\cC'\subset\cC$ and moreover
%$ D\leq\diam\cC'\leq2D$.
\end{lemma}
The hypothesis on the initial $\STC$ present in $\xi$ implies that
\[
\max \bigl\{ \diam\cC\dvtx \cC\in\STC(\xi) \bigr\} \leq (d-n+1)\ln\ln\beta.
\]
Therefore, if
\[
f( \tau_\beta) \geq \exp(\beta L_n),
\]
then, by Lemma \ref{reason}, there exists a random time $T\leq
\tau_\beta$ and
$\cC\in\STC_\xi(0,T)$ such that
%such that
%2\exp(\beta L_n)+1.
%Let $\cC\in\STC_\xi(0,T)$ such that
%
\[
\exp(\beta L_n) \leq \diam\cC \leq 2\exp(\beta L_n).
\]
Let $\Phi$ be the smallest $n$-dimensional box such that
\[
\cC \subset \bigl(\Phi\times \L^{d-n}(\ln\beta) \bigr)\times[0,T].
\]
With the help of Lemma \ref{cresSTC}, we observe that
the box $\Phi$ is crossed by a $\STC$ before time
$\tau_\beta$, where the meaning of
``crossed'' is explained next.

%de6.6 #&#
\begin{definition}
\label{cro}
An $n$-dimensional box $\Phi$ is said to be crossed by a $\STC$
before time $t$
%if there exists $\cC$ a $\STC$ in $\STC_\xi(0,t)$ for the dynamics
%restricted to
if,
for the dynamics restricted to
$\Phi\times
\L^{d-n}(\ln\beta)$
with initial configuration $\xi$ and
$n\pm$ boundary condition,
there exists $\cC$ in $\STC_\xi(0,t)$
whose projection on the first
$n$ coordinates intersects two opposite faces of~$\Phi$.
%The box is said to be crossed vertically in case the opposite faces
%which are touched
%by the same STC
%are the bottom and the top faces of $\Phi$.
\end{definition}

With this definition, we have
\begin{eqnarray*}
\hspace*{-8pt}&&
P \bigl(\exists \cC\in\STC_\xi(0, \tau_\beta )\mbox{ with
}\diam\cC\geq \exp(\beta L_n) \bigr)
\\
&&\qquad\leq P \pmatrix{ %\mbox{there exists }
\exists \Phi\mbox{ $n$-dimensional
box $\subset \L^{n}(L_\beta)$,}
\cr
\exp(\beta
L_n)\leq \diam\Phi\leq 2\exp(\beta L_n),
\cr
\mbox{$\Phi$
is crossed by a $\STC$ before time } \tau_\beta %\mbox{ and such that}\cr
} %of sidelength $\exp(\beta L)$
\\
&&\qquad\leq \bigl|\L^n( L_\beta)\bigr|\times 2\exp(\beta L_n)
\times \max_{x,k} %{ \exp(\beta L_n)\leq k\leq
%2\exp(\beta L_n)}
P \pmatrix{
\mbox{the box } \bigl(x+\L^n(k) \bigr)\times \L^{d-n}(\ln
\beta)
\cr
\mbox{is } \mbox{crossed by a $\STC$ before time } \tau_\beta
}
\end{eqnarray*}
where the maximum is taken over $x,k$ such that
\[
\exp(\beta L_n)\leq k \leq 2\exp(\beta L_n),\qquad \bigl(x+
\L^n(k) \bigr) \subset \L^n( L_\beta).
\]
Let us now fix $x,k$ as above. For simplicity we take $x=0$,
and let us suppose that
$\L^n(k)\times
\L^{d-n}(\ln\beta)$
is crossed
by a $\STC$
before time
$\tau_\beta$
for the process with
initial configuration $\xi$ and $n\pm$ boundary condition.
We can suppose, for instance, that
$\L^n(k)\times
\L^{d-n}(\ln\beta)$
is crossed vertically, that is, that the crossing occurs along the
$n$th coordinate.
Using the monotonicity with respect to the boundary conditions,
we observe that,
for any $i,j$ such that $-k/2\leq i\leq j\leq k/2$,
the parallelepiped
\[
R_{i,j} = \L^{n-1}(k)\times[i,j]\times \L^{d-n}(\ln
\beta)
\]
is also crossed
vertically before time
$\tau_\beta$
for the process
with
initial configuration
$\xi|_{ R_{i,j}} $ and $n-1\pm$ boundary condition
on $R_{i,j}$.\vspace*{9pt}
%P

\textit{Step} 2: \textit{Reduction to a box $S_{i}$
of side length of order $\exp(\beta L_n)/\ln\beta$
devoid of bad events.}
With $k$ defined above,
we consider next the collection of the sets
\[
S_i = \L^{n-1}(k)\times \biggl[ \frac{(2i)k}{4\ln\beta},
\frac{(2i+1)k}{4\ln\beta} \biggr]\times \L^{d-n}(\ln\beta),\qquad |i|<\ln\beta-
\frac{1}{2}.
\]
%
%From the previous discussion, we see that all the sets $S_i$ are
%crossed
%vertically before time
%$\exp(\beta\k)$ for the process with
%initial configuration $\xi$ and $(n-1)\pm$ boundary condition.
These sets are pairwise at distance larger than $\ln\beta$.
By Proposition \ref{nuclei}, up to a $\SES$ event, there exists a set
$S_i$ in which
the event
\[
\cR(S_i) %\exp(\beta\k),D)
= \mathop{\bigcap
_{ Q \ n\mbox{-}\mathrm{small}}}_
{ Q\subset S_i} \cR(Q )
\]
occurs.
This means that
neither local nucleation nor the creation of a large $\STC$
occurs before time
$\tau_\beta$
for the process in $S_i$ with
initial configuration
$\xi|_{ S_{i}} $ and $n\pm$ boundary condition.
From now onward, we will study what is happening in this particular set $S_i$.
Let us define
\begin{eqnarray*}
\mbox{bottom} &=& \L^{n-1}(k)\times \biggl\{ \frac{(2i)k}{4\ln\beta} \biggr\}
\times \L^{d-n}(\ln\beta),
\\
\mbox{top} &=& \L^{n-1}(k)\times \biggl\{ \frac{(2i+1)k}{4\ln\beta} \biggr\}
\times \L^{d-n}(\ln\beta).
\end{eqnarray*}
By Lemma \ref{cresSTC},
any $\STC$ of the process
\[
\bigl(\s^{n-1\pm,\xi}_{S_i,t}, 0 \leq t\leq \tau_\beta
\bigr),
\]
which intersects neither top nor bottom is also a $\STC$
of the process
\[
\bigl(\s^{n\pm,\xi}_{S_i,t}, 0 \leq t\leq \tau_\beta
\bigr),
\]
because it has not been ``helped'' by the
$n-1\pm$ boundary condition.\vspace*{9pt}

\textit{Step} 3: \textit{Control of the diameters of the
$\STC$ born in $S_i$ with $n\pm$ boundary conditions.}

%le6.7 #&#
\begin{lemma}\label{essx}
On the event
$\cR(S_i)$, any $\STC$
in $\STC_\xi
(\s^{n\pm,\xi}_{S_i,t}, 0
\leq t\leq
\tau_\beta
)$
has a diameter at most
$(d-n+2)\ln\ln\beta$.
\end{lemma}
%
%On the event
%$\cR(S_i)$, any such $\STC$ has a diameter at most
%$(d-n+2)\ln\ln\beta$.
%
\begin{pf}
Indeed, suppose that there exists $\cC$ in
$\STC_\xi
(\s^{n\pm,\xi}_{S_i,t}, 0
\leq t\leq
\tau_\beta
)$
with
\[
\diam\cC > (d-n+2)\ln\ln\beta.
\]
By Lemma \ref{reason}, there exists $T\leq\tb$ and
%$\cC'\subset\cC$ a $\STC$ in
$\cC'$ in
$\STC_\xi
(\s^{n\pm,\xi}_{S_i,t}, 0
\leq t\leq
T
)$
such that
\[
(d-n+2)\ln\ln\beta \leq \diam\cC' \leq \tfrac{1}{3}\ln\beta.
\]
Let $Q'$ be a box of side length $\ln\beta$ included in $S_i$
and centered
on a point of $\cC'$.
By Lemma \ref{cresSTC},
$\cC'$ is also a $\STC$ of the process
$(\s^{n\pm,\xi}_{Q',t}, 0
\leq t\leq
\tau_\beta
)$.
Yet
\begin{eqnarray*}
\diam\cC' &\leq& \diam \STC_\xi \bigl(
\s^{n\pm,\xi}_{Q',t}, 0 \leq t\leq T \bigr)
\\
&\leq& \mathop{\sum_{\cC\in\STC(\xi)}}_
{\cC\cap Q'\neq\varnothing}\diam\cC +
\diam \STC\bigl( \s^{n\pm,\xi}_{Q',t}, 0\leq t\leq T\bigr)
\\
&\leq& (d-n+1)\ln\ln\beta+ \diam \STC\bigl( \s^{n\pm,\xi}_{Q',t}, 0
\leq t\leq T\bigr).
\end{eqnarray*}
We have used the hypothesis on the initial clusters present in $\xi$
to bound
the sum.
This inequality implies that
\[
\diam \STC\bigl( \s^{n\pm,\xi}_{Q',t}, 0\leq t\leq T\bigr) \geq \ln
\ln\beta.
\]
Hence the events $\cR(Q')$ and $\cR(S_i)$ would not occur.
\end{pf}

\textit{Step} 4: \textit{Reduction to a flat box $\Delta_{i,j}\subset S_i$
crossed vertically in a time $(\ln\beta)^2$.}
By Lemma \ref{cresSTC},
%Therefore
any $\STC$ in
$(\s^{n-1\pm,\xi}_{S_i,t}, 0
\leq t\leq
\tau_\beta
)$
of diameter strictly larger than
$(d-n+2)\ln\ln\beta$
intersects top or bottom.
Since $S_i$ is vertically crossed by time $\tb$, the middle set,
defined by
\[
\mbox{middle} = \L^{n-1}(k)\times \biggl\{ \frac{(2i+{1/2})k}{4\ln\beta} \biggr\}
\times \L^{d-n}(\ln\beta)
\]
is hit before time $\tb$ by a $\STC$ emanating either from the bottom or
from the top of~$S_i$.
Let us define
\begin{eqnarray*}
\tbot(h) &=& \inf \bigl\{ u\geq0\dvtx  \exists \cC\in\STC_\xi\bigl(
\s^{n-1\pm,\xi}_{S_i,t}, 0 \leq t\leq u\bigr),
\\
&&\hspace*{17.2pt}\cC\cap\mbox{bottom}\neq\varnothing, \exists x=(x_1,
\ldots,x_d)\in\cC, x_n=h \bigr\}.
\end{eqnarray*}
Suppose, for instance, that the first $\STC$ hitting middle emanates
from the bottom.
We have then
\[
\tbot \biggl( \frac{(2i+{1/2})k}{4\ln\beta} \biggr) \leq \tb.
\]
Moreover, setting $h=2i/(4\ln\beta)$, we have
\begin{eqnarray*}
&&\tbot \biggl( h+\frac{k}{8\ln\beta} \biggr)
\\
&&\qquad\geq\sum_{1\leq j\leq J} \bigl( \tbot ( h+j\ln\beta )- \tbot
\bigl( h+(j-1)\ln\beta \bigr) \bigr),
\end{eqnarray*}
where
\[
J = \frac{k}{8(\ln\beta)^2}.
\]
Therefore there exists an index $j\leq J$ such that
\[
\tbot ( h+j\ln\beta )- \tbot \bigl( h+(j-1)\ln\beta \bigr) \leq
\frac{\tb}{J}.
\]
Let
%$T_{i,j}$ and
$\Delta_{i,j}$
be the set
\[
%T_{i,j} =
%h+4(j-1)\ln\beta,
%h+4j\ln\beta
\Delta_{i,j}= \L^{n-1}(k)\times
\bigl[ h+(j-1)\ln\beta, %+\ln\beta,
h+j\ln\beta %-\ln\beta
\bigr] \times
\L^{d-n}(\ln\beta).
\]
%
%Notice that $\Delta$ depends on $i,j$, however we drop $i,j$ from the
%notation.
The set $\Delta_{i,j}$ is isometric to a set of the form
\[
\L^{n-1}(k) \times\L^{d-n+1}(\ln\beta).
\]
We conclude that there exist two indices $i,j$ and two times $a,b$ such that:

$\bullet$ $i,j$ are integers and satisfy $0\leq|i|\leq\ln\beta$,
$0\leq j\leq J$.

$\bullet$ $a,b$ are integers and satisfy $0\leq b-a\leq
\tb/J+2$.

$\bullet$ The event $\cR(S_i)$ occurs.

$\bullet$ At time $a$, the set
$\Delta_{i,j}$ has not been touched by a $\STC$ emanating from top or
bottom of $S_i$
in the process
$(\s^{n-1\pm,\xi}_{S_i,t}, t\geq0)$.
We denote this event by $\cT_a$.

$\bullet$
At time $b$,
the set
$\Delta_{i,j}$ is vertically crossed
%before time $b$
in the process
$(\s^{n-1\pm,\xi}_{S_{i},t},t\geq0)$.
%$(\s^{n-1\pm,\xi}_{\Delta_{i,j},t},t\geq0)$.
%by a $\STC$ born after time $a$.
We denote this event by $\cV_b$.

%%Since the event $\cR(S_i)$ occurs as well, there is no large $\STC$
%intersecting $T_{i,j}$,
%hence we have
From the previous discussion, we see that
\[
P \pmatrix{ \L^n(k)\times \L^{d-n}(\ln\beta) \mbox{ is
crossed}
\cr
\mbox{vertically before time } \tau_\beta } \leq \sum
_{i,j} \sum_{a,b} P
\bigl(\cR(S_i), \cT_a, \cV_b\bigr)
\]
with the summation running over indices $i,j,a,b$ satisfying the above
conditions.\vspace*{9pt}

\textit{Step} 5: \textit{Conditioning on the configuration at the time of arrival
of the large $\STC$.}

We next estimate the probability appearing in the summation.
To alleviate the formulas, we drop $i,j$ from the notation, writing
$S,\Delta, \zeta$ instead of
$S_i,\Delta_{i,j},
\zeta_{i,j}$.
For $Q$ an $n$-small parallelepiped, we denote by
$\cR(Q,a)$ the event
\[
\cR(Q,a) = \pmatrix{ \mbox{neither local nucleation nor creation}
\cr
\mbox{of a
large $\STC$ occurs before time $a$}
\cr
\mbox{%$\exp(\beta\k)$
for the process
%in the set $S$ starting from
$\bigl(\s^{n\pm,\xi}_{Q,t},{t\geq0}\bigr)$}}.
%$(\s^{n\pm,\xi}_{Q,t})_{t\geq0}$
\]
%
%For $S$ a parallelepiped,
We define the event
$\cR(S,a)$ as
\[
\cR(S,a) %\exp(\beta\k),D)
= \mathop{\bigcap_{ Q\ n\mbox{-}\mathrm{small}}}_
{ Q\subset S}
\cR(Q,a ),
\]
and we estimate its probability as in
Proposition \ref{nuclei}. For $a\leq\tb$, we obtain that
\begin{eqnarray*}
&&
P\bigl(\cR(S,a)^c\bigr) \\
&&\qquad\leq \SES
+|S | (n\ln\beta)^d \bigl(2^{(n\ln\beta)^d} +{4\beta(m_n+2)^3(n
\ln\beta)^{2dm_n+3d}} \tau_\beta \bigr) \\
&&\qquad\quad\hspace*{28pt}{}\times\phi_n(\beta)
\exp(-\beta\Ga_n) %\big|R_\beta(n\ln\beta)^d\big|^{|I|}
.
\end{eqnarray*}
Since
\[
\limsup_{\b\to\infty} \frac{1}{\beta} \ln \bigl( |S|
\tau_\beta \exp(-\beta\Ga_n) \bigr) < nL_n+
\k_n-\G_n = 0,
\]
we conclude that
\[
\lim_{\b\to\infty} P\bigl(\cR(S,a)\bigr) = 1.
\]
We will next condition on the configuration at time $a$ in $\Delta$
in order to estimate the probability of the event $\cV_b$,
\begin{eqnarray*}
P\bigl(\cR(S), \cT_a, \cV_b\bigr) &=& \sum
_\zeta P\bigl(\cR(S), \cT_a,
\cV_b, \s^{n-1\pm,\xi}_{S,a}|_{\Delta}=\zeta
\bigr)
\\
&\leq& \sum_\zeta P\bigl(\cR(S,a),
\cT_a, \cV_b, \s^{n-1\pm,\xi}_{S,a}|_{\Delta}=
\zeta \bigr). %\sum_\zeta
%P(\cR(S), \cT_a,  \cV_b,
% \leq
%P(\cR(S,a),  \cV_b,
\end{eqnarray*}
Yet the knowledge of the configuration at time $a$ is not enough
to decide whether the event
$\cV_b$ will occur: we need also to take into account the $\STC$
present at time $a$
in $\Delta$
to determine whether a vertical crossing occurs
in $\Delta$
before time $b$.
Thus we record the $\STC$ which are present in
the configuration
$\s^{n-1\pm,\xi}_{S,a}|_{\Delta}$.
%To avoid introducing another heavy notation, we simply write
We write
\[
\s^{n-1\pm,\xi}_{S,a}|_{\Delta}=\zeta,\qquad \STC_\xi
\bigl(\s^{n-1\pm,\xi}_{S,t},0\leq t\leq a \bigr)|_{\Delta
\times\{a\}}=
\STC(\zeta) %,
\]
to express that the configuration in $\Delta$ at time $a$ is $\zeta$
and that the trace
at time $a$
of the $\STC$ created before time $a$ in $\zeta$ is given by $\STC
(\zeta)$.
We condition next on the following information:
\begin{eqnarray*}
&&
\sum_\zeta P\bigl(\cR(S,a), \cT_a,
\cV_b, \s^{n-1\pm,\xi}_{S,a}|_{\Delta}=\zeta
\bigr)
\\
&&\qquad= \sum_{\zeta,\STC(\zeta)} P \pmatrix{ \cR(S,a),
\cT_a, \cV_b, \s^{n-1\pm,\xi}_{S,a}|_{\Delta}=
\zeta,
\cr
%,
\STC_\xi \bigl(\s^{n-1\pm,\xi}_{S,t},0
\leq t\leq a \bigr)|_{\Delta
\times\{a\}}=\STC(\zeta) }
\\
%P
%%,
%$$
&&\qquad = \sum_{\zeta,\STC(\zeta)} P \left(
\cV_b \bigg| \matrix{ \cR(S,a), \cT_a,
\s^{n-1\pm,\xi}_{S,a}|_{\Delta}=\zeta,
\cr
\STC_\xi
\bigl(\s^{n-1\pm,\xi}_{S,t},0\leq t\leq a \bigr)|_{\Delta
\times\{a\}}=
\STC(\zeta) } \right)
\\
&&\hspace*{36pt}\qquad\quad{}\times P \left( \matrix{ \cR(S,a), \cT_a, \s^{n-1\pm,\xi}_{S,a}|_{\Delta}=
\zeta,
\cr
\STC_\xi \bigl(\s^{n-1\pm,\xi}_{S,t},0\leq t
\leq a \bigr)|_{\Delta
\times\{a\}}=\STC(\zeta) } \right). %\cR(S,a), \cT_a,
%in the process
%$\s^{n-1\pm,\xi}_{S_i,.}$.
\end{eqnarray*}
On the event
$\cT_a$,
by Lemma \ref{cresSTC},
\begin{eqnarray*}
\s^{n-1\pm,\xi}_{S,a}|_{\Delta} &=& \s^{n\pm,\xi}_{S,a}|_{\Delta},
\\
\STC_\xi \bigl(\s^{n-1\pm,\xi}_{S,t},0\leq t\leq a
\bigr)|_{\Delta
\times\{a\}}&=& \STC_\xi \bigl(\s^{n\pm,\xi}_{S,t},0
\leq t\leq a \bigr)|_{\Delta
\times\{a\}},
\end{eqnarray*}
whence
\begin{eqnarray*}
&&P\pmatrix{ \cR(S,a), \cT_a, \s^{n-1\pm,\xi}_{S,a}|_{\Delta}=
\zeta,
\cr
\STC_\xi \bigl(\s^{n-1\pm,\xi}_{S,t},0\leq t
\leq a \bigr)|_{\Delta
\times\{a\}}=\STC(\zeta) }
\\
&&\qquad\leq P\pmatrix{ \cR(S,a), %\cT_a,
\s^{n\pm,\xi}_{S,a}|_{\Delta}=
\zeta,
\cr
\STC_\xi \bigl(\s^{n\pm,\xi}_{S,t},0\leq t
\leq a \bigr)|_{\Delta
\times\{a\}}=\STC(\zeta)}.
\end{eqnarray*}
Let us set
\[
\nu(\zeta) = P \bigl( \s^{n\pm,\xi}_{S,a}|_{\Delta}=\zeta |
\cR(S,a) \bigr).
\]
Thus $\nu$ is the law of
the configuration
$\s^{n\pm,\xi}_{S,a}|_{\Delta}$
conditioned on the event
$\cR(S,a)$.
This configuration, denoted by $\zeta$, comes equipped with the trace
of the $\STC$
created before time $a$, which are denoted by $\STC(\zeta)$.
Formally, the law $\nu$ should be a law on the trace of the $\STC$ at
time $a$;
however, to alleviate the text, we make a slight abuse of notation, and
we deal with $\nu$
as if it was a law on the configurations.
%On the event
%$\cR(S,a)\cap\cT_a$, the event $\cV_b$ has to occur between times $a$
%and $b$.
With this convention and using the Markov property,
we rewrite the previous inequalities as
\begin{eqnarray*}
&&
P\bigl(\cR(S), \cT_a, \cV_b\bigr)
\\
&&\qquad
\leq\sum_{\zeta,\STC(\zeta)} %\sum_{\zeta}
P \bigl(
\cV_b | \s^{n-1\pm,\xi}_{S,a}|_{\Delta}=\zeta,
\STC(\zeta) \bigr) % \nu(\zeta,\STC(\zeta))
\nu(\zeta) P\bigl( \cR(S,a)\bigr)
\\
&&\qquad\leq \sum_{\zeta,\STC(\zeta)} P\left(\matrix{ \mbox{there is a
vertical crossing between}
\vspace*{1pt}\cr
\mbox{times $a$ and $b$ in $\bigl(
\s^{n-1\pm,\xi}_{\Delta,t}, t\geq0\bigr)$ %$\exp(\beta\k)$
%before time $b$ by a STC born after time $a$
}} \bigg|
\matrix{ \s^{n-1\pm,\xi}_{S,a}|_{\Delta}=\zeta
\cr
\STC(\zeta)}
\right) % \nu(\zeta,\STC(\zeta))
\nu(\zeta)
\\
&&\qquad\leq %\Delta_{i,j} is vertically crossed
%in the process
%$\s^{(n-1)\pm,\xi}_{S_i,.}$.
\sum
_{\zeta} P\pmatrix{ \mbox{there exists a vertical
crossing in}
\cr
\STC_\zeta \bigl(\s^{n-1\pm,\zeta}_{\Delta,t}, 0
\leq t\leq b-a \bigr)} %\begin{matrix}
%(\s^{(n-1)\pm,\zeta}_{\Delta,t}, t\geq0)
%is vertically
%}\\
%crossed before time $b-a$}
% \nu(\zeta,\STC(\zeta)).
\nu(\zeta).
%in the process
%$\s^{(n-1)\pm,\xi}_{S_i,.}$.
\end{eqnarray*}
We check next that
%$\nu$ satisfies
the hypothesis on the initial law at
rank $n-1$ is satisfied by the law $\nu$ of $\zeta$
and that the
hypothesis on the initial clusters is satisfied
by $\STC(\zeta)$, the $\STC$ present in $\zeta$.\vspace*{9pt}
%by the trace of the $\STC$ at time $a$.
%We distinguish the cases $n=1$ and $n\geq2$.
%If $n=1$, then the event
%$\cR(S,a)$ implies that
%$\s^{1\pm,\xi}_{S,a}$ is equal to $-1$ everywhere on $S$ and
%there is nothing more to check.
%Suppose next that $n\geq2$.

\textit{Step} 6: \textit{Check of the hypothesis on the initial $\STC$ at rank $n-1$}.
Let $\cC$ belong to
$\STC_\xi(\s^{n\pm,\xi}_{S,t},0\leq t\leq a)$. Then $\cC$
is the union of $\STC$ belonging to
$\STC(\s^{n\pm,\xi}_{S,t},0\leq t\leq a)$
and to $\STC(\xi)$.
Since the event $\cR(S,a)$ occurs, any $\cC$ in
$\STC(\s^{n\pm,\xi}_{S,t},0\leq t\leq a)$
has diameter at most $\ln\ln\beta$. Thus any path in $\cC$
having diameter strictly larger than $\ln\ln\beta$
has to meet a $\STC$ of
$\STC(\xi)$.
Suppose there exists $\cC$ in
$\STC_\xi(\s^{n\pm,\xi}_{S,t},0\leq t\leq a)$
such that
\[
\diam\cC \geq \tfrac{1}{4}\ln\beta.
\]
By Lemma \ref{reason}, there would exist $\cC'\subset\cC$
and $a'\leq a$ such that
\[
\cC'\in\STC_\xi\bigl(\s^{n\pm,\xi}_{S,t},0
\leq t\leq a'\bigr),\qquad \tfrac{1}{4}\ln\beta \leq \diam
\cC' \leq \tfrac{1}{2}\ln\beta.
\]
Let $Q'$ be an $n$-small box containing $\cC'$. The previous discussion
implies that $\cC'$ would meet at least
$\frac{1}{4}(\ln\beta)/\ln\ln\beta$ elements of $\STC(\xi)$.
Thus we would have
\[
\mathop{\sum_{\cC\in\STC(\xi)}}_
{\cC\cap Q'\neq\varnothing}\diam\cC \geq
\frac{\ln\beta}{4\ln\ln\beta} > (d-n+1) \ln\ln\beta
\]
and this would contradict the hypothesis on the initial $\STC$ present
in $\xi$.
Therefore any $\STC$ in
%element of
$\STC_\xi(\s^{n\pm,\xi}_{S,t},0\leq t\leq a)$
has a diameter less than
$\frac{1}{4}\ln\beta$.
Let now $Q$ be an
$(n-1)$-small
parallelepiped included in $\Delta$.
Let $Q'$
be an
$n$-small
parallelepiped containing $Q$ and such that
\[
d\bigl(S\setminus Q',Q\bigr) > \tfrac{1}{3}\ln\beta.
\]
From the previous discussion, we see that
%and Lemma \ref{cresSTC},
a $\STC$
of $\STC_\xi(\s^{n\pm,\xi}_{S,t},0\leq t\leq a)$
which intersects the box $Q$ does not meet the inner boundary of $Q'$.
%if $\cC$ is a $\STC$
%of diameter at most $\ln\ln\beta$ intersecting the box $Q$,
%then $\cC$
By Lemma \ref{cresSTC}, such a $\STC$
is also a $\STC$ of the process
$\STC_\xi (\s^{n\pm,\xi}_{Q',t},0\leq t\leq a )$.
It follows that
\[
\mathop{\sum_{\cC\in\STC(\zeta)}}_{
{\cC\cap Q\neq\varnothing}}\diam\cC \leq
\diam \STC_\xi \bigl(\s^{n\pm,\xi}_{Q',t},0\leq t\leq a
\bigr). % \leq
%{\cC\cap Q'\neq\varnothing}}\diam\cC +
% \leq
%(d-n+1)\ln\ln\beta+\ln\ln\beta
% =
%(d-n+2)\ln\ln\beta
\]
Since the event $\cR(S,a)$ occurs, any $\cC$ in
$\STC(\s^{n\pm,\xi}_{S,t},0\leq t\leq a)$
has diameter at most $\ln\ln\beta$.
From the
hypothesis on the initial $\STC$ at rank $n$, we have
%any
%$\STC$ in
%$\STC_\xi(\s^{n\pm,\xi}_{Q',t},0\leq t\leq a)$
%has diameter at most
%
\begin{eqnarray*}
&&\diam\STC_\xi\bigl(\s^{n\pm,\xi}_{Q',t},0\leq t\leq a
\bigr)
\\
&&\qquad\leq\mathop{\sum_{\cC\in\STC(\xi)}}_
{\cC\cap Q'\neq\varnothing}\diam\cC + \diam
\STC \bigl(\s^{n\pm,\xi}_{Q',t},0\leq t\leq a \bigr)
\\
&&\qquad\leq (d-n+1)\ln\ln\beta+\ln\ln\beta = (d-n+2)\ln\ln\beta,
\end{eqnarray*}
and the hypothesis on the initial $\STC$ present in $\zeta$ is
fulfilled.\vspace*{9pt}
%Of course in the last inequalities we relied on
%the hypothesis on the $\STC$ at rank $n$.

\textit{Step} 7: \textit{Check of the hypothesis on the initial law at rank $n-1$.}
Let $(Q_i,\break i\in I)$ be a family of
$(n-1)$-small
parallelepipeds included in $\Delta$ such that
\[
\forall i,j\in I\qquad i\neq j \quad\Rightarrow\quad d(Q_i,Q_j)
> 5(d-n+2)\ln\ln\beta,
\]
and let $(\sigma_i,i\in I)$ be
a family of configurations
in the parallelepipeds $(Q_i,\break i\in I)$.
For $i\in I$, let $R_i$ be the box $Q_i$ enlarged by a distance
$2\ln\ln\beta$
%$\frac{1}{3}\ln\beta$
along the first $n$ axis.
The boxes
$(R_i,i\in I)$ are $n$-small and satisfy
\[
\forall i,j\in I\qquad i\neq j \quad\Rightarrow\quad d(R_i,R_j)
> 5(d-n+1)\ln\ln\beta.
\]
On the event $\cR(S,a)$, we have by Lemma \ref{cresSTC}
\[
\forall i\in I\qquad \s^{n\pm,\xi}_{S,a}|_{Q_i} =
\s^{n\pm,\xi}_{R_i,a}|_{Q_i}.
\]
Therefore
\begin{eqnarray*}
\nu ( \forall i\in I, \sigma|_{Q_i}=\sigma_i ) &=& P \bigl(
\forall i\in I, \s^{n\pm,\xi}_{S,a}|_{Q_i} =
\sigma_i | \cR(S,a) \bigr)
\\
&=& P \bigl( \forall i\in I, \s^{n\pm,\xi}_{R_i,a}|_{Q_i} =
\sigma_i | \cR(S,a) \bigr). %\forall(\sigma_i,i\in I)\in\prod_{i\in I}\{ -1,+1 \}^{Q_i}
\end{eqnarray*}
%
%so that
%$d(Q_i,(R_i)^c) \geq2
We condition next on the initial configurations in the boxes $R_i, i\in I$,
%Next, %for $i\in I$,
%
\begin{eqnarray*}
&&
P \bigl( \cR(S,a), \forall i\in I, \s^{n\pm,\xi}_{R_i,a}|_{Q_i}
=\sigma_i % |
\bigr)
\\
&&\qquad= %\frac{1}{P( \cR(S,s))}
\sum_{\zeta_i, i\in I} P \bigl( \cR(S,a),
\forall i\in I, \s^{n\pm,\xi}_{R_i,a}|_{Q_i} =
\sigma_i, %\s^{n\pm,\xi}_{R_i,0}|_{Q_i}
\xi|_{R_i} =\zeta_i
\bigr)
\\
&&\qquad\leq %\frac{1}{P( \cR(S,s))}
\sum_{\zeta_i, i\in I} P \bigl(\forall i\in I,
\cR(R_i,a), \s^{n\pm,\xi}_{R_i,a}|_{Q_i} =
\sigma_i, \xi|_{R_i} =\zeta_i \bigr)
%i\in I
%=\sigma_i,
%=\zeta_i
\\
&&\qquad= %\frac{1}{P( \cR(S,s))}
\sum_{\zeta_i, i\in I} P \bigl( \forall i\in I,
\cR(R_i,a), \s^{n\pm,\xi}_{R_i,a}|_{Q_i} =
\sigma_i | \forall i\in I, \xi|_{R_i} =\zeta_i
\bigr)
\\
&&\hspace*{22pt}\qquad\quad{}\times P ( \forall i\in I, \xi|_{R_i} =\zeta_i ).
\end{eqnarray*}
We next use the hypothesis on the law of $\xi$ and the fact
that, once the initial configurations in the boxes $R_i$ are fixed, the
dynamics in these boxes with $n\pm$ boundary conditions are
independent. We obtain
\begin{eqnarray*}
&&P \bigl( \cR(S,a), \forall i\in I, \s^{n\pm,\xi}_{R_i,a}|_{Q_i}
=\sigma_i % |
\bigr)
\\
&&\qquad\leq %\frac{1}{P( \cR(S,s))}
\sum_{\zeta_i, i\in I} \prod
_{i\in I} P \bigl( \cR(R_i,a), \s^{n\pm,\xi}_{R_i,a}|_{Q_i}
=\sigma_i | \xi|_{R_i} =\zeta_i \bigr)
\phi_n(\beta) \rho_{R_i}^{n\pm}(\zeta_i).
\end{eqnarray*}
We recall\vspace*{1pt} that
$(\st^{n\pm,\xi}_{R_i,t})_{t\geq0}$
is the process conditioned to stay in $\cR_n(R_i)$.
On the event
$\cR(R_i,a)$, the initial configuration $\zeta_i$ belongs to
$\cR_n(R_i)$ and
\[
\s^{n\pm,\xi}_{R_i,a}|_{Q_i} = \st^{n\pm,\xi}_{R_i,a}|_{Q_i},
\qquad \rho_{R_i}^{n\pm}(\zeta_i) \leq
(m_n+1)|R_i|^{m_n} \mut_{R_i}^{n\pm}(
\zeta_i). %
%|R_i|
% \leq
%(n\ln\beta)^{d},
%(m_n+1)(n\ln\beta)^{dm_n}
\]
Moreover
$|R_i|
\leq
(n\ln\beta)^{d}$ and
$P(\cR(S,a))\geq1/2$ for $\beta$ large enough.
Thus
\begin{eqnarray*}
&&\nu ( \forall i\in I, \sigma|_{Q_i}=\sigma_i )
\\
&&\qquad\leq\frac{1}{P(\cR(S,a))} P \bigl( \cR(S,a), \forall i\in I, \s^{n\pm,\xi}_{R_i,a}|_{Q_i}
=\sigma_i % |
\bigr)
\\
&&\qquad\leq %\frac{1}{P( \cR(S,s))}
2\prod_{i\in I} \biggl( \sum
_{\zeta_i\in
\cR_n(R_i)
} P \bigl( \st^{n\pm,\xi}_{R_i,a}|_{Q_i}
=\sigma_i %, \cR(Q'_i,s)
| \xi|_{R_i} =\zeta_i
\bigr) \phi_n(\beta) \rho_{R_i}^{n\pm}(
\zeta_i) \biggr)
\\
&&\qquad\leq %\frac{1}{P( \cR(S,s))}
2 \prod_{i\in I} \bigl(
(m_n+1) %|R_i|^{m_n}
(n\ln\beta)^{dm_n} \phi_n(
\beta) \mut_{R_i}^{n\pm} ( \sigma|_{Q_i} =
\sigma_i %, \cR(Q'_i,s)
) \bigr). %\rho_{Q'_i}^{n\pm}(\xi_i)
% =
%P\big( \forall i\in I
%=\sigma_i |
%=\zeta_i
% \big)\\
%P\big( \forall i\in I
%=\zeta_i
% \big)
\end{eqnarray*}
%
%Yet $|R_i|\leq(n\ln\beta)^d$;

\vspace*{9pt}

\textit{Step} 8: \textit{Comparison of
$\mt_{R}^{n\pm}|_{Q}$
and $\rho_{Q}^{n-1\pm}$.}
To conclude we need to prove that
if $Q,R$ are two parallelepipeds which are $n$-small and such that
$Q\subset R$,
then for any configuration $\eta$ in $Q$,
\[
\mt_{R}^{n\pm} ( \sigma|_{Q} =\eta ) \leq
\phi_{n-1}(\beta) \rho_{Q}^{n-1\pm}(\eta),
\]
where $\phi_{n-1}(\beta)$
is a function depending only upon $\beta$
which is $\exp o(\beta)$. This is the purpose of the next three lemmas.
%
%le6.8 #&#
\begin{lemma}
\label{control}
Let $R$ be an $n$-small parallelepiped.
There exists $h_0>0$ such that, for $h\in\ ]0,h_0[$, the following
result holds.
If $\sigma$ is a configuration in $R$ satisfying
\[
|\sigma| \leq m_{n},\qquad H^{n\pm}_{R}(\sigma)
%-H^{n\pm}_{Q}(\underline{\minus})
\leq \Ga_{n-1},
\]
then $|\sigma|
\leq
m_{n-1}$.
\end{lemma}
\begin{pf}
Let $\sigma$ be a configuration satisfying the hypothesis of the lemma,
and let us set
$m=|\sigma|$.
By Lemma \ref{prok}, there exists an $n$-dimensional configuration
$\rho$ such that
$|\rho|=m$ and
$H_{\Z^n}(\rho)
=H^{n\pm}_{R}(\sigma)$. We apply next the simplified
isoperimetric
inequality
stated in Section~\ref{sie},
\begin{eqnarray*}
H_{\Z^n}(\rho) &=& \per(\rho)-h|\rho|
\\
&\geq& \inf \bigl\{ \per(A)\dvtx A\mbox{ is the finite union of $m$ unit cubes} \bigr\}
- hm
\\
&\geq& 2n m^{(n-1)/{n}}-hm.
\end{eqnarray*}
Therefore the number $m$ of pluses in $\sigma$ satisfies
\[
m \leq m_{n},\qquad 2n m^{(n-1)/{n}}-hm \leq \Ga_{n-1}.
\]
Thus, for $h\leq1$,
\[
m \leq \bigl(l_c(n)+1 \bigr)^n \leq \biggl(
\frac{2(n-1)}{h}+1 \biggr)^n \leq \biggl(\frac{2n-1}{h}
\biggr)^n,
\]
whence
\[
2n m^{(n-1)/{n}}-hm = m^{(n-1)/{n}} \bigl(2n-hm^{1/n} \bigr) \geq
m^{(n-1)/{n}},
\]
and we conclude that
\[
m^{(n-1)/{n}} \leq \Ga_{n-1}.
\]
We have the following expansions as
$h\to0$:
\[
m_n \sim \biggl( \frac{2 (n-1)}{h} \biggr)^n,\qquad
\G_n \sim2 \biggl( \frac{2 (n-1)}{h} \biggr)^{n-1}.
\]
Thus,
%there exists a constant $c>0$ such that,
for $h$ small enough,
\[
m^{(n-1)/{n}} \leq \Ga_{n-1} \leq (2n)^{n-1}h^{-(n-2)},
\]
whence
\[
m \leq (2n)^nh^{- {n(n-2)}/({n-1})} \leq m_{n-1},
\]
the last inequality being valid for $h$ small enough, since
$n(n-2)<(n-1)^2$ and
$m_{n-1}$ is of order $h^{-(n-1)}$ as $h$ goes to $0$.
\end{pf}
%
%le6.9 #&#
\begin{lemma}
\label{cutting}
Let $Q\subset R$ be two $n$-small parallelepipeds.
%and let $\pi$ be a half-space.
%orthogonal to one of the first $n$ axis.
If $\eta$ is a configuration in $R$ satisfying
$|\eta|\leq m_{n-1}$, then
\[
H^{n\pm}_{R}(\eta) \geq H^{n\pm}_{Q}(
\eta|_Q).
\]
\end{lemma}
\begin{pf}
We will prove the following intermediate result.
If $\pi$ is a half-space, then
\[
H^{n\pm}_{R}(\eta) \geq H^{n\pm}_{R\cap\pi}(
\eta\cap\pi).
\]
Re\-pea\-ted applications of the above inequality will yield the result stated
in the lemma.
We consider first the case where
$\pi$ is
orthogonal to one of the first $n$ axis, say the $n$th, and it has for equation
\[
\pi = \bigl\{ x=(x_1,\ldots,x_n,\ldots,x_d)\dvtx x_n
\leq h+1/2 \bigr\},
\]
where $h\in\Z$.
We think of $\eta$ as the union of $(d-1)$-dimensional configurations
which are obtained by intersecting $\eta$ with the layers
\[
L_i = \bigl\{ x=(x_1,\ldots,x_d)\in
\Z^d\dvtx  i-\tfrac{1}{2}\leq x_{n}<i+\tfrac{1}{2}
\bigr\},\qquad i\in\Z.
\]
Let us define the hyperplanes
\[
P_i = \bigl\{ x=(x_1,\ldots,x_d)\in
\Z^d\dvtx  x_{n}=i+ \tfrac{1}{2} \bigr\},\qquad i\in\Z.
\]
We have
\[
H^{n\pm}_{R}(\eta) = \sum_{i}
H_{\Z^{d-1}}^{n-1\pm}(\eta\cap L_i) + \sum
_{i} \operatorname{area}( \partial\eta\cap P_i).
%2 \mbox{area}(\mbox{proj}_{n}(\eta))
\]
%
%where
%$\mbox{proj}_{n}(\eta)$
%is the orthogonal projection of $\eta$ on
%$\Z^{n-1}\times\{ 0 \}\times\Z^{d-n}$.
Yet, for any $i>h$, we have
$ |\eta\cap L_i |\leq m_{n-1}$
whence
\[
H_{\Z^{d-1}}^{n-1\pm}(\eta\cap L_i) \geq0.
\]
Moreover
\[
\sum_{i\geq h} \operatorname{area}( \partial\eta\cap
P_i) \geq | \eta\cap L_h|. % \mbox{area}( \eta\cap P_h)-
\]
%
% \subset
%.\]
This is because the boundary conditions are minus on the faces
orthogonal to the
$n$th axis, hence there must be at least one unit interface above each
plus site of
the layer $L_h$.
We conclude that
\begin{eqnarray*}
H^{n\pm}_{R}(\eta) &\geq& \sum_{i\leq h}
H_{\Z^{d-1}}^{n-1\pm}(\eta\cap L_i) + \sum
_{i< h} \operatorname{area}( \partial\eta\cap P_i) +
% \mbox{area}( \eta\cap P_h)
| \eta\cap L_h| %-\mbox{area}( \partial\eta\cap P_h)\\
\\
&=& %\sum_{i\leq h}
% \mbox{area}( \partial\eta\cap P_i)
% +
% \mbox{area}( \eta\cap P_h)-
% =
H^{n\pm}_{R\cap\pi}(\eta\cap\pi)
\end{eqnarray*}
as requested.
The case where $\pi$ is orthogonal to one of the last $d-n$ axis
can be handled similarly.
This case is even easier because the boundary conditions become plus
along $\pi$ and
contribute to lowering the energy.
\end{pf}
%
%le6.10 #&#
\begin{lemma}
%Let $Q\subset R$ be two parallelepipeds which are $n$-small.
Let $Q,R$ be two parallelepipeds which are $n$-small and such that
$Q\subset R$.
If $\eta\in\cR_{n-1}(Q)$,
then
\[
\mt_{R}^{n\pm} ( \sigma|_{Q} =\eta ) \leq
(m_n+1) (n\ln\beta)^{dm_n}\exp \bigl(-\beta
H^{n-1\pm}_{Q}(\eta) %-H^{(n-1)\pm}_{Q}(\underline{\minus})
\bigr).
\]
If $\eta\notin
%{\cal R}_{n}(Q)\setminus
{\cal R}_{n-1}(Q)$, then
\[
\mt_{R}^{n\pm} ( \sigma|_{Q} =\eta ) \leq
(m_n+1) (n\ln\beta)^{dm_n} \exp(-\beta\Ga_{n-1}).
\]
%
% \mbox{ otherwise }
\end{lemma}
\begin{pf}
For any configuration $\eta$ in $Q$,
%${\cal R}_{n-1}(Q)$,
%
\begin{eqnarray*}
&&
\mt_{R}^{n\pm} ( \sigma|_{Q} =\eta ) \\
&&\qquad\leq
\mathop{\sum_{
\rho\in\cR_n(R)}}_{\rho|_{Q}=\eta}
\mt_{R}^{n\pm}(\rho)
\\
&&\qquad\leq \bigl|\cR_n(R) \bigr| \max \bigl\{ \mt_{R}^{n\pm}(
\rho)\dvtx  \rho\in\cR_n(R), \rho|_{Q}=\eta \bigr\}
\\
&&\qquad\leq (m_n+1) (n\ln\beta)^{dm_n} \exp \bigl(-\beta \min \bigl
\{ H_{R}^{n\pm}(\rho)\dvtx  \rho\in\cR_n(R),
\rho|_{Q}=\eta \bigr\} \bigr).
\end{eqnarray*}
%
%where the sum and the max are taken over the configurations $\eta$
%in $Q'$ such that
%$\eta|_{Q}=\s$.
If the minimum in the exponential
is larger than or equal to $\Ga_{n-1}$, then we have the desired
inequality. Suppose that
the minimum is less than $\Ga_{n-1}$.
Let $\rho\in\cR_n(R)$ be such that
$H^{n\pm}_{R}(\rho)\leq\Ga_{n-1}$
and
$\rho|_{Q}=\eta$. By Lemma \ref{control}, we have then
also $|\rho| \leq m_{n-1}$.
Let ${\mathcal C}(\rho)$ be the set of the connected components of
$\rho$.
Since
$\rho\in\cR_n(R)$, we have
\[
\forall C \in{\mathcal C}(\rho)\qquad H^{n\pm}_{R}(C)\geq0
\]
hence
\[
H^{n\pm}_{R}(\rho) \geq \mathop{\sum
_{ C\in{\mathcal C}(\rho)}}_
{ C\cap Q\neq\varnothing} H^{n\pm}_{R}(C).
\]
Let $C\in{\mathcal C}(\rho)$ be such that
${C\cap Q\neq\varnothing}$.
%Re\-pea\-ted applications of
Since $|\rho|\leq m_{n-1}$,
Lemma \ref{cutting} yields that
\[
H^{n\pm}_{R}(C) \geq H^{n\pm}_{Q}(C\cap
Q),
\]
whence
\begin{eqnarray*}
H^{n\pm}_{R}(\rho) &\geq& \mathop{\sum
_{C\in{\mathcal C}(\rho)}}_
{ C\cap Q\neq\varnothing} H^{n\pm}_{Q}(C\cap
Q)
\\
&=& H^{n\pm}_{Q}(\rho\cap Q) = H^{n\pm}_{Q}(
\eta) \geq H^{n-1\pm}_{Q}(\eta).
\end{eqnarray*}
The last inequality is a consequence of the attractivity of the
boundary conditions.
It follows that
$H^{n-1\pm}_{Q}(\eta)\leq\Ga_{n-1}$ so that $\eta$ belongs to $\cR
_{n-1}(Q)$. In addition,
we conclude that
\[
\min \bigl\{ H_{R}^{n\pm}(\rho)\dvtx  \rho\in
\cR_n(R), \rho|_{Q}=\eta \bigr\} \geq H^{n-1\pm}_{Q}(
\eta).
\]
%
%Moreover
%the minimum in the exponential
%is larger than or equal to
%$H^{n-1\pm}_{Q}(\eta)$
which yields the desired inequality.
%=\eta  \big) \leq
%(m_n+1)
%(n\ln\beta)^{dm_n}\exp\big(-\beta
%H^{n-1\pm}_{Q}(\eta)
%%-H^{(n-1)\pm}_{Q}(\underline{\minus})
%as claimed.
\end{pf}

\textit{Step} 9:
\textit{Reduction to a box $\Phi$
of side length of order $\exp(\beta L_{n-1})$.}
Thus the measure $\nu$ on the configurations in $\Delta$ satisfies
the initial hypothesis at rank $n-1$.
Let us set
$\tau'_\beta=b-a$. We have then
\[
\limsup_{\b\to\infty} \frac{1}{\beta} \ln \tau'_\beta
< \k_n-L_n = \k_{n-1}.
\]
We are in position to apply
%the result of Theorem \ref{T2}
the induction hypothesis
at rank $n-1$.
We define
%the length
%l_n(\beta) =
% \mbox{if }n=1 \\
% \mbox{if }n\geq2
%and
the box
\[
\Phi_0 = %\L^{n-1}\big(2l_n(\beta)\big)
%2\exp({\b L_{n-1}})
%) \times\L^{d-n+1}(\ln\beta)
%
\cases{ \L^{1}(2\ln\ln\beta) \times\L^{d-1}(\ln\beta), &\quad
if $n=1$,
\vspace*{1pt}\cr
\L^{n-1}\bigl( 2\exp({\b L_{n-1}}) \bigr) \times
\L^{d-n+1}(\ln\beta), &\quad if $n\geq2$.}
\]
Up to a $\SES$ event, there is
no space--time cluster of diameter
larger than
%$l_n(\beta)$
%
\[
%l_n(\beta) =
%
\cases{ \ln\ln\beta, &\quad if $n=1$,
\vspace*{1pt}\cr
\exp({\b L_{n-1}}), &\quad if $n\geq2$,}
\]
%
%before time
in
\[
\STC_{\zeta}\bigl(\s^{n-1\pm,\zeta}_{\Delta,t}, 0\leq t\leq
\tau'_\beta\bigr).
\]
It follows that
any $\STC$ of the above process is included
in a translate of the box
$\Phi_0$
and the vertical crossing of $\Delta$ can only occur in such a set.
Thus
\begin{eqnarray*}
&&\sum_\zeta P\pmatrix{ \mbox{there is a vertical
crossing in}
\vspace*{1pt}\vspace*{1pt}\cr
\STC_{\zeta}\bigl(\s^{n-1\pm,\zeta}_{\Delta,t}, 0
\leq t\leq\tau '_\beta\bigr) %\begin{matrix}
%(\s^{(n-1)\pm,\zeta}_{\Delta,t}, t\geq0)
%is vertically
%}\\
%crossed before time $\tau'_\beta$}
} \nu(\zeta)
\\
&&\qquad\leq\sum_\Phi \sum_\zeta
P\pmatrix{ %\begin{matrix}
%(\s^{(n-1)\pm,\zeta}_{\Phi,t}, t\geq0)
%is vertically
%}\\
%%$\exp(\beta\k)$
%crossed before time $\tau'_\beta$}
\mbox{there is a vertical
crossing in}
\vspace*{1pt}\vspace*{1pt}\cr
\STC_{\zeta}\bigl(\s^{n-1\pm,\zeta}_{\Phi,t}, 0
\leq t\leq\tau'_\beta\bigr)} %\begin{matrix}
%(\s^{(n-1)\pm,\zeta}_{\Delta,t}, t\geq0)
%is vertically
%}\\
%crossed before time $\tau'_\beta$}
\nu(\zeta) +\SES,
\end{eqnarray*}
where the sum over $\Phi$ runs over the translates of $\Phi_0$
included in $\Delta$.
%Let us fix one such translate.
We estimate
\[
\sum_\zeta P\pmatrix{ %\begin{matrix}
%(\s^{(n-1)\pm,\zeta}_{\Phi,t}, t\geq0)
%is vertically
%}\\
%%$\exp(\beta\k)$
%crossed before time $\tau'_\beta$}
\mbox{there is a vertical crossing in}
\vspace*{1pt}\vspace*{1pt}\cr
\STC_{\zeta}\bigl(
\s^{n-1\pm,\zeta}_{\Phi,t}, 0\leq t\leq\tau'_\beta
\bigr) } \nu(\zeta)
\]
for $\Phi=x+\Phi_0$ a fixed translate of $\Phi_0$.\eject%\vspace*{9pt}

\textit{Step} 10:
\textit{Reduction to boxes $\Phi_i\subset\Phi$
of vertical side length of order $\ln\beta/\ln\ln\beta$.}
We consider the following subsets of
$\Phi$.
Let us set $I=\ln\ln\beta$.
If $n=1$, then we define
for $1\leq i\leq I$
\[
\Phi_i =
x+\L^{1}( 2\ln\ln\beta) \times
\biggl[ \frac{i\ln\beta}{2\ln\ln\beta}, \frac{ (i+1/2)
\ln\beta}{2\ln\ln\beta} \biggr] \times\L^{d-2}(\ln
\beta).
\]
%
% \mbox{if }n=1 \\
If $n=2$, then we define
for $1\leq i\leq I$
\[
\Phi_i = x+\L^{n-1}\bigl( 2\exp({\b L_{n-1}})
\bigr) \times \biggl[ \frac{i\ln\beta}{2\ln\ln\beta}, \frac{ (i+1/2)
\ln\beta}{2\ln\ln\beta} \biggr] \times
\L^{d-n}(\ln\beta). % \mbox{if }n\geq2
\]
%
%$I=\frac{1}{2}\ln\ln\beta$.
These sets are pairwise disjoint and satisfy, for $\beta$ large enough,
\[
\forall i,j\leq I\qquad %\frac{1}{2}\ln\ln\beta
i\neq j \quad\Rightarrow\quad d(\Phi_i,
\Phi_j) \geq \frac{\ln\beta}{4\ln\ln\beta} > 5(d-n+2)\ln\ln\beta.
\]
If the set $\Phi$
endowed with $n-1\pm$ boundary conditions
is vertically crossed
before time $\tau'_\beta$,
so are the sets
$\Phi_i,  1\leq i\leq I$.
The vertical side of $\Phi_i$ is
%larger than
%
\[
\frac{\ln\beta}{4\ln\ln\beta} > (d-n+3)\ln\ln\beta,
\]
and hence the vertical crossing of $\Phi_i$ implies that
a $\STC$ of diameter larger than
$(d-n+3)\ln\ln\beta$
has been created in $\Phi_i$.\vspace*{9pt}

\textit{Step} 11: \textit{Conclusion of the induction step.}
By Lemma \ref{reason}, there exists
an $(n-1)$-small box $Q_i$ included in $\Phi_i$
and
a $\STC$
$\cC'_i$ in
$\STC_\zeta(\s^{n-1\pm,\zeta}_{Q_i,t}, 0
\leq t\leq
\tau'_\beta
)$
such that
\[
\diam\cC'_i \geq (d-n+3)\ln\ln\beta.
\]
Taking into account the
hypothesis on the initial $\STC$ present in $\zeta$,
%$(\s^{n\pm,\xi}_{Q',t}, 0
%)$.
%Yet
%
\begin{eqnarray*}
\diam\cC'_i &\leq& \diam \STC_\zeta \bigl(
\s^{n-1\pm,\zeta}_{Q_i,t}, 0 \leq t\leq \tau'_\beta
\bigr)
\\
&\leq& \mathop{\sum_{\cC\in\STC(\zeta)}}_
{\cC\cap Q_i\neq\varnothing}\diam\cC +
\diam \STC\bigl( \s^{n-1\pm,\zeta}_{Q_i,t}, 0\leq t\leq
\tau'_\beta \bigr)
\\
&\leq& (d-n+2)\ln\ln\beta+ \diam \STC\bigl( \s^{n-1\pm,\zeta}_{Q_i,t}, 0
\leq t\leq \tau'_\beta \bigr).
\end{eqnarray*}
Therefore
\[
\diam \STC\bigl( \s^{n-1\pm,\zeta}_{Q_i,t}, 0\leq t\leq
\tau'_\beta \bigr) \geq \ln\ln\beta
\]
and a large STC is formed in the process
$(\s^{n-1\pm,\zeta}_{\Phi_i,t}, 0\leq t\leq\tau'_\beta)$.
%If we fix the initial configuration $\zeta$ in $\Phi$ and
%If we set
%$\zeta_i=\zeta|_{\Phi_i}$ for $i\in I$,
We have thus
\begin{eqnarray*}
&&
P\pmatrix{ %\begin{matrix}
%(\s^{(n-1)\pm,\zeta}_{\Phi,t}, t\geq0)
%is vertically
%}\\
%%$\exp(\beta\k)$
%crossed before time $\tau'_\beta$}
\mbox{there is a
vertical crossing in}
\vspace*{1pt}\vspace*{1pt}\cr
\STC_{\zeta}\bigl(\s^{n-1\pm,\zeta}_{\Phi,t},
0\leq t\leq\tau'_\beta\bigr) }
\\
&&\qquad\leq % \nu(\zeta)
P\pmatrix{ \mbox{each set $\Phi_i$ is vertically
crossed}
\vspace*{1pt}\vspace*{1pt}\cr
\mbox{ %$\exp(\beta\k)$
in $\bigl(\s^{n-1\pm,\zeta}_{\Phi,t}, 0
\leq t\leq\tau'_\beta\bigr)$} }
\\
&&\qquad\leq P\pmatrix{ \mbox{for each $i\in I$, a large $\STC$ is formed}
\vspace*{1pt}\vspace*{1pt}\cr
\mbox{%
%$\exp(\beta\k)$
in the process $\bigl(\s^{n-1\pm,\zeta}_{\Phi_i,t}, 0\leq t\leq
\tau'_\beta\bigr)$}}
\\
&&\qquad\leq P\pmatrix{ \mbox{for the process $\bigl(\s^{n-1\pm,\zeta}_{\Phi,t},
0\leq t\leq \tau'_\beta\bigr)$}
\vspace*{1pt}\vspace*{1pt}\cr
\mbox{$\ln\ln\beta$
large $\STC$ are created in $(n-1)$-small}
\vspace*{1pt}\vspace*{1pt}\cr
\mbox{%$\exp(\beta\k)$
% before time $\tau_\beta$
parallelepipeds which are pairwise at} %included in $\Phi$
\vspace*{1pt}\vspace*{1pt}\cr
\mbox{distance larger
than $5(d-n+2)\ln\ln\beta$}}.
\end{eqnarray*}
Since $\nu$ satisfies the hypothesis on the initial law at rank $n-1$ and
the volume $\Phi$ and the time $\tau'_\beta$ satisfy
\[
\limsup_{\b\to\infty} \frac{1}{\beta} \ln |\Phi| \leq
(n-1)L_{n-1},\qquad \limsup_{\b\to\infty} \frac{1}{\beta} \ln
\tau'_\beta < \k_{n-1},
\]
we can apply Proposition \ref{nuclei}
to conclude that
\[
\sum_\zeta P\pmatrix{ \mbox{there is a vertical
crossing in}
\vspace*{1pt}\vspace*{1pt}\cr
\STC_{\zeta}\bigl(\s^{n-1\pm,\zeta}_{\Phi,t}, 0
\leq t\leq\tau'_\beta\bigr) } %\begin{matrix}
%(\s^{n-1\pm,\zeta}_{\Phi,t}, t\geq0)
%is vertically
%}\\
%%$\exp(\beta\k)$
%crossed before time $\tau'_\beta$}
\nu(\zeta)
\]
is $\SES$. Coming back along the chain of inequalities, we see that
\[
P\bigl(\cR(S), \cT_a, \cV_b\bigr)
\]
is also $\SES$, as well as
\[
\sum_{i,j} \sum_{a,b}
P\bigl(\cR(S_i), \cT_a, \cV_b\bigr)
\]
since the number of terms in the sums is of order exponential in $\beta$.
Coming back one more step, we obtain that
\[
P \bigl(\exists \cC\in\STC_\xi(0, \tau_\beta )\mbox{ with
}\diam\cC\geq \exp(\beta L_n) \bigr)
\]
is also $\SES$, as required.
%since the processes in disjoint regions are independent, once the
%initial configuration
%is fixed.

%We check that
%H^{n\pm}_{Q}(\sigma)
%%-H^{n\pm}_{Q}(\underline{\minus})
%:
%|\eta|> m_{n-1}
% >
%s6.5 #&#
\subsection{\texorpdfstring{Proof of the lower bound in Theorem \protect\ref{mainfv}}
{Proof of the lower bound in Theorem 1.2}}
\label{concl}
For technical convenience, we consider here boxes of
side length $c\exp(\beta L)$. The statement of Theorem \ref{mainfv}
corresponds to the special case where $c=1$.
Let $L,c>0$ and let $\L_\beta=\L(c\exp(\beta L))$ be a cubic
box of side length $c\exp(\beta L)$.
Let $\kappa$ be such that
\[
\kappa<\max(\Gamma_d-dL,\kappa_d)
\]
and let $\tau_\beta=\exp(\beta\k)$.
We have
\[
\P \bigl( \s_{\L_\beta,\tau_\beta}^{-,\minus}(0)=1 \bigr) = \P\pmatrix{ ( 0,\tb)
\mbox{ belongs to a nonvoid $\STC$}
\vspace*{1pt}\cr
\mbox{of the process $\bigl(
\s_{\L_\beta,t}^{-,\minus}, 0\leq t\leq\tau_\beta\bigr)$}}.
\]
Let us denote by $\cC^*$ the $\STC$
of the process
$(\s_{\L_\beta,t}^{-,\minus}, 0\leq t\leq\tau_\beta)$
containing the space--time point
$( 0,\tb)$.
In case
$\s_{\L_\beta,\tau_\beta}^{-,\minus}(0)=-1$, then
$\cC^*=\varnothing$.
We write then
\begin{eqnarray*}
&&\P \bigl( \s_{\L_\beta,\tau_\beta}^{-,\minus}(0)=1 \bigr)
\\
&&\qquad=\P \bigl(\cC^*\neq\varnothing, \diam\cC^*<\ln\ln\beta \bigr) + \P \bigl(\diam
\cC^*\geq\ln\ln\beta \bigr).
\end{eqnarray*}
By Lemma \ref{cresSTC},
if $\diam\cC^*<\ln\ln\beta$, then
%the STC
$\cC^*$
is also a $\STC$
of the process
$(\s_{\L(\ln\beta),t}^{-,\minus}, 0\leq t\leq\tau_\beta)$.
Thus
\[
\P \bigl(\cC^*\neq\varnothing, \diam\cC^*<\ln\ln\beta \bigr) \leq \P \bigl(
\s_{\L(\ln\beta),\tau_\beta}^{-,\minus}(0)=1 \bigr).
\]
We use
the processes
$(\s^{d\pm,\mt}_{ \L(\ln\beta),t},t\geq0)$ and
$(\st^{d\pm,\mt}_{ \L(\ln\beta),t},t\geq0)$
to estimate the last quantity,
\begin{eqnarray*}
\P \bigl(\s_{\L(\ln\beta),\tau_\beta}^{-,\minus}(0)=1 \bigr) &\leq& P \pmatrix{
\mbox{nucleation occurs before $\tau_\beta$}
\vspace*{1pt}\vspace*{1pt}\cr
\mbox{in the process }
\bigl(\s_{\L(\ln\beta),t}^{-,\minus}, t\geq0\bigr)}
\\
&&{}+ P \pmatrix{ \s_{\L(\ln\beta),\tau_\beta}^{-,\minus}(0)=1\mbox{,  nucleation does
not occur}
\vspace*{1pt}\vspace*{1pt}\cr
\mbox{before $\tau_\beta$ in the process } \bigl(
\s_{\L(\ln\beta),t}^{-,\minus}, t\geq0\bigr) }
\\
&\leq& P \pmatrix{ \mbox{nucleation occurs before $\tau_\beta$}
\vspace*{1pt}\vspace*{1pt}\cr
\mbox{in the process } \bigl(\s_{\L(\ln\beta),t}^{d\pm,\mt}, t\geq0\bigr)} \\
&&{} + P
\bigl(\st_{\L(\ln\beta),\tb}^{d\pm,\mt}(0)=1 \bigr).
\end{eqnarray*}
Thanks to Lemma \ref{fugaup}, the first term is
exponentially small in $\beta$. The second term is less than
$\mt_{\L(\ln\beta)}^{d\pm}(\s(0)=1)$
which is also
exponentially small in $\beta$.
%$m_d(\ln\beta)^{dm}\exp(-\beta(2d-h))$.
It remains to estimate
\[
\P \bigl(\diam\cC^*\geq\ln\ln\beta \bigr).
\]
We distinguish two cases.

$\bullet$ $L>L_d$. In this case, we write
\begin{eqnarray*}
&&\P \bigl(\diam\cC^*\geq\ln\ln\beta \bigr)
\\
&&\qquad=\P \bigl(\ln\ln\beta\leq\diam\cC^*\leq \exp(\beta L_d) \bigr) + \P
\bigl(\diam\cC^*> \exp(\beta L_d) \bigr).
\end{eqnarray*}
We estimate separately each term. First
\begin{eqnarray*}
&&\P \bigl(\diam\cC^*> \exp(\beta L_d) \bigr)
\\
&&\qquad\leq P \pmatrix{ %\mbox{for }i\in I,
\mbox{the process } \bigl(\s_{\L_\beta,t}^{-,\minus},
0\leq t\leq\tb\bigr) \mbox{ creates}
\vspace*{1pt}\cr
\mbox{a $\STC$ of diameter larger than $
\exp({\b L_d})$}},
\end{eqnarray*}
which is $\SES$ by
Theorem \ref{T2}.
Second,
%letting
%$S_\beta=\L(3\exp(\beta L_d))$,
we have by Lemma \ref{cresSTC},
\begin{eqnarray*}
&&\P \bigl(\ln\ln\beta\leq\diam\cC^*\leq \exp(\beta L_d) \bigr)
\\
&&\qquad\leq P \pmatrix{ %\mbox{for }i\in I,
\mbox{a large $\STC$ is created before time $\tb$ in}
\vspace*{1pt}\cr
\mbox{the process } \bigl(\s_{
%S_\beta
\L(3\exp(\beta L_d)),t}^{-,\minus}, t\geq0\bigr)}.
\end{eqnarray*}
We have reduced the problem to the second case, which we handle next.

$\bullet$ $L\leq L_d$. In this case, we write, with the help of
Lemma \ref{cresSTC},
\begin{eqnarray*}
\P \bigl(\diam\cC^*\geq\ln\ln\beta \bigr) &\leq& %\\
P \pmatrix{
\mbox{a large $\STC$ is created before $\tb$}
\vspace*{1pt}\cr
\mbox{in the
process } \bigl(\s^{-,\minus}_{
%S_\beta
\L_\beta,t}, t\geq0\bigr) }
\\
&\leq& \mathop{\sum_{ Q \ d\mbox{-}\mathrm{small}}}_
{ Q\subset\L_\beta} P
\pmatrix{ %\mbox{for }i\in I,
\mbox{a large $\STC$ is created before $\tb$}
\vspace*{1pt}\cr
\mbox{in
the process } \bigl(\s^{-,\minus}_{
%S_\beta
Q,t}, t\geq0\bigr)}.
\end{eqnarray*}
This inequality holds because the first large STC has to be created in
a $d$-small box,
by Lemma \ref{reason}.
Finally, the term inside the summation is estimated as follows:
\begin{eqnarray*}
&&P \pmatrix{ %\mbox{for }i\in I,
\mbox{a large $\STC$ is created before $\tb$}
\vspace*{1pt}\cr
\mbox{in the process } \bigl(\s^{-,\minus}_{
%S_\beta
Q,t}, t\geq0\bigr)}
\\
&&\qquad\leq P \pmatrix{ %\mbox{for }i\in I,
\mbox{a large $\STC$ is created before nucleation}
\vspace*{1pt}\cr
\mbox{in the process } \bigl(\s^{-,\minus}_{
%S_\beta
Q,t}, t\geq0\bigr)}
\\
&&\qquad\quad{} + P \pmatrix{ %\mbox{for }i\in I,
\mbox{nucleation occurs before $\tb$}
\vspace*{1pt}\cr
\mbox{in the
process } \bigl(\s^{-,\minus}_{
%S_\beta
Q,t}, t\geq0\bigr)}.
\end{eqnarray*}
By Theorem \ref{totcontrole} applied with $\cD=\cR_d(Q)$,
the first term of the right-hand side is $\SES$.
By Lemma \ref{fugaup}, the second term is less than
\[
4\b(m_d+2)^2 |Q|^{2m_d+2}\tau_\beta
\exp(-\b\G_{d}) + \SES,
\]
%
%|Q|^{m_d+1}\tau_\beta
% \exp(-\b\G_{d})
% + \SES
whence
\[
\P \bigl(\diam\cC^*\geq\ln\ln\beta \bigr) \leq |\L_\beta | 4\b
%|Q|^{m_d+1}
(d\ln\beta)^{d(2m_d+4)} \tau_\beta \exp(-\b
\G_{d}) + \SES.
\]
It follows that
\[
\limsup_{\b\to\infty} \frac{1}{\beta} \ln \P \bigl(\diam\cC^*\geq
\ln\ln\beta \bigr) \leq dL+\k-\G_d < 0,
\]
and we are done!
%where the sum runs over all the possible choices of boxes

%s7 #&#
\section{The relaxation regime}
\label{relaxa}
In this section, we prove the upper bound on the relaxation time stated
in Theorem \ref{mainfv}.
This part is considerably easier than the lower bound. The argument
relies on the construction
of an infection process, as done by
Dehghanpour and Schonmann
\cite{DS1} in dimension two, together with an induction on the
dimension and
a simple computation involving the associated growth model \cite{CM2}.
Let us give a quick outline of the structure of the proof.
To each site of the lattice, we associate the box of side length
$\ln\beta$ centered at $x$. A site becomes infected once all the spins
in the associated box are equal to $+1$. The site remains infected as
long as the associated box contains less than $2\ln\ln\beta$ minus
spins (Section~\ref{inucl}).
We give a lower bound for the probability of a site becoming infected, and
this corresponds to a nucleation event. We estimate the probability
that a neighbor of an infected site becomes infected, and
this corresponds to the spreading of the infection
(Section~\ref{ispread}).
Finally,
we define a simple scenario for the invasion of a box of side length
$\exp(\beta L)$, starting from a single infected site
(Section~\ref{iinva}).
We combine all these estimates, and we obtain the required upper bound
on the relaxation time.
%s7.1 #&#
\subsection{The infection process}
\label{inucl}
%Let $\L_\beta=\L(\exp(\beta L))$ be a cubic box of side length
Let
$\L(\exp(\beta L))$
be a cubic box of side length
$\exp(\beta L)$.
Following the strategy of
Dehghanpour and Schonmann
\cite{DS1}, we define a renormalized process $(\mu_t)_{t\geq0}$ on
%$\L_\beta$
$\L(\exp(\beta L))$
as follows.
For $x\in
\L(\exp(\beta L))
$, we set
\[
\L_x = x+\L^d(\ln\beta),
\]
and we define $T_x$ to be the first time when all the spins
of the sites of the box $\L_x$ are equal to $+1$ in the process
$(\s_{
\L(\exp(\beta L)),t}^{-,\minus})_{t\geq0}$,
\[
T_x = \inf \bigl\{ t\geq0\dvtx  \forall y\in\L_x,
\s_{
\L(\exp(\beta L)),t}^{-,\minus}(y)=+1 \bigr\}.
\]
For $\L$ a box, we define the set $\cE(\L)$ to be the set of the
configurations in $\L$
having at most $\ln\ln\beta$ minus spins,
\[
{\cE}(\L) = \biggl\{ \eta\in \{ -1,+1 \}^\L\dvtx  \sum
_{x\in\L}\eta(x)\geq|\L|-2\ln\ln\beta \biggr\}.
\]
We set finally
\[
T'_x = \inf \bigl\{ t\geq T_x\dvtx
\s_{
\L(\exp(\beta L)),t}^{-,\minus}|_{\L_x}\notin {\cE}(\L_x) \bigr\}.
\]
The infection process $(\mu_t)_{t\geq0}$ is given by
\[
\forall x\in \L\bigl(\exp(\beta L)\bigr)\qquad \mu_t(x) = %
\cases{0, &\quad if $t< T_x$,
\vspace*{1pt}\cr
1, &\quad if $T_x\leq t< T'_x$,
\vspace*{1pt}\cr
0, &\quad if $t\geq T'_x$.}
\]
We first show that, once a site is infected, with very high
probability, it remains
infected until time $\tb$.

%le7.1 #&#
\begin{lemma}\label{eaz}
%Let $C>0$.
For any $x$ in $
\L(\exp(\beta L))
$,
\[
\forall C>0\qquad P \bigl(T'_x-T_x\leq\exp(\beta
C) \bigr) = \SES.
\]
\end{lemma}
\begin{pf}
From the Markov property and the monotonicity with respect to the
boundary conditions, we have
\begin{eqnarray*}
&&P \bigl(T'_x-T_x\leq\exp(\beta C) \bigr)
\\
% =
%P\Big(
%$(\s_{\L_\beta,t}^{-,\minus}|_{\L(x)})_{t\geq0}$}
% \tau(\cE(\L_x))\leq
&&\qquad\leq P \bigl( \mbox{for the process $\bigl(\s_{\L_x,t}^{-,\plus}
\bigr)_{t\geq0}$}, \tau\bigl(\cE(\L_x)\bigr)\leq \exp(\beta C)
\bigr).
\end{eqnarray*}
We consider the dynamics in $\L_x$ starting from $\plus$ and
restricted to the set
$\cE(\L_x)$, with $-$ boundary conditions on $\L_x$. We denote
by
$(\widehat{\s}_{\L_x,t}^{-,\plus})_{t\geq0}$
the corresponding process.
The invariant measure of this process is the Gibbs measure restricted to
$\cE(\L_x)$, which we denote by
$\widehat{\mu}_{\L_x}$,
\[
\forall\s\in\cE(\L_x)\qquad \widehat{\mu}_{\L_x} (\s) =
\frac{
\mu_{\L_x}^{-}(\s)} {
\mu_{\L_x}^{-}(\cE(\L_x))}.
\]
We use the graphical construction described in Section~\ref{sigm}
to couple the processes
\[
\bigl({\s}_{\L_x,t}^{-,\plus}\bigr)_{t\geq0},\qquad \bigl(
\widehat{\s}_{\L_x,t}^{-,\widehat{\mu}}\bigr)_{t\geq0}.
\]
We define
\[
\partial^{\inn}\cE(\L_x) = \bigl\{ \sigma\in \cE(
\L_x)\dvtx \exists y\in\L_x, \s^y\notin \cE(
\L_x) \bigr\}.
\]
Proceeding as in Lemma \ref{fugaup}, we obtain that
\begin{eqnarray*}
&&
P \bigl( \mbox{for the process $\bigl(\s_{\L_x,t}^{-,\plus}
\bigr)_{t\geq0}$}, \tau\bigl(\cE(\L_x)\bigr)\leq \exp(\beta C)
\bigr)
\\
&&\qquad\leq P \bigl(\exists t\leq\exp(\beta C), \widehat{\s}_{\L_x,t}^{-,\widehat{\mu}}
\in \partial^{\inn}\cE(\L_x) \bigr)
\\
&&\qquad\leq 4\b\lambda \mh_{\L_x}^{-} \bigl( \partial^{\inn}
\cE \bigr) +\exp(-\b\lambda\ln\b),
\end{eqnarray*}
where $\lambda=(\ln\beta)^d
\exp(\beta C)$.
Next, if $\eta\in
\partial^{\inn}\cE(\L_x)$, then
\[
\sum_{y\in\L_x}\eta(y)\leq|\L_x|-2\ln\ln
\beta+1
\]
and
\[
H_{\L_x}^{-}(\eta) - H_{\L_x}^{-}(\plus)
\geq h(\ln\ln\beta-1)
\]
so that
\[
\mh_{\L_x}^{-}(\eta) \leq \exp \bigl(-\beta h(\ln\ln\beta-1)
\bigr).
\]
Thus
\begin{eqnarray*}
\mh_{\L_x}^{-} \bigl( \partial^{\inn} \cE \bigr)
&\leq& \bigl| \partial^{\inn} \cE \bigr|\min \bigl\{ \mh_{\L_x}^{-}(
\eta)\dvtx  \eta\in \partial^{\inn}\cE(\L_x) \bigr\}
\\
&\leq& \bigl((\ln\beta)^d \bigr)^{\ln\ln\beta} \exp \bigl(-\beta( h
\ln\ln\beta-1) \bigr).
\end{eqnarray*}
This last quantity is $\SES$ and the lemma is proved.
\end{pf}
%
%s7.2 #&#
\subsection{Spreading of the infection}
\label{ispread}
We show first that any configuration in
$\cE(\L_x)$ can reach the configuration $\plus$ through a downhill path.

%le7.2 #&#
\begin{lemma}\label{down}
Let $\eta$ belong to
$\cE(\L_x)$. There exists a sequence of $r\leq\ln\ln\beta$
distinct sites $x_1,\ldots,x_r$ such that, if
we set
$\s_0=\eta$ and
\[
\forall i\in\{ 1,\ldots,r \}\qquad \s_i = \s_{i-1}^{x_i},
\]
then we have $\s_r=\plus$ and
for $i\in\{ 1,\ldots,r \}$,
$\eta(x_i)=\s_{i-1}(x_i)= -1$
and
$x_i$ has at least $d$ plus neighbors in $\s_{i-1}$.
\end{lemma}
\begin{pf}
We prove the result by induction over the dimension $d$.
Suppose first that $d=1$.
Let $\eta$ be a configuration in
$\cE(\L^1(\ln\beta))$.
Let $x_0\in\L^1(\ln\beta)$ such that $\eta(x_0)=1$.
We define then
\begin{eqnarray*}
x_1 &=& \max \bigl\{ y< x_0\dvtx \eta(y)=-1 \bigr\},
\\
&&\vdots
\\
x_k &=& \max \bigl\{ y< x_{k-1}\dvtx \eta(y)=-1 \bigr\},
\\
x'_1 &=& \min \bigl\{ y> x_0\dvtx \eta(y)=-1
\bigr\},
\\
&&\vdots
\\
x'_l &=& \min \bigl\{ y> x'_{l-1}\dvtx
\eta(y)=-1 \bigr\}.
\end{eqnarray*}
The sequence of sites
$x_1,\ldots,x_k,x'_1,\ldots,x'_l$ answers the problem.
Suppose that the result has been proved at rank $(d-1)$.
Let $\eta$ be a configuration in
$\cE(\L^d(\ln\beta))$.
We consider the hyperplanes
\[
P_i = \bigl\{ x=(x_1,\ldots,x_d)\in
\Z^d\dvtx  x_{d}=i \bigr\},\qquad i\in\Z,
\]
and we denote by $\eta_i$ the restriction of $\eta$ to $P_i$.
The configuration $\eta_i$ can naturally be identified with a
$(d-1)$-dimensional configuration. Since there is at most
$\ln\ln\beta$ minuses in the configuration $\eta$, there exists an index
$i^*$ such that
$\eta_{i^*}=\plus$. We apply next the induction result at rank $d-1$ to
$\eta_{i^*+1}$. This way, we can fill
$P_i\cap\L^d(\ln\beta)$
with a sequence of positive spin flips which never increase the $(d-1)$-dimensional energy. Each site
which is flipped in
$\eta_{i^*+1}$ has at least $d-1$ plus neighbors in $P_{i^*+1}$, hence
at least
$d$ plus neighbors in
$\L^d(\ln\beta)$. Thus no spin flip of this sequence
increases the $d$-dimensional energy. We iterate the argument, filling successively the
sets
$P_i\cap\L^d(\ln\beta)$ above and below $i^*$ until the box
$\L^d(\ln\beta)$ is completely filled.
\end{pf}

This result leads directly to a lower bound on the time needed to reach the
configuration $\plus$ starting from a configuration of
$\cE(\L^d(\ln\beta))$.
%
%co7.3 #&#
\begin{corollary}\label{zddz}
For any configuration $\eta$ in
$\cE(\L_x)$, we have
%$\cE(\L^d(\ln\beta))$, we have
%
\[
P \bigl( %\mbox{for the process
%$(\s_{\L_x,t}^{-,\eta})_{t\geq0}$}
% \tau(\cE(\L_x))\leq
\s_{\L_x,
\ln\ln\beta
}^{-,\eta}=\plus \bigr) \geq 7^{-|\L_x|\ln\ln\beta}.
\]
\end{corollary}
\begin{pf}
Let $\eta\in
\cE(\L_x)$, and
let $x_1,\ldots,x_r$,
$r\leq\ln\ln\beta$,
be a sequence of sites as given by
Lemma \ref{down}. We evaluate the probability that,
starting from $\eta$, the successive spin flips at
$x_1,\ldots,x_r$ occur.
For $i\in\{ 1,\ldots,r \}$, let $E_i$ be the event:
during the time interval $[i-1,i]$, there is a time arrival
for the Poisson process associated to the site $x_i$, and none
for the other sites of the box $\L_x$.
Let $F$ be the event that there is no arrival for
the Poisson processes in the box $\L_x$ during $[r,\ln\ln\beta]$.
We have then
\begin{eqnarray*}
P(F) &\geq& \biggl(1- \frac{1}{e} \biggr)^{|\L_x|\ln\ln\beta},
\\
\forall i \in\{ 1,\ldots,r \}\qquad P(E_i) &\geq& \frac{1}{e}
\biggl(1- \frac{1}{e} \biggr)^{|\L_x|} % \geq
\end{eqnarray*}
and
\[
P \biggl(F\cap\bigcap_{1\leq i\leq r}E_i \biggr)
= P(F)\times \prod_{i\in I}P(E_i) \geq
% \geq
7^{-|\L_x|\ln\ln\beta}.
\]
Yet the event $E_1\cap\cdots\cap E_r\cap F$ implies that, at time
$r$, the process
starting from $\eta$ has reached the configuration $\plus$ and that
it does not move until
time $\ln\ln\beta$.
\end{pf}

For $x\in
\L(\exp(\beta L))
$, we define the enlarged neighborhood $\L'_x$ of $\L_x$
as
\[
\L'_x = \bigcup_{y:|y-x|=1}
\L_y.
\]

%pr7.4 #&#
\begin{proposition}\label{nucl}
Let $n\in\{ 1,\ldots,d \}$.
Let $\eta$ be a configuration in
$\L'_x$ such that there exist $d-n$ neighbors $y_1,\ldots,y_{d-n}$
of $x$ in $d-n$ distinct directions
for which the restriction
$\eta|_{\L_{y_i}}$
is in $\cE(\L_{y_i})$ for
$i\in\{ 1,\ldots,d-n \}$.
We have the following estimates:\vspace*{9pt}

\textit{Nucleation}:
For any $\kappa$ such that
$\G_{n-1}<\kappa<\G_{n}$ and $\varepsilon>0$,
we have for $\beta$ large enough
\[
P\pmatrix{ \mbox{in the process }\bigl(
\s_{\L'_x,t}^{-,\eta}\bigr)_{t\geq0}\mbox{, the site $x$}
\vspace*{1pt}\cr
\mbox{becomes infected before time $\exp(\beta\kappa)$}} \geq \exp \bigl(\beta(
\kappa-\G_{n}-\varepsilon) \bigr).
\]
%
%P\left(
%site $x$}\\

\vspace*{9pt}

\textit{Spreading}: For any $\kappa$ such that
$\kappa>\G_{n}$, we have
\[
P\pmatrix{ \mbox{in the process }\bigl(\s_{\L'_x,t}^{-,\eta}
\bigr)_{t\geq0}\mbox{, the site $x$ has}
\vspace*{1pt}\cr
\mbox{not become infected by
time $\exp(\beta\kappa)$}} =\SES.
\]
\end{proposition}

%We consider next a finite box $Q$ satisfying the following hypothesis.
%
%{\bf Hypothesis on $Q$.} The box $Q$ is such that
%$|Q| = \exp o(\ln\beta)$, which means that
%$$\lim_{\beta\to\infty}
%Let $\kappa$ be such that
%$\G_{d-1}<\kappa<\G_d$.
%We have
%$$\lim_{\beta\to\infty}
% =
%For any $\e>0$, we have for $\b$ sufficiently large
% \le
%where $Q = \L^d(\log\b)$.
%
\begin{pf}
%The upper bound follows from reversibility
%and is a very well-known and very general result
%(see e.g. [MO2] Lemma 4.1).
% The upper bound is a consequence of Lemma \ref{fugaup}.
%We prove now the lower bound.
We consider the process
$(\s_{\L'_x,t}^{-,\eta})_{t\geq0}$,
and we set
\[
\tau_\plus = \t\bigl({\{ \plus|_{\L'_x} \}^c}
\bigr) = \inf \bigl\{ t\geq0\dvtx \forall y\in\L'_x,
\s_{\L'_x,t}^{-,\eta}(y)=+1 \bigr\}
\]
the hitting time of the configuration equal to $+1$ everywhere in $\L'_x$.
Let
\[
I = \exp(\b\kappa)- \exp(\b\G_{n-1})
\]
and let $\theta$ be the time of the last visit to $\minus|_{\L_x}$
before reaching $\plus|_{\L'_x}$,
\[
\theta = \sup \bigl\{ t\leq \tau_{\plus}\dvtx \forall y\in\L_x,
\s_{\L'_x,t}^{-,\eta}(y)=-1 \bigr\}.
\]
In case the process does not visit
$\minus|_{\L_x}$ before
$\tau_{\plus}$,
we set $\theta=0$.
Let $\alpha$ be the configuration in $\L'_x$ such that
\[
\forall y\in\L'_x\qquad \alpha(y) = %
\cases{+1, &\quad if $\displaystyle y\in\bigcup_{1\leq i\leq d-n}
\L_{y_i}$,
\vspace*{1pt}\cr
-1, &\quad if $y\in\L_x$.}
\]
We write, using the Markov property,
\begin{eqnarray*}
&&\P \bigl(\t_\plus< \exp(\b\kappa) \bigr)
\\
&&\qquad\geq
\sum_{0\leq i \leq I} \P \bigl(\sigma_{\L'_x,i}^{-,\eta}=
\alpha, i\leq\theta<i+1, \t_\plus< i+ \exp(\b\G_{n-1}) \bigr)
\\
&&\qquad\geq \biggl(\sum_{0\leq i \leq I} \P \bigl(
\sigma_{\L'_x,i}^{-,\eta}=\alpha, \t_\plus> i \bigr) \biggr)
\P \pmatrix{ \mbox{for the process } \bigl(\s_{\L'_x,t}^{-,\alpha}
\bigr)_{t\geq0}
\vspace*{1pt}\cr
0\leq\theta<1, \t_\plus< \exp(\b
\G_{n-1})}.
\end{eqnarray*}
By Proposition \ref{dee}, the maximal depth in the reference
cycle path in the box $\L_x$ with $n\pm$ boundary conditions
is strictly less than $\G_{n-1}$, so that we have
for $\e>0$ and $\b$ large enough
\[
\P \pmatrix{ \mbox{for the process } \bigl(\s_{\L'_x,t}^{-,\alpha}
\bigr)_{t\geq0}
\vspace*{1pt}\cr
0\leq\theta<1, \t_\plus< \exp(\b
\G_{n-1})} %\P\big(0\leq\theta<1,
\geq \exp \bigl(-\b(
\G_n+\e) \bigr).
\]
This estimate is a continuous-time analog of Theorem $5.2$ and
Proposition $10.9$ of \cite{CaCe}.
It relies on a continuous-time formula
giving the expected exit time given the exit point, which
is the analog of Lemma 10.2 of \cite{CaCe}.
Let $\cC_n^\alpha$ be the largest cycle included in
$\{ -1,+1 \}^{\L'_x}$
containing $\alpha$ and not $\plus$.
For $i\leq I$, we have
\begin{eqnarray*}
\P \bigl(\sigma_{\L'_x,i}^{-,\eta}=\alpha, \t_\plus> i
\bigr) &\geq& \P \bigl(\sigma_{\L'_x,i}^{-,\eta}=\alpha, \t\bigl(
\cC_n^\alpha\bigr)> i \bigr)
\\
&\geq& \P \bigl( \mbox{for the process } \bigl(\s_{\L'_x,t}^{-,\alpha}
\bigr)_{t\geq0}, \t\bigl(\cC_n^\alpha\bigr)> I \bigr)\\
&&{}\times
\P \bigl(\sigma_{\L'_x,i}^{-,\alpha}=\alpha | \t\bigl(
\cC_n^\alpha\bigr)>i \bigr).
\end{eqnarray*}
Since $\kappa<\G_n$, then
\[
\lim_{\b\to\infty}\P \bigl( \mbox{for the process } \bigl(
\s_{\L'_x,t}^{-,\alpha}\bigr)_{t\geq0}, \t\bigl(
\cC_n^\alpha\bigr)> I \bigr) = 1.
\]
This follows from the continuous-time analog of corollary $10.8$
of \cite{CaCe}.
We compare next the process starting from $\alpha$
with the process starting from
$\mt_{\cC_n^\alpha}$,
the Gibbs measure
restricted
to the metastable cycle $\cC_n^\alpha$. We have
\begin{eqnarray*}
\mt_{\cC_n^\alpha}(\alpha) &=& \sum_{\eta\in
\cC_n^\alpha}
\mt_{\cC_n^\alpha}(\eta) \P \bigl(\sigma_{\L'_x,i}^{-,\eta}=\alpha |
\t\bigl(\cC_n^\alpha\bigr)>i \bigr)
\\
&\leq& \mt_{\cC_n^\alpha}(\alpha) \P \bigl(\sigma_{\L'_x,i}^{-,\alpha}=
\alpha | \t\bigl(\cC_n^\alpha\bigr)>i \bigr) + \sum
_{\eta\in
\cC_n^\alpha,\eta\neq\alpha} \mt_{\cC_n^\alpha}(\eta).
\end{eqnarray*}
The configuration $\alpha$ is the bottom of the cycle $\cC_n^\alpha$.
Thus there exists $\delta>0$ such that
\[
\forall\eta\in\cC_n^\alpha\qquad \eta\neq\alpha
\quad\Longrightarrow\quad
\mt_{\cC_n^\alpha}(\eta) \leq \mt_{\cC_n^\alpha}(\alpha) \exp(-\beta\delta).
\]
For $\beta$ large enough, we have also
$|\cC_n^\alpha | \leq \exp(\beta\delta/2)$.
We conclude that
\[
\P \bigl(\sigma_{\L'_x,i}^{-,\alpha}=\alpha | \t\bigl(
\cC_n^\alpha\bigr)>i \bigr) \geq %\frac{\mt_{\cC_n^\alpha}(\alpha)}
%{\mt_{\cC_n^\alpha}(\alpha)-
%{\sum_{\eta\in
%%\cC_n^\alpha\cap\bigcap_{1\leq i\leq d-n} \cE(\L_{y_i})
%}\exp\big(-\b H_{\L'_x}^-(\eta)\big)}\\
% \geq
\frac{1} {
1
-\exp(-\beta\delta/2)}. %-\big|\cC_n^\alpha\big|\exp(-\beta\delta)}.
\]
%
%The last inequality is a consequence of Lemma \ref{eaz}: up to a $
%$t\leq\tb$, for any
%$i\in\{ 1,\ldots,d-n \}$, the configuration
%$\s_{\L'_x,t}|_{\L_{y_i}}$ belongs to $\cE(L_{y_i})$.
Combining these estimates, we conclude that for $\b$ large enough,
\[
\P \bigl(\t_\plus< \exp(\b\kappa) \bigr) \geq I %\frac{1}{2} \big|Q\big|^{-m_n}
\exp
\bigl(-\b(\G_n+\e) \bigr).
\]
Sending successively $\b$ to $\infty$ and $\e$ to $0$, we
obtain the desired lower bound.
The second estimate stated in the proposition is a standard consequence
of the first.
\end{pf}
%
%s7.3 #&#
\subsection{Invasion}
\label{iinva}
We denote by $e_1,\ldots,e_d$ the canonical orthonormal basis of
$\mathbb R^d$.
We will prove the following result by induction over $n$.
%
%pr7.5 #&#
\begin{proposition}\label{inva}
Let $n\in\{ 0,\ldots,d \}$, and let $L\geq0$.
Let $\L_\beta^n$ be the parallelepiped
\[
\L_\beta^n = \L^{n}\bigl(\exp(\beta L)\bigr)
\times\L^{d-n}(1).
\]
For any $s\geq0$
and any
$\kappa>
\max ( \Gamma_{n}-nL, \kappa_{n} )$, we have
\[
P \left(\matrix{ \mbox{all the sites of $\L_\beta^n$ are}
\vspace*{1pt}\cr
\mbox{infected at time} %$s +\exp\Big(\beta\max\big( \Gamma_{n}-nL, \kappa_{n}\big)\Big) $
\vspace*{1pt}\cr
\mbox{$s +\exp(\beta
\kappa)$}} \right|\left. \matrix{ \mbox{all the sites of}
\vspace*{1pt}\cr
\mbox{$e_{n+1}+
\L_\beta^n,\ldots, e_{d}+\L_\beta^n$}
\vspace*{1pt}\cr
\mbox{are infected at time $s$}} \right) = 1-\SES.
\]
\end{proposition}
%
%1The process in $\Z^d$, with void boundary conditions
%1is denoted by $\s^\a_t$.
%1If $\a=\zero$, we drop the initial configuration from the notation.
%
\begin{pf}
Thanks to the Markovian character of the process, we need only to consider
the case where $s=0$.
Let us consider first the case $n=0$.
We have then $\kappa_0=\Gamma_0=0$.
The box
$\L_\beta^0$ is reduced to the singleton $\{ 0 \}$.
The result is an immediate consequence of Proposition \ref{nucl}.
%Let $\kappa>0$ and let $\tau_\beta=\exp(\beta\k)$.
We suppose now that $n\geq1$ and that the result has been proved
at rank $n-1$.
Let $L>0$, let $\L_\beta^n$ be a parallelepiped as in the
statement of the proposition and let
$\kappa> \max( \Gamma_{n}-nL, \kappa_{n})$.
We define the nucleation time
$\tau_{\mathrm{nucleation}}$
in $\Lambda_\beta^n$ as
\[
\tau_{\mathrm{nucleation}} = \inf \bigl\{ t\geq0\dvtx \exists
x\in\Lambda^n_\beta,
\mu_{t}(x)=1 \bigr\}.
\]
Let
$c> \max( \Gamma_{n}-nL, \Gamma_{n-1})$.
Let $(x_i)_{i\in I}$ be a family of sites of
$\L_\beta^n$ which are pairwise at distance larger than
$4\ln\beta$ and such that
\[
|I| \geq \frac{\exp(\beta L n)}{(6\ln\beta)^{n}}.
\]
We can, for instance, consider the sites
of the sublattice
$(5\ln\beta){\mathbb Z^n}\times\L^{d-n}(1)$
which are included in
$\L_\beta^n$.
For $i\in I$,
let $\eta_i$ be the initial configuration restricted to the box
${\L'_{x_i}}$. We write
\begin{eqnarray*}
\P \bigl( \tau_{\mathrm{nucleation}} > \exp(\beta c) \bigr) &\leq& P\pmatrix{ \mbox{no
site $x$ in $\L_\beta^n$ has become}
\vspace*{1pt}\cr
\mbox{infected by
time $\exp(\beta c)$} %\mbox{in the process }(\s_{\L'_x,t}^{-,\eta})_{t\geq0}
}
\\
&\leq& P\pmatrix{ \mbox{for any $i$ in $I$, the site $x_i$ has
not}
\vspace*{1pt}\cr
\mbox{become infected by time $\exp(\beta c)$}
\vspace*{1pt}\cr
\mbox{in the process }
\bigl(\s_{\L'_{x_i},t}^{-,\eta_i}\bigr)_{t\geq0}}
\\
&\leq& \prod_{i\in I} P\pmatrix{ \mbox{the site
$x_i$ has not become infected by}
\vspace*{1pt}\cr
\mbox{time $\exp(\beta c)$
in the process }\bigl(\s_{\L'_{x_i},t}^{-,\eta_i}
\bigr)_{t\geq0}}.
\end{eqnarray*}
Since all the sites of
$e_{n+1}+\L_\beta^n,\ldots, e_{d}+\L_\beta^n$
are initially infected,
by Proposition~\ref{nucl}
we have for any $\varepsilon>0$,
\[
\P \bigl( \tau_{\mathrm{nucleation}} > \exp(\beta c) \bigr) \leq \bigl(1-\exp \bigl(
\beta(c-\G_{n}-\varepsilon) \bigr) \bigr)^{
{\exp(\beta L n)}/({(6\ln\beta)^{n}})
}.
\]
Therefore,
up to a SES event, the first infected site in the box
$\L^n_\beta$
appears before time
$\exp(\beta c)$.
%$$\exp\big(\beta
%+\beta\varepsilon\big)
%.
%$$
For $i\geq1$, we define
the first time $\tau^i$ when there is a $n$-dimensional
parallelepiped
of infected sites
of diameter
larger than or equal to $i$ in $\Lambda_\beta^n$, that is,
\begin{eqnarray*}
\tau^i &=& \inf \bigl\{ t\geq0\mbox{: there is a $n$-dimensional parallelepiped}
\\
&&\hspace*{17pt}\mbox{of infected sites included in $
\Lambda_\beta^n$ whose}\\
&&\hspace*{44pt}
\mbox{$\dinf$ diameter is larger
than or equal to $i$} \bigr\}.
\end{eqnarray*}
%
%$\dinf C\geq i$}
%
% \big\}.
The face of an $n$-dimensional parallelepiped is an $n-1$-dimensional
parallelepiped.
The sites of a face of an infected parallelepiped
in $\L_\beta^n$ have already $d-n+1$ infected neighbors.
From the induction hypothesis, up to a SES event, an
$n-1$-dimensional box
of side length $\exp(\beta K)$
whose sites have already
$d-n+1$ infected neighbors
is fully infected at
a time
\[
\exp \bigl(\beta \bigl( \max\bigl(\Gamma_{n-1}-(n-1)K,
\kappa_{n-1}\bigr)+\varepsilon \bigr) \bigr).
\]
This implies that, up to a SES event,
the box $\L_\beta^n$ is fully occupied
at time
\begin{eqnarray*}
&&
\tau^ {\exp(\beta L)} \\
&&\qquad\leq \tau_{\mathrm{nucleation}} + \sum
_{1\leq i<
{\exp(\beta L)}
} \bigl(\tau^{i+1} - \tau^{i}\bigr)
\\
&&\qquad\leq \exp(\beta c) + \sum_{1\leq i< \exp(\beta L)
}\!\! 2n\exp \biggl(\beta
\biggl( \max\biggl(\Gamma_{n-1}-\frac{n-1}{\beta}\ln i,
\kappa_{n-1}\biggr)+\varepsilon \biggr)\!\biggr). %\\
% \leq
%+
%}
\end{eqnarray*}
We consider two cases.

$\bullet$ First case: $L\leq L_{n-1}$.
Notice that $L_0=0$, hence this case can happen only whenever $n\geq2$.
In this case, we have
\[
\forall i< \exp(\beta L)\qquad \kappa_{n-1} \leq \Gamma_{n-1}-
\frac{n-1}{\beta}\ln i
\]
and
\begin{eqnarray*}
&&
\sum_{1\leq i< \exp(\beta L)} \exp \biggl(\beta \max\biggl(
\Gamma_{n-1}-\frac{n-1}{\beta}\ln i,\kappa_{n-1}\biggr) \biggr)
\\
&&\qquad
\leq \exp(\beta \Gamma_{n-1}) \sum_{1\leq i< \exp(\beta L)}
\frac{1}{i^{n-1}}
\\
&&\qquad
\leq \exp(\beta \Gamma_{n-1}) \sum_{1\leq i< \exp(\beta L)}
\frac{1}{i} \leq \beta L \exp(\beta \Gamma_{n-1}).
\end{eqnarray*}

$\bullet$ Second case: $L> L_{n-1}$. We have then
\begin{eqnarray*}
&&\sum_{
\exp(\beta L_{n-1})
\leq i< \exp(\beta L)} \exp \biggl(\beta \max\biggl(
\Gamma_{n-1}-\frac{n-1}{\beta}\ln i,\kappa_{n-1}\biggr) \biggr)
\\
&&\qquad\leq \bigl( \exp(\beta L)- \exp(\beta L_{n-1}) \bigr) \exp(\beta
\kappa_{n-1})
\\
&&\qquad\leq \exp \bigl(\beta(L+\kappa_{n-1}) \bigr).
\end{eqnarray*}
We conclude that, in both cases, for any
$\varepsilon>0$, up to a SES event,
the box $\L_\beta^n$ is fully occupied
at time
\[
2n\beta L\exp(\beta\varepsilon) \bigl( \exp \bigl(\beta( \Gamma_{n}-nL)
\bigr) + \exp(\beta \Gamma_{n-1}) + \exp \bigl(\beta(L+
\kappa_{n-1}) \bigr) \bigr).
\]
Therefore, for any $\kappa$ such that
\[
\kappa > \max ( \Gamma_{n}-nL, \Gamma_{n-1}, L+
\kappa_{n-1} )
\]
the probability that
the box $\L_\beta^n$ is not fully occupied
at time
$\exp(\beta\kappa)$
is SES.
If $L\leq L_n$, then
\[
\max ( \Gamma_{n}-nL, \Gamma_{n-1}, L+\kappa_{n-1}
) = \Gamma_{n}-nL
\]
and we have the desired estimate.
Suppose next that
$L> L_n$. By the previous result applied with $L=L_n$,
we know that,
%, we use first the result for boxes
%of sidelength $L_d$ to get that,
for any $\kappa>\kappa_n$,
up to a SES event,
%the probability that
a box of side length
$\exp(\beta L_n)$
%the box $\L_\beta$
is fully occupied
at time
$\exp(\beta\kappa)$.
We cover $\L_\beta^n$ by boxes
of sidelength
$\exp(\beta L_n)$. Such a cover contains at most
$\exp(\beta nL)$ boxes, thus
\begin{eqnarray*}
&&
\P \bigl( \mbox{$\L_\beta^n$ is not fully occupied at time
$\tau_\beta$} \bigr)
\\
&&\qquad
\leq \P\pmatrix{ \mbox{there exists a box included in $\L_\beta^n$
of side length}
\vspace*{1pt}\cr
\mbox{$\exp(\beta L_n)$ which is not fully
occupied at time $\tau_\beta$}}
\\
&&\qquad\leq \exp(\beta nL) \P\lleft(
\matrix{ \mbox{the box $\L^n\bigl( \exp(\beta L_n)\bigr)$
is not}
\vspace*{1pt}\cr
\mbox{fully occupied at time $\tau_\beta$}} \rright).
\end{eqnarray*}
The last probability being SES, we are done.
\end{pf}

Proposition \ref{inva} with $n=d$ readily yields the upper bound of
the relaxation time stated in Theorem \ref{mainfv}.

\makeatletter\write@toc@ignorecontentsline\makeatother
\section*{Acknowledgments}

Rapha\"el Cerf warmly thanks Roberto Schonmann for discussions on this
problem while he visited UCLA in 1995. We thank two anonymous Referees
for their careful reading of the manuscript and their numerous
comments.
\makeatletter\write@toc@restorecontentsline\makeatother

%suskaldyti doi

% imsref loaded by lrinkeviciute, 2013-07-09 13:40:40

\printaddresses

\end{document}